\documentclass{birkmult}

\usepackage{amsmath}
\usepackage{amssymb}
 \newtheorem{thm}{Theorem}[section]
 \newtheorem{cor}[thm]{Corollary}
 \newtheorem{lem}[thm]{Lemma}
 
 \theoremstyle{definition}
 
 \theoremstyle{remark}
 \newtheorem{rem}[thm]{Remark}
 \newtheorem{ex}{Example}
 \numberwithin{equation}{section}

\begin{document}
%--------------------------------------------------------------------------
% editorial commands: to be inserted by the editorial office
%
%\firstpage{1}
%\issuenumber{4}
%\Volumeandyear{41 (2001)}
%\Copyrightyear{2002}
%\submitted{March 14, 2000}
\title[Matrix-$J$-unitary Rational Formal Power Series]
 {Matrix-$J$-unitary Non-commutative Rational\\ Formal Power Series}

%----------Author 1
\author{D. Alpay}

\address{%
Department of Mathematics\\
Ben-Gurion University of the Negev\\
Beer-Sheva 84105, Israel}

\email{dany@math.bgu.ac.il}

%----------Author 2
\author{D.~S.~Kalyuzhny\u{\i}-Verbovetzki\u{\i}}
\address{%
Department of Mathematics\\
Ben-Gurion University of the Negev\\
Beer-Sheva 84105, Israel}
\email{dmitryk@math.bgu.ac.il}

\thanks{The second author was supported by the Center for Advanced Studies in Mathematics, Ben-Gurion University of the Negev.}


\date{}

\begin{abstract}
Formal power series in $N$ non-commuting indeterminates can be considered as a counterpart of functions of one variable holomorphic at $0$, and some of their properties are described in terms of coefficients.
However, really fruitful analysis  begins when one considers for them evaluations on $N$-tuples of $n\times n$ matrices (with $n=1,2,\ldots$) or operators on an infinite-dimensional separable Hilbert space. Moreover, such evaluations appear in control, optimization and stabilization problems of modern system engineering.

In this paper, a theory of realization and minimal factorization of rational matrix-valued functions which are $J$-unitary on the imaginary line or on the unit circle is extended to the setting of non-commutative rational formal power series. The property of $J$-unitarity holds on $N$-tuples of $n\times n$ skew-Hermitian versus unitary matrices ($n=1,2,\ldots$), and a rational formal power series is called \emph{matrix-$J$-unitary} in this case. The close relationship between minimal realizations and structured Hermitian solutions $H$ of the Lyapunov or Stein equations is established. The results are specialized for the case of \emph{matrix-$J$-inner} rational formal power series. In this case $H>0$, however the proof of that is more elaborated than in the one-variable case and involves a new technique. For the rational \emph{matrix-inner} case, i.e., when $J=I$, the theorem of Ball, Groenewald and Malakorn on unitary realization of a formal power series from the non-commutative Schur--Agler class admits an improvement: the existence of a minimal (thus, finite-dimensional) such unitary realization and its uniqueness up to a unitary similarity is proved. A version of the theory for \emph{matrix-selfadjoint} rational formal power series is also presented. The concept of non-commutative formal reproducing kernel Pontryagin spaces is introduced, and in this framework the backward shift realization of a matrix-$J$-unitary rational formal power series in a finite-dimensional non-commutative de Branges--Rovnyak space is described.
\end{abstract}

%%% ----------------------------------------------------------------------
\maketitle
%%% ----------------------------------------------------------------------
\tableofcontents
\section{Introduction}\label{sec:intr}
In the present paper we study a non-commutative analogue of rational matrix-valued functions which are $J$-unitary on the imaginary line or on the unit circle and, as a special case, $J$-inner ones.
 Let $J\in{\mathbb C}^{q\times q}$ be a signature
matrix, i.e., a matrix which is both self-adjoint and unitary.
A ${\mathbb C}^{q\times q}$-valued rational function $F$ is
\emph{$J$-unitary on the imaginary line} if
\begin{equation}
\label{mairie de Montreuil} {F}(z)J{F}(z)^*=J
\end{equation}
at every point of holomorphy of ${F}$ on the imaginary line. It is
called \emph{$J$-inner} if moreover
\begin{equation}
\label{croix de Chavaux}
{F}(z)J{F}(z)^*\le J
\end{equation}
at every point of holomorphy of ${F}$ in the open right
half-plane $\Pi$. Replacing the imaginary line by the unit circle
${\mathbb T}$ in \eqref{mairie de Montreuil} and the open right
half-plane $\Pi$ by the open unit disk ${\mathbb D}$ in \eqref{croix
de Chavaux}, one defines $J$-unitary functions  on the unit circle
(resp., $J$-inner functions in the open unit disk). These classes of rational functions were studied in \cite{AG1} and \cite{AG2}
using the theory of realizations of rational matrix-valued
functions, and in \cite{AD1} using the theory of reproducing kernel
Pontryagin spaces. The circle and line cases were studied in a
unified way in \cite{AD2}. We mention also the earlier papers \cite{KP,EP} that inspired much of ivestigation of these and other classes of rational matrix-valued functions with symmetries.

We now recall some of the arguments in
\cite{AG1}, then explain the difficulties appearing in the
several complex variables setting, and why the arguments of
\cite{AG1} extend to the non-commutative framework. So let $F$
be a rational function which is $J$-unitary on the imaginary line, and assume
that $F$ is holomorphic in a neighbourhood of the origin. It then
admits a minimal realization
$$F(z)=D+C(I_\gamma-zA)^{-1}zB$$
where $D={F}(0)$, and $A,B,C$ are matrices of appropriate sizes (the size $\gamma\times\gamma$ of the square matrix $A$ is minimal possible for such a realization). Rewrite \eqref{mairie de Montreuil} as
\begin{equation}
\label{jean yanne} {F}(z)=J{F}(-\overline{z})^{-*}J,
\end{equation}
where $z$ is in the domain of holomorphy of both ${F}(z)$ and
${F}(-\overline{z})^{-*}$. We can rewrite \eqref{jean yanne} as
$$D+C(I_\gamma-zA)^{-1}zB=J\left(D^{-*}+D^{-*}B^*(I_\gamma +z(A-BD^{-1}C)^*)^{-1}zC^{*}D^{-*}\right)J.$$
The above equality gives two minimal realizations of a given rational
matrix-valued function. These realizations are therefore similar, and there
is a uniquely defined matrix (which, for convenience, we denote by $-H$) such that
\begin{equation}
\label{washington}
\begin{pmatrix}
-H&0\\
0&I_q\end{pmatrix}
\begin{pmatrix}A&B\\C&D
\end{pmatrix}=
\begin{pmatrix}
-(A^*-C^*D^{-*}B^*)&C^*D^{-*}J\\
JD^{-*}B^*&JD^{-*}J\end{pmatrix}
\begin{pmatrix}
-H&0\\
0&I_q\end{pmatrix}.
\end{equation}
The matrix $-H^*$ in the place of $-H$ also satisfies \eqref{washington}, and by
uniqueness of the similarity matrix we have $H=H^*$, which leads to the following theorem.
\begin{thm}
\label{th1}
Let $F$ be a rational matrix-valued function holomorphic in a
neighbourhood of the origin and let $F(z)=D+C(I_\gamma-zA)^{-1}zB$ be a
minimal realization of $F$. Then $F$ is $J$-unitary on the imaginary
line if and only if the following conditions hold:

$(1)$ $D$ is $J$-unitary, that is, $DJD^*=J$;

$(2)$ there exists an Hermitian invertible matrix $H$ such that
\begin{eqnarray}
\label{st-ambroise}
A^*H+HA&=&-C^*JC,\\
B&=&-H^{-1}C^*JD.
\end{eqnarray}
The matrix $H$ is uniquely determined by a given minimal
realization (it is called the associated Hermitian matrix to this realization). It holds that
\begin{equation}
\frac{J-F(z)JF(z')^*}{z+\overline{z'}}=C(I_\gamma-zA)^{-1}H^{-1}(I_\gamma-z' A)^{-*}C^*.
\end{equation}
In particular, $F$ is $J$-inner if and only if $H>0$.
\end{thm}
The finite-dimensional reproducing kernel Pontryagin  space $\mathcal{K}({F})$ with reproducing
kernel $$K^F(z,z')=\frac{J-{F}(z)J{F}(z')^*}{(z+\overline{z'})}$$ provides a minimal
state space realization for ${F}$: more precisely (see \cite{AD1}),
$${F}(z)=D+C(I_\gamma-zA)^{-1}zB,$$
where
$$
\begin{pmatrix}A&B\\C&D\end{pmatrix}:\ \begin{pmatrix}{\mathcal K}({F})\\{\mathbb C}^q\end{pmatrix}
\rightarrow
\begin{pmatrix}{\mathcal K}({F})\\{\mathbb C}^q\end{pmatrix}
$$
is defined by
\begin{eqnarray*}
(Af)(z)=(R_0f)(z):=\frac{f(z)-f(0)}{z}, &
Bu=\frac{F(z)-F(0)}{z}u,\\
Cf=f(0), &
Dx=F(0)x.
\end{eqnarray*}

Another topic considered in \cite{AG1} and \cite{AD1} is
$J$-unitary factorization. Given a matrix-valued function $F$ which is
$J$-unitary on the imaginary line one looks for all minimal
factorizations of $F$ (see \cite{BGK}) into factors which are themselves
$J$-unitary on the imaginary line. There are two equivalent
characterizations of these factorizations: the first one uses the
theory of realization and the second one uses the theory of
reproducing kernel Pontryagin spaces.
\begin{thm}
\label{voltaire}
Let $F$ be a rational matrix-valued function which is $J$-unitary on the
imaginary line and holomorphic in a neighbourhood of the origin, and
let $F(z)=D+C(I_\gamma-zA)^{-1}zB$ be a minimal realization of $F$, with
the associated Hermitian matrix $H$. There is a one-to-one correspondence
between minimal $J$-unitary factorizations of $F$ (up to a multiplicative
$J$-unitary constant) and $A$-invariant subspaces which are
non-degenerate in the (possibly, indefinite) metric induced by $H$.
\end{thm}

In general, $F$ may fail to have non-trivial $J$-unitary
factorizations.

\begin{thm}
Let $F$ be a rational matrix-valued function which is $J$-unitary on the
imaginary line and holomorphic in a neighbourhood of the origin.
There is a one-to-one correspondence between minimal $J$-unitary
factorizations of $F$ (up to a multiplicative $J$-unitary constant) and
$R_0$-invariant non-degenerate subspaces of ${\mathcal K}(F)$.
\end{thm}

The arguments in the proof of Theorem \ref{th1} do not go through
in the several complex variables context. Indeed, uniqueness,
up to a similarity, of minimal realizations doesn't hold anymore (see, e.g., \cite{GR2,FM,K}). On the other hand, the
notion of realization still makes sense in the non-commutative setting, namely for non-commutative rational \emph{formal power series} (\emph{FPSs} in short), and there is a uniqueness result for minimal
realizations in this case (see \cite{B,M,BGM1}). The latter allows us to extend the notion and study of $J$-unitary matrix-valued functions to the non-commutative case. We introduce the notion of a \emph{matrix-$J$-unitary} rational FPS as a formal power series in $N$ non-commuting indeterminates which is $J\otimes I_n$-unitary on $N$-tuples of $n\times n$ skew-Hermitian versus unitary matrices for $n=1,2,\ldots$. We extend to this case the theory of minimal realizations, minimal $J$-unitary factorizations, and backward shift models in finite-dimensional de Branges--Rovnyak spaces. We also introduce, in a similar way, the notion of matrix-selfadjoint rational formal power series, and show how to deduce the related theory for them from the theory of matrix-$J$-unitary ones.

We now turn to the outline of this paper. It consists of eight
sections.  Section~\ref{sec:intr} is this introduction. In Section~\ref{sec:prelim} we
review various results in the theory of FPSs. Let us note that the theorem on null spaces for matrix substitutions and its corollary, from our paper \cite{AK}, which are recollected in the end of Section~\ref{sec:prelim}, become an important tool in our present work on FPSs.
In Section~\ref{sec:ocm} we study the properties of
observability, controllability and minimality of Givone-Roesser nodes in the
non-commutative setting and give the corresponding criteria in terms of matrix evaluations for their ``formal transfer functions". We also formulate  a theorem on minimal factorizations of a rational FPS. In Section~\ref{sec:line} we define the
non-commutative analogue of the imaginary line and study
matrix-$J$-unitary FPSs for this case. We in particular obtain a
non-commutative version of Theorem~\ref{th1}. We obtain a counterpart of the
Lyapunov equation \eqref{st-ambroise} and of Theorem~\ref{voltaire} on minimal $J$-unitary factorizations. The unique solution of the Lyapunov equation has in this case a block diagonal structure: $H={\rm diag}(H_1,\ldots,H_N)$, and is said to be \emph{the associated structured Hermitian matrix} (associated with a given minimal realization of a matrix-$J$-unitary FPS). Section~\ref{sec:circle} contains the analogue of the previous section for the case of a
non-commutative counterpart of the unit circle. These two sections do not take into account
a counterpart of condition \eqref{croix de Chavaux}, which is considered
in Section~\ref{sec:inner} where we study matrix-$J$-inner rational FPSs. In particular, we show that the associated structured Hermitian matrix $H={\rm diag}(H_1,\ldots,H_N)$ is strictly positive in this case, which generalizes the statement in Theorem~\ref{th1} on $J$-inner functions.  We define non-commutative counterparts of the right halfplane and the unit disk, and formulate our results for both of these domains. The second one is the disjoint union of the products of $N$ copies of $n\times n$ matrix unit disks, $n=1,2,\ldots$, and plays a role of a ``non-commutative polydisk". In Theorem~\ref{thm:nc-schur} we show that any (not necessarily rational) FPS with operator coefficients, which takes contractive values in this domain, belongs to the non-commutative Schur--Agler class, defined by J.~A.~Ball, G.~Groenewald and T.~Malakorn in \cite{BGM2}. (The opposite is trivial: any function from this class has the above-mentioned property.) In other words, the contractivity of values of a FPS on $N$-tuples of strictly contractive $n\times n$ matrices,  $n=1,2,\ldots$, is sufficient for the contractivity of its values on $N$-tuples of strictly contractive operators in an infinite-dimensional separable Hilbert space. Thus, matrix-inner rational FPSs (i.e., matrix-$J$-inner ones for the case $J=I_q$) belong to the non-commutative Schur--Agler class. For this case, we recover the theorem on unitary realizations for FPSs from the latter class which was obtain in \cite{BGM2}. Moreover, our Theorem~\ref{thm:consprime} establishes the existence of a \emph{minimal}, thus \emph{finite-dimensional}, unitary Givone--Roesser realization of a rational matrix-inner FPS and \emph{the uniqueness of such a realization up to a unitary similarity}. This implies, in particular, non-commutative Lossless Bounded Real Lemma (see \cite{Ran,AG1} for its one-variable counterpart). A non-commutative version of standard Bounded Real Lemma (see \cite{ZDG}) has been presented recently in \cite{BGM3}. In Section~\ref{sec:sa} we study matrix-selfadjoint rational FPSs. In
Section~\ref{sec:spec} we introduce non-commutative formal reproducing kernel Pontryagin spaces in a way which extends one that  J.~A.~Ball and V.~Vinnikov  have introduced in \cite{BV} non-commutative formal reproducing kernel Hilbert spaces. We describe minimal backward shift realizations in
non-commutative formal reproducing kernel Pontryagin spaces which serve as a counterpart of finite-dimensional de Branges--Rovnyak spaces. Let us note that we derive an explicit formula \eqref{kfk-formal} for the corresponding reproducing kernels. In the last subsection of Section~\ref{sec:spec} we present  examples of matrix-inner rational FPSs with scalar coefficients, in two non-commuting indeterminates,  and the corresponding reproducing kernels computed by formula \eqref{kfk-formal}.

\section{Preliminaries}\label{sec:prelim}
In this section we introduce the notations which will be used throughout this paper and review
some definitions from the theory of formal power series. The symbol
${\mathbb C}^{p\times q}$ denotes the set of $p\times q$ matrices with
complex entries, and $\left({\mathbb C}^{r\times s}\right)^{p\times q}$ is the
space of $p\times q$ block matrices with block entries in
${\mathbb C}^{r\times s}$. The tensor product $A\otimes B$ of  matrices
$A\in{\mathbb C}^{r\times s}$ and $B\in \mathbb{C}^{p\times q}$ is the
element of $\left({\mathbb C}^{r\times s}\right)^{p\times q}$ with $(i,j)$-th
block entry equal to $Ab_{ij}$. The tensor product ${\mathbb C}^{r\times s}
\otimes {\mathbb C}^{p\times q}$ is the linear span of finite sums of the
form $C=\sum_{k=1}^nA_k\otimes B_k$ where $A_k\in{\mathbb C}^{r\times s}$ and
$B_k\in{\mathbb C}^{p\times q}$. One identifies ${\mathbb C}^{r\times s}
\otimes {\mathbb C}^{p\times q}$ with $\left({\mathbb C}^{r\times s}\right)^{p
\times q}$. Different representations for an element $C\in
{\mathbb C}^{r\times s}\otimes {\mathbb C}^{p\times q}$ can be reduced to a
unique one:
$$
C=\sum_{\mu=1}^r\sum_{\nu=1}^s\sum_{\tau=1}^p\sum_{\sigma=1}^q
c_{\mu\nu\tau\sigma}E'_{\mu\nu}\otimes E''_{\tau\sigma},
$$
where the matrices $E'_{\mu\nu}\in{\mathbb C}^{r\times s}$ and
$E''_{\tau\sigma}\in{\mathbb C}^{p\times q}$ are given by
$$
\left(E'_{\mu\nu}\right)_{ij}=\begin{cases}1&{\rm if}\quad(i,j)=
(\mu,\nu)\\
0&{\rm if}\quad(i,j)\not =
(\mu,\nu)
\end{cases},\qquad \mu,i=1,\ldots, r\quad{\rm and}\quad
\nu,j=1,\ldots s,
$$
$$
\left(E''_{\tau\sigma}\right)_{k\ell}=
\begin{cases}1&{\rm if}\quad(k,\ell)=
(\tau,\sigma)\\
0&{\rm if}\quad(k,\ell)\not = (\tau,\sigma)
\end{cases},
\qquad \tau,k=1,\ldots, p\quad{\rm and}\quad
\sigma,\ell=1,\ldots q.
$$
We denote by $\mathcal{F}_N$ the free semigroup with $N$ generators $g_1,
\ldots, g_N$ and the identity element $\emptyset$ with respect to the concatenation
product. This means that the generic element of $\mathcal{F}_N$ is a word
$w=g_{i_1}\cdots g_{i_n}$, where $i_\nu\in\left\{1,\ldots, N\right\}$
for $\nu=1,\ldots, n$, the identity element $\emptyset$ corresponds to the empty word,
and for another word $w'=g_{j_1}\cdots g_{j_m}$, one defines the product as
$$
ww'=g_{i_1}\cdots g_{i_n}g_{j_1}\cdots g_{j_m},\quad w\emptyset=\emptyset w=w.$$
We denote by $w^T=g_{i_n}\cdots g_{i_1}\in\mathcal{F}_N$ the
\emph{transpose} of $w=g_{i_1}\cdots g_{i_n}\in\mathcal{F}_N$
and by
$|w|=n$ the \emph{length} of the word $w$. Correspondingly, $\emptyset^T=\emptyset$, and $|\emptyset |=0$.

A  \emph{formal power series} (\emph{FPS} in short) in non-commuting indeterminates $z_1,\ldots,z_N$ with
coefficients in a linear space $\mathcal{E}$ is given by
\begin{equation}
\label{FPS}
f(z)=\sum_{w\in\mathcal{F}_N}f_wz^w,\quad f_w\in\mathcal{E},
\end{equation}
where for $w=g_{i_1}\cdots g_{i_n}$ and $z=(z_1,\ldots,z_N)$ we set $z^w=z_{i_1}\cdots z_{i_n}$, and $z^{\emptyset}=1$.
 We denote by $\mathcal{E}\left\langle \left\langle z_1,\ldots,z_N\right\rangle\right\rangle$ the linear space of FPSs in non-commuting indeterminates $z_1,\ldots,z_N$ with coefficients in $\mathcal{E}$.
A series $f\in\mathbb{C}^{p\times q}\left\langle \left\langle z_1,\ldots,z_N\right\rangle\right\rangle$ of the form \eqref{FPS} can also be viewed as a $p\times q$ matrix whose entries
are formal power series with coefficients in ${\mathbb C}$, i.e., belong to the space $\mathbb{C}\left\langle \left\langle z_1,\ldots,z_N\right\rangle\right\rangle$, which has an additional structure of non-commutative ring (we assume that the indeterminates $z_j$ formally commute with the coefficients $f_w$). The \emph{support}
of a FPS $f$ given by \eqref{FPS} is the set
$$
{\rm supp}\,f=\left\{w\in\mathcal{F}_N:\ f_w\not=0\right\}.
$$
\emph{Non-commutative polynomials} are formal power series with finite
support. We denote by $\mathcal{E}\left\langle z_1,\ldots,z_N\right\rangle$ the subspace in the space $\mathcal{E}\left\langle \left\langle z_1,\ldots,z_N\right\rangle\right\rangle$ consisting of non-commutative polynomials. Clearly, a FPS is determined by its coefficients $f_w$.
Sums and products of two FPSs $f$ and $g$  with matrix coefficients of compatible sizes (or with operator coefficients) are given by
\begin{equation}
\label{OPS}
\left(f+g\right)_w=f_w+g_w,\qquad (fg)_w=
\sum_{w'w''=w} f_{w'}g_{w''}.
\end{equation}

A FPS $f$ with coefficients in ${\mathbb C}$ is invertible if
and only if $f_\emptyset\not=0$. Indeed, assume that $f$ is invertible. From
the definition of the product of two FPSs in \eqref{OPS} we get
$f_\emptyset (f^{-1})_\emptyset=1,$
and hence $f_\emptyset\not=0$. On the other hand, if $f_\emptyset\not=0$ then
$f^{-1}$ is given by
$$
f^{-1}(z)=\sum_{k=0}^\infty \left(1-
f_\emptyset^{-1}f(z)\right)^kf_\emptyset^{-1}.
$$
The formal power series in the right-hand side is well defined since the
expansion of $\left(1-f_\emptyset^{-1}f\right)^k$ contains words of length
at least $k$, and thus the coefficients  $(f^{-1})_w$ are finite
sums.

A FPS with coefficients in ${\mathbb C}$ is called
\emph{rational} if it can be expressed as a finite number of sums, products and
inversions of non-commutative polynomials.
A formal power series with coefficients in ${\mathbb C}^{p\times q}$ is called
\emph{rational} if it is a $p\times q$ matrix whose all entries are rational
FPSs with coefficients in ${\mathbb C}$. We will denote by
${\mathbb C}^{p\times q}\left\langle \left\langle z_1,\ldots, z_N\right\rangle\right\rangle_{\rm rat}$
the linear space of rational
FPSs with coefficients in ${\mathbb C}^{p\times q}$. 
Define the product of $f\in{\mathbb C}^{p\times q}\left\langle \left\langle z_1,\ldots, z_N\right\rangle\right\rangle_{\rm rat}$ and
$p\in{\mathbb C}\left\langle z_1,\ldots, z_N\right\rangle$ as follows:
\begin{enumerate}
\item $f\cdot 1=f$ for every $f\in{\mathbb C}^{p\times q}\left\langle \left\langle z_1,\ldots, z_N\right\rangle\right\rangle_{\rm rat}$;
\item For every word $w'\in\mathcal{F}_N$ and every
$f\in{\mathbb C}^{p\times q}\left\langle \left\langle z_1,\ldots, z_N\right\rangle\right\rangle_{\rm rat}$,
$$f\cdot z^{w'}=\sum_{w\in\mathcal{F}_N}f_w z^{ww'}
=\sum_w f_vz^w$$
where the last sum is taken over all $w$ which can be written as
$w=vw'$ for some $v\in\mathcal{F}_N$;
\item For every $f\in{\mathbb C}^{p\times q}\left\langle \left\langle z_1,\ldots, z_N\right\rangle\right\rangle_{\rm rat}$,
$p_1,p_2\in{\mathbb C}\left\langle z_1,\ldots, z_N\right\rangle$ and $\alpha_1,\alpha_2
\in{\mathbb C}$,
$$f\cdot(\alpha_1 p_1+\alpha_2 p_2)=
\alpha_1 (f\cdot p_1)+\alpha_2 (f\cdot p_2).$$
\end{enumerate}
The space ${\mathbb C}^{p\times q}\left\langle \left\langle z_1,\ldots, z_N\right\rangle\right\rangle_{\rm rat}$ is a right
module over the ring ${\mathbb C}\left\langle z_1,\ldots, z_N\right\rangle$ with respect to this
product. A structure of left ${\mathbb C}\left\langle z_1,\ldots, z_N\right\rangle$-module
 can be defined in a similar
way since the indeterminates commute with coefficients.

Formal power series are used in various branches of mathematics, e.g., in abstract algebra, enumeration problems and
combinatorics; rational
formal power series have been extensively used in theoretical
computer science, mostly in automata theory and language theory (see
\cite{BR}). The Kleene--Sch\"{u}tzenberger theorem \cite{Kle,Sch} (see also \cite{F}) says
that a FPS $f$ with coefficients in ${\mathbb
C}^{p\times q}$ is rational if and only if it is \emph{recognizable}, i.e., there exist $r\in{\mathbb
N}$ and matrices
$C\in{\mathbb C}^{p\times r},\ A_1,\ldots, A_N\in
{\mathbb C}^{r\times r}$ and $B\in{\mathbb C}^{r\times q}$
such that for every word $w=g_{i_1}\cdots g_{i_n}\in\mathcal{F}_N$
one has
\begin{equation}
\label{rec}
f_w=CA^w B,\quad{\rm where}\quad A^w=A_{i_1}\ldots A_{i_n}.
\end{equation}
Let $\mathcal{H}_f$ be the \emph{Hankel matrix} whose rows and columns are indexed
by the words of $\mathcal{F}_N$ and defined by
$$
\left(\mathcal{H}_f\right)_{w,w'}=f_{ww'^T},\quad w,w'\in\mathcal{F}_N.$$

It follows from \eqref{rec} that if the FPS $f$ is recognizable then
$\left(\mathcal{H}_f\right)_{w,w'}=CA^{ww'^T}B$ for all $w, w'\in
\mathcal{F}_N$. M.~Fliess has shown in \cite{F} that a FPS $f$ is rational (that is, recognizable) if and only if
$$\gamma:={\rm rank}\,\mathcal{H}_f<\infty.$$
In this case the number $\gamma$ is the smallest possible $r$ for a representation \eqref{rec}.

In control theory, rational FPSs appear as the input/output
mappings of linear systems with structured uncertainties. For instance, in \cite{BD} a system matrix is given by
$$M=\begin{pmatrix}
A&B\\ C&D\end{pmatrix}\in{\mathbb C}^{(r+p)\times (r+q)},
$$
and the uncertainty operator is given by
$$
\Delta(\delta)={\rm diag}(\delta_1I_{r_1},\ldots,\delta_NI_{r_N}),$$
where $r_1+\cdots +r_N=r$. The uncertainties $\delta_k$ are linear
operators on $\ell^2$ representing disturbances or small perturbation
parameters which enter the system at different locations. Mathematically,
they can be interpreted as non-commuting indeterminates. The input/output
map is a linear fractional transformation
\begin{equation}
\label{lft}
LFT(M,\Delta(\delta))=D+C(I_r-\Delta(\delta)A)^{-1}\Delta(\delta)B,
\end{equation}
which can be interpreted as a non-commutative
transfer function $T^{\rm nc}_\alpha$ of a linear system $\alpha$ with
evolution on $\mathcal{F}_N$:
\begin{equation}
\label{voltaire1}
\alpha:\,\begin{cases}
x_j(g_jw)&=A_{j1}x_1(w)+\cdots+A_{jN}x_N(w)+B_ju(w), \quad j=1,\ldots, N,\\
y(w)&=C_{1}x_1(w)+\cdots+C_{N}x_N(w)+Du(w),
\end{cases}
\end{equation}
where $x_j(w)\in{\mathbb C}^{r_j}$ ($j=1,\ldots, N$), $u(w)\in
{\mathbb C}^q$, $y(w)\in{\mathbb C}^p$, and the matrices $A_{jk},B$ and $C$ are of
appropriate sizes along the decomposition
${\mathbb C}^{r}={\mathbb C}^{r_1}\oplus\cdots\oplus{\mathbb C}^{r_N}$. Such a system appears in \cite{M,BGM1,BGM2,BGM3} and
is known as the \emph{non-commutative Givone--Roesser model} of
multidimensional linear system; see \cite{GR1,GR2,R} for its commutative counterpart.

In this paper we do not consider system evolutions (i.e., equations
\eqref{voltaire1}). We will use the terminology
\emph{$N$-dimensional Givone--Roesser operator node}  (for brevity, \emph{GR-node})
for the collection of data
\begin{equation}
\label{GR}
\alpha=(N;A,B,C,D;{\mathbb C}^r=\bigoplus_{j=1}^N{\mathbb C}^{r_j},
{\mathbb C}^q,{\mathbb C}^p).
\end{equation}
Sometimes instead of spaces $\mathbb{C}^r,\mathbb{C}^{r_j}\ (j=1,\ldots,N),\mathbb{C}^q$ and $\mathbb{C}^p$ we shall consider abstract finite-dimensional linear spaces $\mathcal{X}$ (the \emph{state space}), $\mathcal{X}_j\ (j=1,\ldots,N)$, $\mathcal{U}$ (the \emph{input space}) and $\mathcal{Y}$ (the \emph{output space}), respectively, and a node
$$\alpha=(N;A,B,C,D;\mathcal{X}=\bigoplus_{j=1}^N\mathcal{X}_j, \mathcal{U},\mathcal{Y}),$$
where $A,B,C,D$ are linear operators in the corresponding pairs of spaces.
The \emph{non-commutative transfer function of a GR-node $\alpha$} is a rational FPS
\begin{equation}
\label{tf}
T^{\rm nc}_\alpha(z)=D+C(I_r-\Delta(z)A)^{-1}\Delta(z)B.
\end{equation}

Minimal GR-realizations \eqref{GR} of non-commutative rational FPSs, that is, representations of them in the form \eqref{tf}, with
minimal possible $r_k$ for $k=1,\ldots, N$ were studied in \cite{BD,B,M,BGM1}. For $k=1,\ldots, N$, the
\emph{$k$-th observability matrix} is
$$
{\mathcal O}_k={\rm col}(C_k,C_1A_{1k},\ldots, C_NA_{Nk}, C_1A_{11}A_{1k},\ldots C_1A_{1N}A_{Nk},
\ldots)$$
and the \emph{$k$-th controllability matrix} is
$$
{\mathcal C}_k={\rm row}(B_k,A_{k1}B_1,\ldots, A_{kN}B_{N},A_{k1}A_{11}B_1,\ldots A_{kN}A_{N1}B_1,
\ldots)$$
(note that these are infinite block matrices).
A GR-node $\alpha$ is called \emph{observable} (resp., \emph{controllable}) if
${\rm rank}\,{\mathcal O}_k=r_k$ (resp., ${\rm rank}\,{\mathcal C}_k=r_k$)
for $k=1,\ldots, N$. A GR-node $\alpha=(N;A,B,C,D;{\mathbb C}^r=\bigoplus_{j=1}^N{\mathbb C}^{r_j},
{\mathbb C}^q,{\mathbb C}^p)$ is observable if and only if its \emph{adjoint GR-node} $\alpha^*=(N;A^*,C^*,B^*,D^*;{\mathbb C}^r=\bigoplus_{j=1}^N{\mathbb C}^{r_j},
{\mathbb C}^p,{\mathbb C}^q)$ is controllable. (Clearly, $(\alpha^*)^*=\alpha$.)

In view of the sequel, we introduce some notations. We set:
\[
\begin{split}
A^{w g_{\nu}}&=A_{j_1j_2}A_{j_2j_3}\cdots A_{j_{k-1}j_k}A_{j_k\nu},\\
(C\flat A)^{g_{\nu} w}&=C_{\nu}
A_{\nu j_1}A_{j_1j_2}\cdots A_{j_{k-1}j_k},\\
(A\sharp B)^{w g_{\nu}}&=
A_{j_1j_2}\cdots A_{j_{k-1}j_k}A_{j_k\nu}B_\nu,\\
(C\flat A\sharp B)^{g_{\mu}  w  g_{\nu}}&=
C_{\mu}
A_{\mu j_1}A_{j_1j_2}\cdots A_{j_{k-1}j_k}A_{j_k\nu}B_\nu,
\end{split}
\]
where $w=g_{j_1}\cdots g_{j_k}\in\mathcal{F}_N$ and $\mu,\nu
\in\left\{1,\ldots, N\right\}$. We also define:
\[
\begin{split}
A^{g_{\nu}}=A^{\emptyset} &=I_{\gamma}\\
(C\flat A)^{g_{\nu}}&=C_{\nu},\\
(A\sharp B)^{g_{\nu}}&=B_{\nu},\\
(C\flat A\sharp B)^{g_{\nu}}&=C_{\nu} B_{\nu},\\
(C\flat A\sharp B)^{g_{\mu} g_{\nu}}&=C_{\mu} A_{\mu\nu}B_{\nu},
\end{split}
\]
and hence, with the lexicographic order of words in $\mathcal{F}_N$,
$${\mathcal O}_k={\rm col}_{w\in\mathcal{F}_N}
(C\flat A)^{w g_k}\quad{\rm and}\quad
{\mathcal C}_k={\rm row}_{w\in\mathcal{F}_N}(A\sharp B)^{g_kw^T },$$
and the coefficients of the FPS $T_\alpha^{{\rm nc}}$
(defined by \eqref{tf}) are given by
$$
(T_\alpha^{\rm nc})_\emptyset=D,\qquad
(T_\alpha^{\rm nc})_w=(C\flat A\sharp B)^w\quad{\rm for}\quad
w=g_{j_1}\cdots g_{j_n}\in\mathcal{F}_N.$$
The \emph{$k$-\emph{th} Hankel matrix associated with a FPS} $f$ is
defined in \cite{M} (see also \cite{BGM1}) as
$$
(\mathcal{H}_{f,k})_{w,w'g_k}=f_{w g_k w'^T}\quad {\rm with}\quad
w,w'\in\mathcal{F}_N,$$
that is, the rows of $\mathcal{H}_{f,k}$ are indexed by all the words of
$\mathcal{F}_N$ and the columns of $\mathcal{H}_{f,k}$ are indexed by all the words of
$\mathcal{F}_N$ ending by $g_k$, provided the lexicographic order is
used. If a GR-node $\alpha$ defines a realization of $f$, that is,
$f=T_\alpha^{{\rm nc}}$, then
$$
(\mathcal{H}_{f,k})_{w,w'g_k}=(C\flat A\sharp B)^{w g_kw'^T}=
(C\flat A)^{w g_k}(A\sharp B)^{g_kw'^T},$$
i.e., $\mathcal{H}_{f,k}={\mathcal O}_k{\mathcal C}_k$. Hence,
the node $\alpha$ is minimal if and only if $\alpha$ is both observable and controllable, i.e.,
$$\gamma_k:={\rm rank}\,\mathcal{H}_{f,k}=r_k\quad{\mbox{\rm for all}}\quad k\in\left\{1,\ldots, N
\right\}.
$$
This last set of conditions is an analogue of the above mentioned result of
Fliess on minimal recognizable representations of rational formal power
series. Every non-commutative rational FPS has a minimal GR-realization.

Finally, we note (see \cite{BD,M}) that two minimal GR-realizations of a given rational
FPS are \emph{similar}: if
$\alpha^{(i)}=(N;A^{(i)},B^{(i)},C^{(i)},D;{\mathbb C}^{\gamma}=
\bigoplus_{k=1}^N{\mathbb C}^{\gamma_k},{\mathbb C}^q,
{\mathbb C}^p)$ (i=1,2) are minimal GR-nodes  such that $T_{\alpha^{(1)}}^{\rm nc}=
T_{\alpha^{(2)}}^{\rm nc}$ then there exists a block diagonal invertible matrix
$T={\rm diag}(T_1,\ldots , T_N)$ (with $T_k\in{\mathbb C}^{\gamma_k\times \gamma_k}$)
such that
\begin{equation}\label{sim}
A^{(1)}=T^{-1}A^{(2)}T,\quad
B^{(1)}=T^{-1}B^{(2)},\quad
C^{(1)}=C^{(2)}T.
\end{equation}
Of course, the converse is also true, moreover, any two similar (not necessarily minimal) GR-nodes have the same transfer functions.

Now we turn to the discussion on substitutions of matrices for indeterminates in formal power series. Many properties of non-commutative FPSs or non-commutative polynomials are described in terms of matrix substitutions, e.g., matrix-positivity of non-commutative polynomials (non-commutative Positivestellensatz) \cite{H,McC,HMcC,HMcCP}, matrix-positivity of FPS kernels \cite{KV}, matrix-convexity \cite{CHSY,H1}. The non-commutative Schur--Agler class, i.e., the class of FPSs with operator coefficients, which take contractive values on all $N$-tuples of strictly contractive operators on $\ell^2$, was studied in \cite{BGM2} \footnote{In fact, a more general class was studied in \cite{BGM2}, however for our purposes it is enough to consider here only the case mentioned above.}; we will show in Section~\ref{sec:inner} that in order that a FPS belongs to this class it suffices to check its contractivity on $N$-tuples of strictly contractive $n\times n$ matrices, for all  $n\in\mathbb{N}$. The notions of matrix-$J$-unitary (in particular, matrix-$J$-inner) and matrix-selfadjoint rational FPS, which will be introduced and studied in the present paper, are also defined in terms of substitutions of matrices (of a certain class) for indeterminates.

Let $p(z)=\sum_{|w|\leq m}p_wz^w\in\mathbb{C}\left\langle z_1,\ldots,z_N\right\rangle$. For $n\in\mathbb{N}$ and an $N$-tuple of matrices $Z=(Z_1,\ldots,Z_N)\in\left(\mathbb{C}^{n\times n}\right)^N$, set
 $$p(Z)=\sum_{|w|\leq m}p_wZ^w,$$
where $Z^w=Z_{i_1}\cdots Z_{i_{|w|}}$ for $w=g_{i_1}\cdots g_{i_{|w|}}\in\mathcal{F}_N$, and $Z^{\emptyset}=I_n$. Then for any rational expression for a FPS $f\in \mathbb{C}\left\langle\left\langle  z_1,\ldots,z_N\right\rangle\right\rangle_{\rm rat}$ its value at $Z\in\left(\mathbb{C}^{n\times n}\right)^N$ is well defined provided all of the inversions of polynomials $p^{(j)}\in\mathbb{C}\left\langle z_1,\ldots,z_N\right\rangle$ in this expression are well defined at $Z$. The latter is the case at least in some neighbourhood of $Z=0$, since $p_{\emptyset}^{(j)}\not =0$.

Now, if $f\in\mathbb{C}^{p\times q}\left\langle\left\langle  z_1,\ldots,z_N\right\rangle\right\rangle_{\rm rat}$ then the value $f(Z)$ at some $Z\in\left(\mathbb{C}^{n\times n}\right)^N$ is well defined whenever the values of matrix entries $\left(f_{ij}(Z)\right)\ (i=1,\ldots,p;j=1,\ldots,q)$ are well defined at $Z$. As a function of matrix entries $(Z_k)_{ij}\ (k=1,\ldots,N;i,j=1,\ldots,n)$, $f(Z)$ is rational $\mathbb{C}^{p\times q}\otimes\mathbb{C}^{n\times n}$-valued function, which is holomorphic on an open and dense set in $\mathbb{C}^{n\times n}$. The latter set contains some neighbourhood \begin{equation}\label{gamma}\Gamma_n(\varepsilon):=\{ Z\in\left(\mathbb{C}^{n\times n}\right)^N:\ \| Z_k\| <\varepsilon,\ k=1,\ldots,N\}\end{equation} of $Z=0$, where $f(Z)$ is given by
$$f(Z)=\sum_{w\in\mathcal{F}_N}f_w\otimes Z^w.$$
The following results from \cite{AK} on matrix substitutions  are used in the sequel.
\begin{thm}\label{thm:kers}
Let $f\in\mathbb{C}^{p\times q}\left\langle\left\langle
  z_1,\ldots,z_N\right\rangle\right\rangle_{\rm rat}$, and $m\in\mathbb{Z}_+$ be such that \[\bigcap_{w\in\mathcal{F}_N:|w|\leq m}\ker f_w=\bigcap_{w\in\mathcal{F}_N}\ker f_w.\] Then there exists $\varepsilon>0$ such that for every $n\in\mathbb{N}:\ n\geq m^m$ (in the case $m=0$, for every $n\in\mathbb{N}$),
\begin{equation}\label{kers}
\bigcap\limits_{Z\in\Gamma_n(\varepsilon)}\ker f(Z)=\left(\bigcap\limits_{w\in\mathcal{F}_N:\ |w|\leq m}\ker f_w\right)\otimes\mathbb{C}^n,
\end{equation}
and moreover, there exist $l\in\mathbb{N}:\ l\leq qn$, and $N$-tuples of matrices $Z^{(1)},\ldots, Z^{(l)}$ from $\Gamma_n(\varepsilon)$ such that
        \[ \bigcap\limits_{j=1}^l\ker f(Z^{(j)})=\left(\bigcap\limits_{w\in\mathcal{F}_N:\ |w|\leq m}\ker f_w\right)\otimes\mathbb{C}^n.\]
\end{thm}
\begin{cor}\label{cor:vanish-fps}
In conditions of Theorem~\ref{thm:kers}, if for some $n\in\mathbb{N}:\ n\geq m^m$ (in the case $m=0$, for some $n\in\mathbb{N}$) one has $f(Z)=0,\ \forall Z\in\Gamma_n(\varepsilon)$, then $f=0$.
\end{cor}

\section{More on observability, controllability, and mi\-nimality in the
non-commutative setting}\label{sec:ocm}
In this section we prove a number of results on observable, controllable and
minimal GR-nodes in the multivariable non-commutative
setting, which generalize some well known statements for one-variable
nodes (see \cite{BGK}).

Let us introduce the \emph{$k$-th truncated observability matrix}
$\widetilde{\mathcal O}_k$ and the \emph{$k$-th truncated controllability
matrix} $\widetilde{\mathcal C}_k$ of a GR-node \eqref{GR}
by
\begin{equation*}
\widetilde{\mathcal O_k}={\rm col}_{|w|< pr}(C\flat A)^{w
g_k},\quad \widetilde{{\mathcal C}_k}={\rm row}_{|w|< rq} (A\sharp
B)^{g_kw^T},
\end{equation*}
with the lexicographic order of words in $\mathcal{F}_N$.
\begin{thm} \label{thm:trunc} For each $k\in\{ 1,\ldots,N\}$:
${\rm rank}\,\widetilde{{\mathcal O}_k}={\rm rank}\,{\mathcal O_k}$ and
${\rm rank}\,\widetilde{{\mathcal C}_k}={\rm rank}\,{\mathcal C}_k$.
\end{thm}
\begin{proof}
Let us show that for every fixed $k\in\left\{1,\ldots, N\right\}$
matrices of the form $(C\flat A)^{w g_k}$ with $|w|\ge pr$ are
representable as linear combinations of matrices $(C\flat
A)^{{\widetilde w}g_k}$ with $|\widetilde{w}|<pr$. First we remark
that if for each fixed $k\in\left\{1,\ldots, N\right\}$ and
$j\in{\mathbb N}$ all matrices of the form $(C\flat A)^{w g_k}$
with $|w|=j$ are representable as linear combinations of matrices
of the form $(C\flat A)^{w'g_k}$ with $|w'|<j$ then the same holds
for matrices of the form $(C\flat A)^{w g_k}$ with $|w|= j+1$.
Indeed, if $w=i_1\cdots i_ji_{j+1}$ then there exist words
$w_1',\ldots, w_s'$ with $ |w_1'|<j,\ldots,|w_s'|<j$ and
$a_1,\ldots, a_s\in {\mathbb C}$ such that $$ (C\flat
A)^w=\sum_{\nu=1}^s a_\nu(C\flat A)^{w_\nu' g_{i_{j+1}}}.$$ Then
for every $k\in\left\{1,\ldots, N\right\}$,
\begin{eqnarray*}
\lefteqn{(C\flat A)^{w g_k}= (C\flat A)^w A_{i_{j+1},k}
=\sum_{\nu=1}^s a_\nu (C\flat A)^{w'_\nu g_{i_{j+1}}}
A_{i_{j+1},k}}\\
&=& \sum_{\nu:\,|w_\nu'|<j-1}a_\nu(C\flat A)^{w_\nu'
g_{i_{j+1}}}A_{i_{j+1},k}+
\sum_{\nu:\,|w_\nu'|=j-1}a_\nu(C\flat A)^{w_\nu'
g_{i_{j+1}}}A_{i_{j+1},k}\\
&=& \sum_{\nu:\,|w_\nu'|<j-1}a_\nu(C\flat A)^{w_\nu'
g_{i_{j+1}}g_k}+
\sum_{\nu:\,|w_\nu'|=j-1}a_\nu(C\flat A)^{w_\nu'
g_{i_{j+1}}g_k}.
\end{eqnarray*}
Consider these two sums separately. All the terms in the first sum
are of the form $a_\nu(C\flat A)^{(w_\nu' g_{i_{j+1}})g_k}$ with
$|w_\nu' g_{i_{j+1}}|<j$. In the second sum, by the assumption,
for each matrix $(C\flat A)^{w'_\nu g_{i_{j+1}}g_k}$ there exist
words $w_{1\nu}'',\ldots, w_{t\nu}''$ of length  strictly less
than $j$ and complex numbers $b_{1\nu},\ldots, b_{t\nu}$ such that
$$(C\flat A)^{w'_\nu g_{i_{j+1}}g_k}=\sum_{\mu=1}^tb_{\mu\nu}
(C\flat A)^{w''_{\mu\nu}g_k}. $$ Hence $(C\flat A)^{w g_k}$ is a
linear combination of matrices of the form $(C\flat
A)^{\widetilde{w} g_k}$ with $|\widetilde{w}|<j$. Reiterating this
argument we obtain that any matrix of the form $(C\flat A)^{w
g_k}$ with $|w|\ge j$ and fixed $k\in\left\{ 1,\ldots, N\right\}$
can be represented as a linear combination of matrices of the form
$(C\flat A)^{\widetilde{w} g_k}$ with $|\widetilde{w}|<j$. In
particular,
\begin{equation}\label{rank}{\rm rank\, col}_{|w|< j}(C\flat A)^{wg_k}={\rm
rank}\,\mathcal{O}_k,\quad k=1,\ldots,N.\end{equation} Since for
any $k\in\{ 1,\ldots,N\}$ one has $(C\flat A)^{wg_k}\in{\mathbb
C}^{p\times r_k}$ and ${\rm dim}\,{\mathbb C}^{p\times r_k}=pr_k$,
we obtain that for some $j\le pr$, and moreover for $j=pr$
\eqref{rank} is true, i.e., ${\rm rank}\,\widetilde{{\mathcal
O}_k}= {\rm rank}\,{{\mathcal O}_k}$.

The second equality is proved analogously.
\end{proof}
\begin{rem}
The sizes of the truncated matrices $\widetilde{{\mathcal O}_k}$
and $\widetilde{{\mathcal C}_k}$ depend only on the sizes of matrices
$A,B$ and $C$, and do not depend on these matrices themselves. Our estimate
for the size of $\widetilde{{\mathcal O}_k}$  is rough, and one could probably
improve it. For our present purposes, only the finiteness of the matrices
$\widetilde{{\mathcal O}_k}$ and $\widetilde{{\mathcal C}_k}$ is important, and
not their actual sizes.\end{rem}
\begin{cor}\label{cor:co}
A GR-node \eqref{GR} is observable (resp., controllable) if and
only if for every $k\in\left\{1,\ldots, N\right\}$:
$${\rm rank}\,\widetilde{{\mathcal O}_k}=r_k\quad
(resp,\ {\rm rank}\,\widetilde{{\mathcal C}_k}=r_k),$$
or equivalently, the matrix ${\mathcal O}_k$ (resp., ${\mathcal C}_k$) is left (resp., right) invertible.
\end{cor}
\begin{rem}
Corollary~\ref{cor:co} is comparable with Theorems~7.4 and 7.7 in \cite{M}, however we note again that the matrices $\widetilde{{\mathcal O}_k}$ and $\widetilde{{\mathcal C}_k}$ here are finite.
\end{rem}
\begin{thm}
\label{thm:sim}
Let
$\alpha^{(i)}=(N;A^{(i)},B^{(i)},C^{(i)},D,{\mathbb C}^{\gamma}=
\oplus_{k=1}^N{\mathbb C}^{{\gamma}_k},{\mathbb C}^q,
{\mathbb C}^p),\ i=1,2$, be minimal GR-nodes with the same
transfer function. Then they are similar, the similarity transform is unique and given by   $T={\rm diag}(T_1,\ldots, T_N)$ where
\begin{equation}\label{t}
T_k=\left(\widetilde{\mathcal O}_k^{(2)}\right)^+
\widetilde{\mathcal O}_k^{(1)}=
\widetilde{\mathcal C}_k^{(2)}
\left(\widetilde{\mathcal C}_k^{(1)}\right)^\dag
\end{equation}
(here ``$+$" denotes a left inverse, while
``$\dag$" denotes a right inverse).
\end{thm}
\begin{proof} We already mentioned in Section~\ref{sec:prelim} that two minimal
nodes with the same transfer function are similar. Let
$T'={\rm diag}\,(T_1',\ldots, T_N')$ and
$T''={\rm diag}\,(T''_1,\ldots,
T''_N)$
be two similarity transforms. Let $x\in{\mathbb C}^{\gamma_k}$. Then, for
every $w\in\mathcal{F}_N$,
$$
(C^{(2)}\flat A^{(2)})^{w g_k}\left(T_k''-
T_k'\right)x=(C^{(1)}\flat A^{(1)})^{w g_k}x-
(C^{(1)}\flat A^{(1)})^{w g_k}x=0.$$
Since $x$ is arbitrary, from the observability of $\alpha^{(2)}$ we get
$T_k'=T_k''$ for $k=1,\ldots, N$, hence the
similarity transform is unique. Comparing the coefficients in the two
FPS representations of the transfer function, we obtain
$$(C^{(1)}\flat A^{(1)}\sharp B^{(1)})^w=
(C^{(2)}\flat A^{(2)}\sharp B^{(2)})^w$$
for all of $w\in\mathcal{F}_N\setminus\left\{\emptyset\right\}$, and therefore
$$\widetilde{{\mathcal O}_k^{(1)}}\widetilde{{\mathcal C}_k^{(1)}}=
\widetilde{{\mathcal O}_k^{(2)}}\widetilde{{\mathcal C}_k^{(2)}},\quad
k=1,\ldots, N.$$
Thus we obtain
$$
\left(\widetilde{{\mathcal O}_k^{(2)}}\right)^+
\widetilde{{\mathcal O}_k^{(1)}}=
\widetilde{{\mathcal C}_k^{(2)}}
\left(\widetilde{{\mathcal C}_k^{(1)}}\right)^\dag,\quad
k=1,\ldots, N.$$
Denote the operators which appear in these equalities by $T_k,\ k=1,\ldots, N$.  A direct computation shows that $T_k$ are invertible with $$T_k^{-1}=\left(\widetilde{{\mathcal O}_k^{(1)}}\right)^+
\widetilde{{\mathcal O}_k^{(2)}}=
\widetilde{{\mathcal C}_k^{(1)}}
\left(\widetilde{{\mathcal C}_k^{(2)}}\right)^\dag.$$ Let us verify that $T={\rm diag}(T_1,\ldots, T_N)\in\mathbb{C}^{\gamma\times\gamma}$ is a similarity transform between $\alpha^{(1)}$ and $\alpha^{(2)}$. It follows from the controllability of $\alpha^{(1)}$ that for arbitrary $k\in\{ 1,\ldots,N\}$ and $x\in\mathbb{C}^{\gamma_k}$ there exist words $w_j\in\mathcal{F}_N$, with $|w_j|<\gamma q$,  scalars $a_j\in\mathbb{C}$ and vectors $u_j\in\mathbb{C}^q,\ j=1,\ldots,s$, such that
$$x=\sum_{\nu=1}^sa_\nu(A^{(1)}\sharp B^{(1)})^{g_kw_\nu^T}u_\nu.$$
Then
\begin{eqnarray*}
\lefteqn{T_kx=\left(\widetilde{{\mathcal O}_k^{(2)}}\right)^+
\widetilde{{\mathcal O}_k^{(1)}}x=\sum_{\nu=1}^sa_\nu \left(\widetilde{{\mathcal O}_k^{(2)}}\right)^+
\widetilde{{\mathcal O}_k^{(1)}}(A^{(1)}\sharp B^{(1)})^{g_kw_\nu^T}u_\nu }\\
&=&\sum_{\nu=1}^sa_\nu \left(\widetilde{{\mathcal O}_k^{(2)}}\right)^+
\widetilde{{\mathcal O}_k^{(2)}}(A^{(2)}\sharp B^{(2)})^{g_kw_\nu^T}u_\nu=\sum_{\nu=1}^sa_\nu(A^{(2)}\sharp B^{(2)})^{g_kw_\nu^T}u_\nu.
\end{eqnarray*}
This explicit formula implies the set of equalities
$$T_kB_k^{(1)}=B_k^{(2)},\quad T_kA_{kj}^{(1)}=A_{kj}^{(2)}T_j,\quad C_k^{(1)}=C_k^{(2)}T_k,\quad k,j=1,\ldots,N,$$
which is equivalent to \eqref{sim}.
\end{proof}
\begin{rem}
Theorem~\ref{thm:sim} is comparable with Theorem~7.9 in \cite{M}. However, we establish in Theorem~\ref{thm:sim} the uniqueness and an explicit formula for the similarity transform $T$.
\end{rem}
Using Theorem~\ref{thm:kers}, we will prove now the following criteria of observability, controllability,
and minimality for GR-nodes analogous to the ones  proven in \cite[Theorem 3.3]{AK}
for recognizable FPS represntations.
\begin{thm}\label{thm:com}
  A GR node $\alpha$ of the form \eqref{GR} is observable (resp., controllable)
  if and only if for every $k\in\{ 1,\ldots,N\}$ and $n\in\mathbb{N}:\ n\geq (pr-1)^{pr-1}$ (resp,
  $n\geq (rq-1)^{rq-1}$), which means in the case of $pr=1$ (resp., $rq=1$):
  ``for every $n\in\mathbb{N}$",
\begin{eqnarray}
\bigcap\limits_{Z\in\Gamma_n(\varepsilon)}
\ker \varphi_k(Z) &=& 0\label{z-o}\\
(resp., \bigvee\limits_{Z\in\Gamma_n(\varepsilon)}{\rm ran}\,\psi_k(Z) &=& \mathbb{C}^{r_k}\otimes\mathbb{C}^n),\label{z-c}
\end{eqnarray}
where the rational FPSs $\varphi_k$ and $\psi_k$ are defined by
  \begin{eqnarray}
\varphi_k(z) &=& C(I_{r}-\Delta(z)A)^{-1}\big|_{\mathbb{C}^{r_k}},\label{phi_k}\\
\psi_k(z) &=& P_k(I_{r}-A\Delta(z))^{-1}B,\label{psi_k}
\end{eqnarray}
with $P_k$ standing for the orthogonal projection onto $\mathbb{C}^{r_k}$ (which is naturally identified here with the subspace in $\mathbb{C}^r$), the symbol ``$\bigvee$" means linear
span,
 $\varepsilon=\| A\|^{-1}$ ($\varepsilon>0$ is arbitrary in the case $A=0$), and
 $\Gamma_n(\varepsilon)$ is defined by \eqref{gamma}.
 This GR-node is minimal if both of conditions
   \eqref{z-o} and \eqref{z-c} are fulfilled.
\end{thm}
\begin{proof}
First, let us remark that for all $k=1,\ldots,N$ the functions $\varphi_k$ and $\psi_k$ are well defined in $\Gamma_n(\varepsilon)$, and holomorphic as
functions of matrix entries $(Z_j)_{\mu\nu},\ j=1,\ldots,N,\ \mu,\nu=1,\ldots,n$. Second, Theorem~\ref{thm:trunc}
implies that in Theorem~\ref{thm:kers} applied to $\varphi_k$ one can choose $m=pr-1$, and then from \eqref{kers} obtain that
observability for a GR-node $\alpha$ is equivalent to condition \eqref{z-o}. Since $\alpha$ is controllable if and only if $\alpha^*$
is observable, controllability for $\alpha$ is equivalent to condition \eqref{z-c}. Since minimality for a GR-node $\alpha$ is equivalent
to controllability and observability together, it is in turn equivalent to conditions \eqref{z-o} and \eqref{z-c} together.
\end{proof}
Let $\alpha'=(N;A',B',C',D';{\mathbb C}^{r'}=\bigoplus_{j=1}^N{\mathbb C}^{r'_j},
{\mathbb C}^s,{\mathbb C}^p)$ and $\alpha''=(N;A'',B'',\\ C'',D'';{\mathbb C}^{r''}=\bigoplus_{j=1}^N{\mathbb C}^{r''_j},
{\mathbb C}^q,{\mathbb C}^s)$ be GR-nodes. For $k,j=1,\ldots,N$ set $r_j=r'_j+r''_j$, and
\begin{equation}\label{prod}
\begin{array}{ll}
        A_{kj}=\begin{pmatrix}A'_{kj} & B'_kC''_j\\0 & A''_{kj}\end{pmatrix}\in\mathbb{C}^{r_k\times r_j}, & B_k=\begin{pmatrix}B'_kD''\\B''_k\end{pmatrix}\in\mathbb{C}^{r_k\times q},\\
        C_j=\begin{pmatrix}C'_j & D'C''_j\end{pmatrix}\in\mathbb{C}^{p\times r_j}, & D=D'D''\in\mathbb{C}^{p\times q}.
\end{array}
\end{equation}
Then $\alpha=(N;A,B,C,D;\mathbb{C}^r=\bigoplus_{j=1}^N\mathbb{C}^{r_j},\mathbb{C}^q,\mathbb{C}^p)$  will be called the \emph{product of GR-nodes} $\alpha'$ and $\alpha''$ and denoted by $\alpha=\alpha'\alpha''$. A straightforward calculation shows that $$T^{\rm nc}_{\alpha}=T^{\rm nc}_{\alpha'}T^{\rm nc}_{\alpha''}.$$

Consider a GR-node
\begin{equation}\label{qq}
\alpha=(N;A,B,C,D;\mathbb{C}^r=\bigoplus_{j=1}^N\mathbb{C}^{r_j},\mathbb{C}^q):=(N;A,B,C,D;\mathbb{C}^r=\bigoplus_{j=1}^N\mathbb{C}^{r_j},\mathbb{C}^q,\mathbb{C}^q)
\end{equation}
 with invertible operator $D$. Then $$\alpha^\times=(N;A^\times,B^\times,C^\times,D^\times;\mathbb{C}^r=\bigoplus_{j=1}^N\mathbb{C}^{r_j},\mathbb{C}^q),$$
 with
\begin{equation}\label{cross}
A^\times=A-BD^{-1}C,\quad B^\times=BD^{-1},\quad C^\times=-D^{-1}C,\quad D^\times=D^{-1},
\end{equation}
will be called the \emph{associated GR-node}, and $A^\times$ the \emph{associated main operator},  of $\alpha$. It is easy to see that, as well as in the one-variable case, $\left(T_\alpha^{\rm nc}\right)^{-1}=T_{\alpha^\times}^{\rm nc}$. Moreover, ${(\alpha^\times)}^\times=\alpha$ (in particular, ${(A^\times)}^\times=A$), and $(\alpha'\alpha'')^\times=\alpha''^\times\alpha'^\times$ up to the natural identification of $\mathbb{C}^{r'_j}\oplus\mathbb{C}^{r''_j}$ with $\mathbb{C}^{r''_j}\oplus\mathbb{C}^{r'_j}$, $j=1,\ldots,N$, which is a similarity transform.
\begin{thm}\label{thm:ass} A GR-node \eqref{qq} with invertible operator $D$ is minimal if and only if its associated GR-node $\alpha^\times$ is minimal.
\end{thm}
\begin{proof}
Let a GR-node $\alpha$ of the form \eqref{qq} with invertible operator $D$ be minimal, and $x\in\ker\mathcal{O}_k^\times$ for some $k\in\{ 1,\ldots,N\}$, where $\mathcal{O}_k^\times$ is the $k$-th observability matrix for the GR-node $\alpha^\times$. Then $x\in\ker(C^\times\flat A^\times)^{wg_k}$ for every $w\in\mathcal{F}_N$. Let us show that $x\in\ker\mathcal{O}_k=\bigcap_{w\in\mathcal{F}_N}\ker(C\flat A)^{wg_k}$, i.e, $x=0$.

For $w=\emptyset$, $C_k^\times x=0$ means $-D^{-1}C_kx=0$ (see \eqref{cross}), which is equivalent to $C_kx=0$.
For $|w|>0,\ w=g_{i_1}\cdots g_{i_{|w|}}$,
\begin{eqnarray*}
(C\flat A)^{wg_k} &=& C_{i_1}A_{i_1i_2}\cdots A_{i_{|w|}k}\\
&=& -DC_{i_1}^\times (A_{i_1i_2}^\times    +B_{i_1}D^{-1}C_{i_2})\cdots (A_{i_{|w|}k}^\times +B_{i_{|w|}}D^{-1}C_k)\\
&=& L_0C_k^\times +\sum_{j=1}^{|w|}L_jC_{i_j}^\times A_{i_ji_{j+1}}^\times\cdots A_{i_{|w|}k}^\times,
\end{eqnarray*}
with some matrices $L_j\in\mathbb{C}^{q\times q},\ j=0,1,\ldots,|w|$. Thus, $x\in\ker(C\flat A)^{wg_k}$ for every $w\in\mathcal{F}_N$, i.e., $x=0$, which means that $\alpha^\times$ is observable.

Since $\alpha$ is controllable if and only if $\alpha^*$ is observable (see Section~\ref{sec:prelim}), and $D^*$ is invertible whenever $D$ is invertible, the same is true for $\alpha^\times$ and $(\alpha^\times)^*=(\alpha^*)^\times$. Thus, the controllability of $\alpha^\times$ follows from the controllability of $\alpha$. Finally, the minimality of $\alpha^\times$ follows from the minimality of $\alpha$. Since $(\alpha^\times)^\times=\alpha$,  the minimality of $\alpha$ follows from the minimality of $\alpha^\times$.
\end{proof}

Suppose that for a GR-node \eqref{qq}, projections $\Pi_k$ on $\mathbb{C}^{r_k}$
are defined such that
$$
A_{kj}\ker\Pi_j\subset \ker \Pi_k,\quad
(A^\times)_{kj}{\rm ran}\,\Pi_j\subset {\rm ran}\,\Pi_k,\quad
k,j=1,\ldots, N.$$ We do not assume that $\Pi_k$ are
orthogonal. We shall call $\Pi_k$ a \emph{$k$-th supporting
projection} for $\alpha$. Clearly, the map $\Pi={\rm
diag}(\Pi_1,\ldots, \Pi_N):\ {\mathbb C}^r\rightarrow{\mathbb C}^r$ satisfies
$$
A\ker\Pi\subset \ker \Pi,\quad\quad A^\times{\rm
ran}\,\Pi\subset {\rm ran}\,\Pi,$$ i.e., it is a \emph{supporting
projection} for the one-variable node $(1;A,B,C,D; {\mathbb
C}^r,{\mathbb C}^q)$ in the sense of \cite{BGK}. If $\Pi$ is a
supporting projection for $\alpha$, then $I_r-\Pi$ is a supporting
projection for $\alpha^\times$.

The following theorem and corollary are analogous to, and are
proved in the same way as Theorem 1.1 and its corollary in
\cite[pp. 7--9]{BGK} (see also \cite[Theorem~2.1]{S}).
\begin{thm}
\label{thm:fact}
Let \eqref{qq} be a GR-node with invertible operator $D$.
Let $\Pi_k$ be a projection on ${\mathbb C}^{r_k}$, and let
$$
A=
\begin{pmatrix} A_{kj}^{(11)}&A_{kj}^{(12)}\\
A_{kj}^{(21)}&A_{kj}^{(22)}\end{pmatrix},\quad
B_j=\begin{pmatrix} B_{j}^{(1)}\\B_{j}^{(2)}\end{pmatrix},\quad
C_k=\begin{pmatrix} C_{k}^{(1)}&C_{k}^{(2)}\end{pmatrix}
$$
be the block matrix representations of the operators $A_{kj}, B_j$ and
$C_k$ with respect to the decompositions ${\mathbb C}^{r_k}=
\ker\Pi_k\dot{+} {\rm ran}\,\Pi_k$, for $k,j\in\left\{1,\ldots, N\right\}$.
Assume that $D=D' D''$, where $D'$ and
$D''$ are invertible operators on ${\mathbb C}^q$, and set
\[
\begin{split}
\alpha'&=(N;A^{(11)}, B^{(1)}(D'')^{-1}, C^{(1)},
D'; \ker\Pi=\bigoplus_{k=1}^N
\ker \Pi_k,{\mathbb C}^q),\\
\alpha''&=(N;A^{(22)}, B^{(2)},(D')^{-1} C^{(2)},
D''; {\rm ran}\,\Pi=\bigoplus_{k=1}^N
{\rm ran}\, \Pi_k,{\mathbb C}^q).
\end{split}
\]
Then $\alpha=\alpha'\alpha''$ (up to a similarity
which maps ${\mathbb C}^{r_k}=\ker\Pi_k\dot{+} {\rm ran}\,\Pi_k$ onto
${\mathbb C}^{{\rm dim}(\ker\Pi_k)}\oplus
{\mathbb C}^{{\rm dim}({\rm ran}\Pi_k)}$ ($k=1,\ldots, N$) such that
$\ker\Pi_k\stackrel{\cdot}{+}\left\{0\right\}$ is mapped onto
${\mathbb C}^{{\rm dim}(\ker\Pi_k)}\oplus\left\{0\right\}$ and
$\left\{0\right\}\stackrel{\cdot}{+}{\rm ran}\Pi_k$ is mapped onto
$\left\{0\right\}\oplus{\mathbb C}^{{\rm dim}({\rm ran}\Pi_k)}$ )  if and only
if $\Pi$ is a supporting projection for $\alpha$.
\end{thm}
\begin{cor}\label{cor:fact}
In the assumptions of Theorem~\ref{thm:fact},
$$T_\alpha^{\rm nc}=F'F'',$$
where
\[
\begin{split}
F'(z)&=D'+C(I_r-\Delta(z)A)^{-1}(I_r-\Pi)\Delta(z)B
(D^{\prime\prime})^{-1},\\
F''(z)&=D^{\prime\prime}+(D^{\prime})^{-1}C\Pi(I_r-\Delta(z)A)^{-1}\Delta(z)B.
\end{split}
\]
\end{cor}
We assume now that the external operator of the GR-node
\eqref{qq} is equal to $D=I_q$ and that we also take
$D^{\prime}=D^{\prime\prime}=I_q$. Then, the GR-nodes
$\alpha'$ and $\alpha^{\prime\prime}$ of Theorem~\ref{thm:fact}
are called \emph{projections of $\alpha$ with respect to the supporting projections
$I_r-\Pi$ and $\Pi$}, respectively, and we use the notations
\[
\begin{split}
\alpha'&={\rm pr}_{I_r-\Pi}(\alpha)=
\left(N;A^{(11)}, B^{(1)}, C^{(1)},
D'; \ker\Pi=\bigoplus_{k=1}^N
\ker \Pi_k,{\mathbb C}^q\right),\\
\alpha^{\prime\prime}&={\rm pr}_\Pi(\alpha)=
\left(N;A^{(22)}, B^{(2)},C^{(2)},
D^{\prime\prime}; {\rm ran\,}\Pi=\bigoplus_{k=1}^N
{\rm ran\,} \Pi_k,{\mathbb C}^q\right).
\end{split}
\]

Let $F',F''$ and $F$ be rational FPSs
with coefficients in ${\mathbb C}^{q\times q}$ such that
\begin{equation}
\label{wagram}
F=F' F^{\prime\prime}.
\end{equation}
The factorization \eqref{wagram} will be said to be \emph{minimal} if whenever
$\alpha'$ and $\alpha^{\prime\prime}$ are minimal GR-realizations of $F'$ and $F^{\prime\prime}$, respectively,
$\alpha'\alpha^{\prime\prime}$ is a minimal GR-realization
of $F$.

In the sequel, we will use the notation
\begin{equation}\label{min}
\alpha=\left(N;A,B,C,D;\mathbb{C}^\gamma=\bigoplus_{k=1}^N\mathbb{C}^{\gamma_k\times\gamma_k},\mathbb{C}^q\right)
\end{equation}
for a minimal GR-realization (i.e., $r_k=\gamma_k$ for $k=1,\ldots,N$) of a rational FPS $F$ in the case when $p=q$.

The following theorem is the multivariable non-commutative version of
\cite[Theorem 4.8]{BGK}. It gives a complete description of all
minimal factorizations in terms of supporting projections.
\begin{thm}
Let $F$ be a rational FPS with a minimal GR-realization \eqref{min}. Then
the following statements hold:
\begin{description}
        \item[(i)] if $\Pi={\rm diag}(\Pi_1,\ldots, \Pi_N)$ is a supporting
projection for $\alpha$, then $F'$ is the transfer function of
${\rm pr}_{I_\gamma-\Pi}(\alpha)$,
$F^{\prime\prime}$ is the transfer function of ${\rm pr}_{\Pi}(\alpha)$,
and $F=F' F^{\prime\prime}$ is a minimal factorization of $F$;
        \item[(ii)]  if $F=F' F^{\prime\prime}$ is a minimal factorization of $F$, then
there exists a uniquely defined supporting projection
$\Pi={\rm diag}(\Pi_1,\ldots, \Pi_N)$ for the GR-node
$\alpha$ such that $F'$ and $F^{\prime\prime}$ are the transfer
functions of ${\rm pr}_{I_\gamma-\Pi}(\alpha)$ and ${\rm pr}_{\Pi}
(\alpha)$, respectively.
\end{description}
\end{thm}
\begin{proof} (i). Let $\Pi$ be a supporting projection for $\alpha$. Then, by
Theorem \ref{thm:fact},
$$\alpha={\rm pr}_{I_\gamma-\Pi}(\alpha){\rm pr}_{\Pi}(\alpha).$$
By the assumption, $\alpha$ is minimal. We now show that
the GR-nodes $\alpha'={\rm pr}_{I_\gamma-\Pi}(\alpha)$ and
$\alpha^{\prime\prime}={\rm pr}_{\Pi}(\alpha)$ are also minimal. To this end, let $x\in {\rm ran}\,\Pi_k$. Then
$$
\left(C^{(2)}\flat A^{(22)}\right)^{w g_k}x=
\left(C\flat A\right)^{w g_k}\Pi_kx=
\left(C\flat A\right)^{w g_k}x.$$
Thus, if ${\mathcal O}_k^{\prime\prime}$ denotes the $k$-th observability
matrix of $\alpha^{\prime\prime}$, then $x\in\ker{\mathcal O}_k^{\prime\prime}$ implies
$x\in\ker{\mathcal O}_k$, and the observability of $\alpha$
implies that $\alpha^{\prime\prime}$ is also observable. Since
$$
\left(A^{(22)}\sharp B^{(2)}\right)^{g_kw^T}=\Pi_k
\left(A\sharp B\right)^{g_kw^T},$$
one has ${\mathcal C}^{\prime\prime}_k=\Pi_k{\mathcal C}_k$, where
${\mathcal C}^{\prime\prime}_k$ is the $k$-th controllability matrix of
$\alpha^{\prime\prime}$. Thus, the controllability of
$\alpha$ implies the controllability of $\alpha^{\prime\prime}$. Hence, we
have proved the minimality of $\alpha^{\prime\prime}$. Note that we have
used that $\ker\Pi={\rm ran}\,(I_\gamma-\Pi)$ is $A$-invariant. Since ${\rm ran}\,\Pi=\ker(I_\gamma-\Pi)$
is $A^\times$-invariant, by Theorem \ref{thm:ass} $\alpha^\times$ is
minimal. Using
$$\alpha^\times=(\alpha'\alpha^{\prime\prime})^\times=
(\alpha^{\prime\prime})^\times(\alpha^{\prime})^\times,$$
we prove the minimality of $(\alpha')^\times$ in the same way as that
of $\alpha''$. Applying once again Theorem \ref{thm:ass}, we obtain the
minimality of $\alpha'$. The dimensions of the state spaces
of the minimal GR-nodes $\alpha', \alpha''$ and $\alpha$
are related by
$$\gamma_k=\gamma_k'+\gamma_k^{\prime\prime},\qquad k=1,\ldots, N.$$
Therefore, given any minimal GR-realizations
$\beta'$ and $\beta^{\prime\prime}$ of $F'$ and
$F^{\prime\prime}$, respectively, the same equalities hold for the state space dimensions of
$\beta'$, $\beta''$ and
$\beta$. Thus, $\beta'\beta^{\prime\prime}$ is a minimal
GR-node, and the factorization $F=F' F^{\prime\prime}$
is minimal.

(ii). Assume that the factorization $F=F' F^{\prime\prime}$
is minimal. Let $\beta'$ and $\beta^{\prime\prime}$ be minimal
GR-realizations of $F'$ and $F^{\prime\prime}$ with $k$-th state
space dimensions equal to $\gamma'_k$ and $\gamma^{\prime\prime}_k$,
respectively ($k=1,\ldots, N$). Then $\beta'\beta^{\prime\prime}$ is a minimal
GR-realization of $F$ and its $k$-th  state space dimension
is equal to $\gamma_k=\gamma_k'+\gamma_k^{\prime\prime}$
($k=1,\ldots, N$). Hence $\beta'\beta^{\prime\prime}$ is similar to
$\alpha$. We denote the corresponding GR-node similarity by
$T={\rm diag}(T_1,\ldots, T_N)$, where
$$
T_k:\, {\mathbb C}^{\gamma^{\prime}}\oplus
{\mathbb C}^{\gamma^{\prime\prime}}\rightarrow
{\mathbb C}^{\gamma},\quad k=1,\ldots N,$$
is the canonical isomorphism.
Let $\Pi_k$ be the projection of ${\mathbb C}^{\gamma_k}$ along
$T_k{\mathbb C}^{\gamma_k'}$ onto $T_k
{\mathbb C}^{\gamma_k^{\prime\prime}}$, $k=1,\ldots, N$, and set
$\Pi={\rm diag}(\Pi_1,\ldots, \Pi_k)$. Then $\Pi$ is a supporting projection
for $\alpha$. Moreover ${\rm pr}_{I_\gamma-\Pi}(\alpha)$ is similar to
$\beta'$, and ${\rm pr}_\Pi(\alpha)$ is similar to $\beta^{\prime\prime}$. The
uniqueness of $\Pi$ is proved in the same way as in \cite[Theorem 4.8]{BGK}.
The uniqueness of the GR-node similarity follows
from Theorem \ref{thm:sim}.
\end{proof}

\section{Matrix-$J$-unitary formal power series: A multivariable non-commutative analogue of the line case}\label{sec:line}
In this section we study a multivariable non-commutative analogue of rational $q\times q$ matrix-valued functions which are $J$-unitary on the imaginary line $i\mathbb{R}$ of the complex plane $\mathbb{C}$.

\subsection{Minimal Givone--Roesser realizations and the Lyapunov equation}
Denote by $\mathbb{H}^{n\times n}$ the set of Hermitian $n\times n$ matrices. Then  $\left(i\mathbb{H}^{n\times n}\right)^N$ will denote the set of $N$-tuples of skew-Hermitian matrices. In our paper, the set $$\mathcal{J}_N=\coprod_{n\in\mathbb{N}}\left(i\mathbb{H}^{n\times n}\right)^N,$$
where ``$\coprod$" stands for a disjoint union, will be a counterpart of the imaginary line $i\mathbb{R}$.

Let $J\in\mathbb{C}^{q\times q}$ be a signature matrix. We will call a rational FPS $F\in\mathbb{C}^{q\times q}\left\langle \left\langle z_1,\ldots,z_N\right\rangle\right\rangle_{\rm rat}$
\emph{matrix-$J$-unitary on $\mathcal{J}_N$} if for every $n\in\mathbb{N}$,
\begin{equation}\label{j-un}
F(Z)(J\otimes I_n)F(Z)^*=J\otimes I_n
\end{equation}
at all points $Z\in\left(i\mathbb{H}^{n\times n}\right)^N$ where it is defined. For a fixed $n\in\mathbb{N}$, $F(Z)$ as a function of matrix entries is rational and holomorphic on some open neighbourhood $\Gamma_n(\varepsilon)$ of $Z=0$, e.g., of the form \eqref{gamma}, and $\Gamma_n(\varepsilon)\cap \left(i\mathbb{H}^{n\times n}\right)^N$ is a uniqueness set in
$(\mathbb{C}^{n\times n})^N$ (see \cite{Sh} for the uniqueness theorem in several complex variables). Thus, \eqref{j-un} implies that
\begin{equation}\label{j-un-z}
F(Z)(J\otimes I_n)F(-Z^*)^*=J\otimes I_n
\end{equation}
at all points $Z\in(\mathbb{C}^{n\times n})^N$ where $F(Z)$ is holomorphic and invertible (the set of such points is open and dense, since $\det F(Z)\not\equiv 0$).

The following theorem is a counterpart of Theorem~2.1 in \cite{AG1}.
\begin{thm}\label{thm:Lyap}
Let $F$ be a rational FPS with a minimal GR-realization \eqref{min}. Then $F$ is matrix-$J$-unitary on $\mathcal{J}_N$ if and only if the following conditions are fulfilled:

a) $D$ is $J$-unitary, i.e., $DJD^*=J$;

b) there exists an invertible Hermitian solution $H={\rm diag}(H_1,\ldots,H_N)$, with $H_k\in\mathbb{C}^{\gamma_k\times\gamma_k},\ k=1,\ldots,N$, of the Lyapunov equation
\begin{equation}\label{L}
A^*H+HA=-C^*JC,
\end{equation}
and
\begin{equation}\label{b}
B=-H^{-1}C^*JD.
\end{equation}
The property b) is equivalent to

b') there exists an invertible Hermitian matrix $H={\rm diag}(H_1,\ldots,H_N)$, with $H_k\in\mathbb{C}^{\gamma_k\times\gamma_k},\ k=1,\ldots,N$, such that
\begin{equation}\label{L'}
H^{-1}A^*+AH^{-1}=-BJB^*,
\end{equation}
and
\begin{equation}\label{c}
C=-DJB^*H.
\end{equation}
\end{thm}
\begin{proof}
Let $F$ be matrix-$J$-unitary. Then $F$ is holomorphic at the point $Z=0$ in $\mathbb{C}^N$, hence $D=F(0)$ is $J$-unitary (in particular, invertible). Equality \eqref{j-un-z} may be rewritten as
\begin{equation}\label{ttf}
F(Z)^{-1}=(J\otimes I_n)F(-Z^*)^*(J\otimes I_n).
\end{equation}
Since \eqref{ttf} holds for all  $n\in\mathbb{N}$, it follows from Corollary~\ref{cor:vanish-fps} that the FPSs corresponding to the left and the right sides of equality \eqref{ttf} coincide. Due to Theorem~\ref{thm:ass}, $\alpha^\times=(N;A^\times,B^\times,C^\times,D^\times;\mathbb{C}^\gamma=\bigoplus_{k=1}^N\mathbb{C}^{\gamma_k},\mathbb{C}^q)$ with $A^\times,B^\times,C^\times,D^\times$ given by \eqref{cross} is a minimal GR-realization of $F^{-1}$. Due to \eqref{ttf}, another minimal GR-realization of $F^{-1}$ is $\tilde{\alpha}=(N;\tilde{A},\tilde{B},\tilde{C},\tilde{D};\mathbb{C}^\gamma=\bigoplus_{k=1}^N\mathbb{C}^{\gamma_k},\mathbb{C}^q)$, where
$$\tilde{A}=-A^*,\quad \tilde{B}=C^*J,\quad \tilde{C}=-JB^*,\quad \tilde{D}=JD^*J.$$
By Theorem~\ref{thm:sim}, there exists unique similarity transform $T={\rm diag}(T_1,\ldots,T_N)$ which relates $\alpha^\times$ and $\tilde{\alpha}$, where $T_k\in\mathbb{C}^{\gamma_k\times\gamma_k}$ are invertible for $k=1,\ldots,N$, and
\begin{equation}\label{eq-sim}
T(A-BD^{-1}C) =-A^*T,\quad TBD^{-1}=C^*J,\quad D^{-1}C=JB^*T.
\end{equation}
Note that the relation $D^{-1}=JD^*J$, which means $J$-unitarity of $D$, has been already established above. It is easy to check that relations \eqref{eq-sim} are also valid for $T^*$ in the place of $T$. Hence, by the uniqueness of similarity matrix, $T=T^*$. Setting $H=-T$, we obtain from \eqref{eq-sim} the equalities \eqref{L} and \eqref{b}, as well as \eqref{L'} and \eqref{c}, by a straightforward calculation.

Let us prove now a slightly more general statement than the converse. Let $\alpha$ be a (not necessarily minimal) GR-realization of $F$ of the form \eqref{qq}, where $D$ is $J$-unitary, and let $H={\rm diag}(H_1,\ldots,H_N)$ with $H_k\in\mathbb{C}^{r_k\times r_k},\ k=1,\ldots,N$,  be an Hermitian invertible matrix satisfying \eqref{L} and \eqref{b}. Then in the same way as in \cite[Theorem~2.1]{AG1} for the one-variable case, we obtain for $Z,Z'\in\mathbb{C}^{n\times n}$:
\begin{eqnarray}
\lefteqn{F(Z)(J\otimes I_n)F(Z')^*=J\otimes I_n-(C\otimes I_n)\left(I_r\otimes I_n-\Delta(Z)(A\otimes I_n)\right)^{-1}}\nonumber\\
&\times  \Delta(Z+Z'^*)(H^{-1}\otimes I_n)\left(I_r\otimes I_n-(A^*\otimes I_n)\Delta(Z'^*)\right)^{-1}(C^*\otimes I_n)\label{L-id}
\end{eqnarray}
(note that $\Delta(Z)$ commutes with $H^{-1}\otimes I_n$). It follows from \eqref{L-id} that $F(Z)$ is $(J\otimes I_n)$-unitary on $(i\mathbb{H}^{n\times n})^N$ at all points $Z$ where it is defined. Since $n\in\mathbb{N}$ is arbitrary, $F$ is matrix-$J$-unitary on $\mathcal{J}_N$. Clearly, conditions a) and b') also imply the matrix-$J$-unitarity of $F$ on $\mathcal{J}_N$.
\end{proof}
Let us make some remarks. First, it follows from the proof of Theorem~\ref{thm:Lyap} that the  structured solution $H={\rm diag}(H_1,\ldots,H_N)$ of the Lyapunov equation \eqref{L} is uniquely determined by a given minimal GR-realization of $F$. The matrix $H={\rm diag}(H_1,\ldots,H_N)$ is called the \emph{associated structured Hermitian matrix}
(associated with this minimal GR-realization of $F$). The matrix $H_k$ will be called the \emph{$k$-th component of the associated Hermitian matrix} ($k=1,\ldots,N$). The explicit formulas for $H_k$ follow from \eqref{t}:
\begin{eqnarray*}
\lefteqn{H_k = -\left[{\rm col}_{|w|\leq qr-1}\left((JB^*)\flat (-A^*)\right)^{wg_k}\right]^+{\rm col}_{|w|\leq qr-1}\left((D^{-1}C)\flat A^\times\right)^{wg_k}}\\
&= -{\rm row}_{|w|\leq qr-1}\left((-A^*)\sharp (C^*J)\right)^{g_kw^T}\left[{\rm row}_{|w|\leq qr-1}\left(A^\times\sharp (BD^{-1})\right)^{g_kw^T}\right]^\dag.
\end{eqnarray*}

Second, let $\alpha$ be a (not necessarily minimal) GR-realization of $F$ of the form \eqref{qq}, where $D$ is $J$-unitary, and let $H={\rm diag}(H_1,\ldots,H_N)$ with $H_k\in\mathbb{C}^{r_k\times r_k},\ k=1,\ldots,N$,  be an Hermitian, not necessarily invertible, matrix satisfying \eqref{L} and \eqref{c}. Then in the same way as in \cite[Theorem~2.1]{AG1} for the one-variable case, we obtain for $Z,Z'\in\mathbb{C}^{n\times n}$:
\begin{eqnarray}
F(Z')^*(J\otimes I_n)F(Z)=J\otimes I_n-(B^*\otimes I_n)\left(I_r\otimes I_n-\Delta(Z'^*)(A^*\otimes I_n)\right)^{-1}\nonumber\\
 \times (H\otimes I_n)\Delta(Z'^*+Z)\left(I_r\otimes I_n-(A\otimes I_n)\Delta(Z)\right)^{-1}(B\otimes I_n)\hspace{1cm} \label{L-id'}
\end{eqnarray}
(note that $\Delta(Z)$ commutes with $H\otimes I_n$). It follows from \eqref{L-id'} that $F(Z)$ is $(J\otimes I_n)$-unitary on $(i\mathbb{H}^{n\times n})^N$ at all points $Z$ where it is defined. Since $n\in\mathbb{N}$ is arbitrary, $F$ is matrix-$J$-unitary on $\mathcal{J}_N$.

Third, if  $\alpha$ is a (not necessarily minimal) GR-realization of $F$ of the form \eqref{qq}, where $D$ is $J$-unitary, and equalities \eqref{L'} and \eqref{c} are valid with $H^{-1}$ replaced by some, possibly not invertible, Hermitian matrix  $Y={\rm diag}(Y_1,\ldots,Y_N)$ with $Y_k\in\mathbb{C}^{r_k\times r_k},\ k=1,\ldots,N$,
then $F$ is matrix-$J$-unitary on $\mathcal{J}_N$. This follows from the fact that \eqref{L-id} is valid with $H^{-1}$ replaced by $Y$.
\begin{thm}\label{thm:inv}
Let $(C,A)$ be an observable pair of matrices $C\in\mathbb{C}^{q\times r},A\in\mathbb{C}^{r\times r}$ in the sense that $\mathbb{C}^r=\bigoplus_{k=1}^N\mathbb{C}^{r_k}$ and $\mathcal{O}_k$ has full column rank for each $k\in\{ 1,\ldots,N\}$, and let $J\in\mathbb{C}^{q\times q}$ be a signature matrix. Then there exists a  matrix-$J$-unitary on $\mathcal{J}_N$ rational FPS $F$ with a minimal GR-realization $\alpha=(N;A,B,C,D;\mathbb{C}^r=\bigoplus_{k=1}^N\mathbb{C}^{r_k},\mathbb{C}^q)$ if and only if the Lyapunov equation \eqref{L} has a structured solution $H={\rm diag}(H_1,\ldots,H_N)$ which is both Hermitian and invertible. If such a solution $H$ exists, possible choices of $D$ and $B$ are
\begin{equation}\label{db}
D_0=I_q,\quad B_0=-H^{-1}C^*J.
\end{equation}
Finally, for a given such $H$, all other choices of $D$ and $B$ differ from $D_0$ and $B_0$ by a right multiplicative $J$-unitary constant matrix.
\end{thm}
\begin{proof}
Let $H={\rm diag}(H_1,\ldots,H_N)$ be a structured solution of the Lyapunov equation \eqref{L} which is both Hermitian and invertible. We first check that the pair $(A,-H^{-1}C^*J)$ is controllable, or equivalently, that the pair $(-JCH^{-1},A^*)$ is observable. Using the Lyapunov equation \eqref{L}, one can see that for any $k\in\{ 1,\ldots,N\}$ and $w=g_{i_1}\cdots g_{i_{|w|}}\in\mathcal{F}_N$ there exist matrices $K_0,\ldots,K_{|w|-1}$ such that
\begin{eqnarray*}
(C\flat A)^{wg_k} &=& (-1)^{|w|-1}J((-JCH^{-1})\flat A^*)^{wg_k}H_k\\
&+& K_0J(-JC_{i_2}H_{i_2}^{-1}(A^*)_{i_2i_3}\cdots(A^*)_{i_{|w|}k})H_k+ \cdots\\
& +& K_{|w|-2}J(-JC_{i_{|w|}}(A^*)_{i_{|w|}k})H_k+K_{|w|-1}J(-JC_kH_k^{-1})H_k.
\end{eqnarray*}
Thus, if $x\in\ker ((-JCH^{-1})\flat A^*)^{wg_k}$ for all of $w\in\mathcal{F}_N$ then $H_k^{-1}x\in\ker\mathcal{O}_k$, and the observability of the pair $(C,A)$ implies that $x=0$. Therefore, the pair $(-JCH^{-1},A^*)$ is observable, and the pair $(A,-H^{-1}C^*J)$ is controllable. By Theorem~\ref{thm:Lyap} we obtain that
\begin{equation}\label{f0}
F_0(z)=I_q-C(I_r-\Delta(z)A)^{-1}\Delta(z)H^{-1}C^*J
\end{equation}
is a matrix-$J$-unitary on $\mathcal{J}_N$ rational FPS, which has a minimal GR-realization $\alpha_0=(N:A,-H^{-1}C^*J,C,I_q;\mathbb{C}^r=\bigoplus_{k=1}^N\mathbb{C}^{r_k},\mathbb{C}^q)$ with the associated structured Hermitian matrix $H$.

Conversely, let $\alpha=(N;A,B,C,D;\mathbb{C}^r=\bigoplus_{k=1}^N\mathbb{C}^{r_k},\mathbb{C}^q)$ be a minimal GR-node. Then by Theorem~\ref{thm:Lyap} there exists an Hermitian and invertible matrix $H={\rm diag}(H_1,\ldots,H_N)$ which solves \eqref{L}.

Given $H={\rm diag}(H_1,\ldots,H_N)$, let $B, D$ be any solution of the inverse problem, i.e., $\alpha=(N;A,B,C,D;\mathbb{C}^r=\bigoplus_{k=1}^N\mathbb{C}^{r_k},\mathbb{C}^q)$ is a minimal GR-node with the associated structured Hermitian matrix $H$. Then for $F_0=T_{\alpha_0}^{\rm nc}$ and $F=T_\alpha^{\rm nc}$ we obtain from \eqref{L-id} that
$$F(Z)(J\otimes I_n)F(Z')^*=F_0(Z)(J\otimes I_n)F_0(Z')^*$$
 for any $n\in\mathbb{N}$, at all points $Z,Z'\in(\mathbb{C}^{n\times n})^N$ where both $F$ and $F_0$ are defined. By the uniqueness theorem in several complex variables (matrix entries for $Z_k$'s and ${Z'}_k^*$'s, $k=1,\ldots,N$), we obtain that $F(Z)$ and $F_0(Z)$ differ by a right multiplicative $(J\otimes I_n)$-unitary constant, which clearly has to be $D\otimes I_n$, i.e.,
$$F(Z)=F_0(Z)(D\otimes I_n).$$
Since $n\in\mathbb{N}$ is arbitrary, by Corollary~\ref{cor:vanish-fps} we obtain
$$F(z)=F_0(z)D.$$
Equating the coefficients of these two FPSs, we easily deduce using the observability of the pair $(C,A)$  that $B=-H^{-1}C^*JD$.
\end{proof}
The following dual theorem is proved analogously.
\begin{thm}\label{thm:inv'}
Let $(A,B)$ be a controllable pair of matrices $A\in\mathbb{C}^{r\times r},B\in\mathbb{C}^{r\times q}$ in the sense that $\mathbb{C}^r=\bigoplus_{k=1}^N\mathbb{C}^{r_k}$ and $\mathcal{C}_k$ has full row rank for each $k\in\{ 1,\ldots,N\}$, and let $J\in\mathbb{C}^{q\times q}$ be a signature matrix. Then there exists a  matrix-$J$-unitary on $\mathcal{J}_N$ rational FPS $F$ with a minimal GR-realization $\alpha=(N;A,B,C,D;\mathbb{C}^r=\bigoplus_{k=1}^N\mathbb{C}^{r_k},\mathbb{C}^q)$ if and only if the Lyapunov equation
$$GA^*+AG=-BJB^*$$ has a structured solution $G={\rm diag}(G_1,\ldots,G_N)$ which is both Hermitian and invertible. If such a solution $G$ exists, possible choices of $D$ and $C$ are
\begin{equation}\label{dc}
D_0=I_q,\quad C_0=-JB^*G^{-1}.
\end{equation}
Finally, for a given such $G$, all other choices of $D$ and $C$ differ from $D_0$ and $C_0$ by a left multiplicative $J$-unitary constant matrix.
\end{thm}
\begin{thm}\label{thm:c=o}
Let $F$ be a matrix-$J$-unitary on $\mathcal{J}_N$ rational FPS, and $\alpha$ be its GR-realization. Let $H={\rm diag}(H_1,\ldots,H_N)$ with $H_k\in\mathbb{C}^{r_k\times r_k},\ k=1,\ldots,N$,  be an Hermitian invertible matrix satisfying \eqref{L} and \eqref{b}, or equivalently, \eqref{L'} and \eqref{c}. Then $\alpha$ is observable if and only if $\alpha$ is controllable.
\end{thm}
\begin{proof}
Suppose that $\alpha$ is observable. Since by Theorem~\ref{thm:Lyap} $D=F_\emptyset$ is $J$-unitary, by Theorem~\ref{thm:inv} $\alpha$ is a minimal GR-node. In particular, $\alpha$ is controllable.

Suppose that $\alpha$ is controllable. Then by Theorem~\ref{thm:inv'} $\alpha$ is minimal, and in particular, observable.
\end{proof}

\subsection{The associated structured Hermitian matrix}\label{sec:Herm-line}
\begin{lem}
\label{lem:h}
Let $F$ be a matrix-$J$-unitary on ${\mathcal J}_N$ rational
FPS, and let
$\alpha^{(i)}=(N;A^{(i)},B^{(i)},C^{(i)},D;{\mathbb C}^\gamma=
\bigoplus_{k=1}^N{\mathbb C}^{\gamma_k},{\mathbb C}^q)$ be
minimal GR-realizations of $F$, with the associated structured Hermitian
matrices $H^{(i)}={\rm diag}(H_1^{(i)},\ldots, H_N^{(i)})$, $i=1,2$. Then $\alpha^{(1)}$ and $\alpha^{(2)}$ are similar, i.e., \eqref{sim} holds with
 a uniquely defined invertible matrix $T={\rm diag}(T_1,\ldots, T_N)$,
and
\begin{equation}
\label{h-sim}
H_k^{(1)}=T_k^*H_k^{(2)}T_k,\qquad k=1,\ldots, N.
\end{equation}
In particular, the matrices $H_k^{(1)}$ and $H_k^{(2)}$ have the same
signature.
\end{lem}
The proof is easy and analogous to the proof of Lemma 2.1 in \cite{AG1}.
\begin{rem}
\label{rem:un}
The similarity matrix $T={\rm diag}(T_1,\ldots, T_N)$ is a unitary mapping
from ${\mathbb C}^\gamma=\bigoplus_{k=1}^N{\mathbb C}^{\gamma_k}$ endowed
with the inner product $[\,\cdot\, ,\, \cdot\,]_{H^{(1)}}$ onto
${\mathbb C}^\gamma=\bigoplus_{k=1}^N{\mathbb C}^{\gamma_k}$ endowed
with the inner product $[\,\cdot\, ,\, \cdot\,]_{H^{(2)}}$, where
$$
[x,y]_{H^{(i)}}=\langle H^{(i)}x,y\rangle_{{\mathbb C}^\gamma},\qquad x,y\in{\mathbb C}^\gamma,\ i=1,2,
$$
that is,
$$
[x,y]_{H^{(i)}}=\sum_{k=1}^N[x_k,y_k]_{H_k^{(i)}},\quad i=1,2,$$
where $x_k,y_k\in{\mathbb C}^{\gamma_k}$,  $x={\rm col}_{k=1,\ldots ,N}(x_k)$,
 $y={\rm col}_{k=1,\ldots ,N}(y_k)$, and
 $$
[x_k,y_k]_{H^{(i)}_k}=\langle H^{(i)}_kx_k,y_k\rangle_{{\mathbb C}^{\gamma_k}},\quad k=1,\ldots,N,\ i=1,2.
$$
\end{rem}
Recall the following definition \cite{KL}. Let $K_{w,w'}$ be a $\mathbb{C}^{q\times q}$-valued function defined for $w$ and $w'$ in some set $E$ and such that $(K_{w,w'})^*=K_{w',w}$. Then $K_{w,w'}$ is called a \emph{kernel with $\kappa$ negative squares} if for any $m\in\mathbb{N}$, any points $w_1,\ldots,w_m$ in $E$, and any vectors $c_1,\ldots,c_m$ in $\mathbb{C}^q$ the matrix $(c_j^*K_{w_j,w_i}c_i)_{i,j=1,\ldots,m}\in\mathbb{H}^{m\times m}$ has at most $\kappa$ negative eigenvalues, and has exactly $\kappa$ negative eigenvalues for some choice of $m,w_1,\ldots,w_m,c_1,\ldots,c_m$. We will use this definition to give a characterization of the number of negative eigenvalues of the $k$-th component  $H_k,\ k=1,\ldots,N$, of the associated structured Hermitian matrix $H$.
\begin{thm}\label{thm:neg}
Let $F$ be a matrix-$J$-unitary on $\mathcal{J}_N$ rational FPS, and let $\alpha$ be its minimal GR-realization of the form \eqref{min}, with the associated structured Hermitian matrix $H={\rm diag}(H_1,\ldots,H_N)$. Then for $k=1,\ldots,N$ the number of negative eigenvalues of the matrix $H_k$ is equal to the number of negative squares of each of the kernels
\begin{eqnarray}
K_{w,w'}^{F,k} &=& (C\flat A)^{wg_k}H_k^{-1}(A^*\sharp C^*)^{g_kw'^T},\quad w,w'\in\mathcal{F}_N,\label{kerns}\\
K_{w,w'}^{F^*,k} &=& (B^*\flat A^*)^{wg_k}H_k(A\sharp B)^{g_kw'^T},\quad w,w'\in\mathcal{F}_N,\label{kerns*}
\end{eqnarray}
For $k=1,\ldots,N$, denote by $\mathcal{K}_k(F)$ (resp., $\mathcal{K}_k(F^*)$) the linear span of the functions $w\mapsto K_{w,w'}^{F,k}c$ (resp., $w\mapsto K_{w,w'}^{F^*,k}c$) where $w'\in\mathcal{F}_N$ and $c\in\mathbb{C}^q$. Then
$$\dim\,\mathcal{K}_k(F)=\dim\,\mathcal{K}_k(F^*)=\gamma_k.$$
\end{thm}
\begin{proof}
Let $m\in\mathbb{N},\ w_1,\ldots,w_m\in\mathcal{F}_N$, and $c_1,\ldots,c_m\in\mathbb{C}^q$. Then the matrix equality
$$(c_j^*K_{w_j,w_i}^{F,k}c_i)_{i,j=1,\ldots,m}=X^*H_k^{-1}X,$$
with
$$X={\rm row}_{1\leq i\leq m}\left((A^*\sharp C^*)^{g_kw_i^T}c_i\right),$$
implies that the kernel $K_{w,w'}^{F,k}$ has at most $\kappa_k$ negative squares, where $\kappa_k$ denotes the number of negative eigenvalues of $H_k$. The pair $(C,A)$ is observable, hence we can choose a basis of $\mathbb{C}^q$ of the form $x_i=(A^*\sharp C^*)^{g_kw_i^T}c_i,\ i=1,\ldots,q$. Since the matrix $X={\rm row}_{i=1,\ldots,q}(x_i)$ is non-degenerate, and therefore the matrix $X^*H_k^{-1}X$ has exactly $\kappa_k$ negative eigenvalues, the
kernel $K_{w,w'}^{F,k}$ has  $\kappa_k$ negative squares. Analogously,
from the controllability of the pair $(A,B)$ one can obtain that the kernel $\mathcal{K}_k(F^*)$ has  $\kappa_k$ negative squares.

Since $\mathcal{K}_k(F)$ is the span of functions (of variable $w\in\mathcal{F}_N$) of the form $(C\flat A)^{wg_k}y,\ y\in\mathbb{C}^{\gamma_k}$, it follows that $\dim\,\mathcal{K}_k(F)\leq\gamma_k$. From the observability of the pair $(C,A)$ we obtain that $(C\flat A)^{wg_k}y\equiv 0$ implies $y=0$, thus $\dim\,\mathcal{K}_k(F)=\gamma_k$. In the same way we obtain that the controllability of the pair $(A,B)$ implies that $\dim\,\mathcal{K}_k(F^*)=\gamma_k$.
\end{proof}
We will denote by $\nu_k(F)$ the number of negative squares of either the
kernel $K_{w,w'}^{F,k}$ or the kernel $K_{w,w'}^{F^*,k}$ defined by
\eqref{kerns} and \eqref{kerns*}, respectively.
\begin{thm}\label{thm:factor}
Let $F^{(i)}$ be matrix-$J$-unitary on $\mathcal{J}_N$ rational FPSs, with minimal GR-realizations $\alpha^{(i)}=(N;A^{(i)},B^{(i)},C^{(i)},D^{(i)};{\mathbb C}^{\gamma^{(i)}}=
\bigoplus_{k=1}^N{\mathbb C}^{\gamma_k^{(i)}},{\mathbb C}^q)$ and the associated structured Hermitian
matrices $H^{(i)}={\rm diag}(H_1^{(i)},\ldots, H_N^{(i)})$, respectively, $i=1,2$. Suppose that the product $\alpha=\alpha^{(1)}\alpha^{(2)}$ is a minimal GR-node. Then the matrix
$H={\rm diag}(H_1,\ldots,H_N)$, with
\begin{equation}\label{hk}
H_k=\begin{pmatrix}
H_k^{(1)} & 0\\
0 & H_k^{(2)}
\end{pmatrix}\in\mathbb{C}^{(\gamma_k^{(1)}+\gamma_k^{(2)})\times (\gamma_k^{(1)}+\gamma_k^{(2)})},\quad k= 1,\ldots,N,
\end{equation}
is the associated structured Hermitian matrix for $\alpha=\alpha^{(1)}\alpha^{(2)}$.
\end{thm}
\begin{proof}
It suffices to check that \eqref{L} and \eqref{b} hold for the matrices $A,B,C,D$ defined as in \eqref{prod}, and $H={\rm diag}(H_1,\ldots,H_N)$ where $H_k,\ k=1,\ldots,N$, are defined in \eqref{hk}. This is an easy computation which is omitted.
\end{proof}
\begin{cor}
Let $F_1$ and $F_2$ be matrix-$J$-unitary on ${\mathcal J}_N$ rational
FPSs, and suppose that the factorization $F=F_1F_2$ is minimal.
Then
\[
\nu_k(F_1F_2)=\nu_k(F_1)+\nu_k(F_2),\quad k=1,\ldots, N.
\]
\end{cor}
\subsection{Minimal matrix-$J$-unitary factorizations}
In this subsection we consider minimal factorizations of rational formal power
series which are matrix-$J$-unitary on ${\mathcal J}_N$ into factors both
of which are also matrix-$J$-unitary on ${\mathcal J}_N$. Such factorizations
will be called \emph{minimal matrix-$J$-unitary factorizations}.

Let $H\in{\mathbb C}^{r\times r}$ be an invertible Hermitian matrix. We denote
by $[\,\cdot\, ,\,\cdot\,]_H$ the Hermitian sesquilinear form
$$[x,y]_H=\langle Hx, y\rangle$$
where $\langle \,\cdot\, ,\,\cdot\,\rangle$ denotes the standard inner product
of ${\mathbb C}^r$. Two vectors $x$ and $y$ in ${\mathbb C}^r$ are called
$H$-orthogonal if $[x,y]_H=0$. For any subspace $M\subset {\mathbb C}^r$
denote
$$
M^{[\perp]}=\left\{y\in{\mathbb C}^r:\,\, \langle y,m\rangle_H=0\quad
\forall m\in M\right\}.$$
The subspace $M\subset {\mathbb C}^r$ is called \emph{non-degenerate} if
$M\cap M^{[\perp]}=\left\{0\right\}$. In this case,
$$
M[\stackrel{\cdot}{+}]M^{[\perp]}={\mathbb C}^r$$
where $[\stackrel{\cdot}{+}]$ denotes the \emph{$H$-orthogonal direct sum}.

In the case when $H={\rm diag}(H_1,\ldots,H_N)$ is the  structured Hermitian matrix associated with a given minimal GR-realization of a matrix-$J$-unitary on $\mathcal{J}_N$ rational FPS $F$, we will call
$[\,\cdot\, ,\,\cdot\,]_H$ the \emph{associated inner product} (associated with the given minimal GR-realization of $F$). In more details,
$$[x,y]_H=\sum_{k=1}^N[x_k,y_k]_{H_k},$$
where $x_k,y_k\in{\mathbb C}^{\gamma_k}$ and $x={\rm col}_{k=1,\ldots, N}(x_k)$,
 $y={\rm col}_{k=1,\ldots ,N}(y_k)$, and
 $$[x_k,y_k]_{H_k}=\left\langle H_kx_k,y_k\right\rangle_{\mathbb{C}^{\gamma_k}},\quad k=1,\ldots,N.$$

The following theorem (as well as its proof) is analogous to its one-variable counterpart, Theorem 2.6 from \cite{AG1} (see also \cite[Chapter~II]{S}).
\begin{thm}\label{thm:min-fact}
Let $F$ be a matrix-$J$-unitary on $\mathcal{J}_N$ rational FPS, and let $\alpha$ be its minimal GR-realization of the form \eqref{min}, with the associated structured Hermitian matrix $H={\rm diag}(H_1,\ldots,H_N)$. Let $\mathcal{M}=\bigoplus_{k=1}^N\mathcal{M}_k$ be an $A$-invariant subspace such that $\mathcal{M}_k\subset\mathbb{C}^{\gamma_k},\ k=1,\ldots,N$, and $\mathcal{M}$ is non-degenerate in the associated inner product $[\,\cdot\, ,\,\cdot\,]_H$. Let
$\Pi={\rm diag}(\Pi_1,\ldots,\Pi_N)$ be the projection defined by
$$\ker\Pi=\mathcal{M},\quad {\rm ran}\,\Pi=\mathcal{M}^{[\perp ]},$$
or in more details,
$$\ker\Pi_k=\mathcal{M}_k,\quad {\rm ran}\,\Pi_k=\mathcal{M}_k^{[\perp ]},\quad k=1,\ldots,N.$$
Let $D=D_1D_2$ be a factorization of $D$ into two $J$-unitary factors. Then
the factorization $F=F_1F_2$ where
\begin{eqnarray*}
F_1(z)&=&D_1+C(I_\gamma-\Delta(z)A)^{-1}\Delta(z)(I_\gamma-\Pi)
BD_2^{-1},\\
F_2(z)&=&D_2+D_1^{-1}C\Pi
(I_\gamma-\Delta(z)A)^{-1}\Delta(z)B,
\end{eqnarray*}
is a minimal matrix-$J$-unitary factorization of $F$.

Conversely, any minimal matrix-$J$-unitary factorization of $F$
can be obtained in such a way. For a fixed $J$-unitary decomposition $D=D_1D_2$, the
correspondence between minimal matrix-$J$-unitary factorizations of $F$
and non-degenerate $A$-invariant subspaces of the form $\mathcal{M}=\bigoplus_{k=1}^N
\mathcal{M}_k$, where $\mathcal{M}_k\subset{\mathbb C}^{\gamma_k}$ for $k=1,\ldots, N$, is
one-to-one.
\end{thm}
\begin{rem} We omit here the proof, which can be easily restored, with making use of Theorem \ref{thm:fact} and
Corollary \ref{cor:fact}.
\end{rem}
\begin{rem}
Minimal matrix-$J$-unitary factorizations do not always exist, even for
$N=1$. Examples of $J$-unitary on $i\mathbb{R}$ rational functions which
have non-trivial minimal factorizations
but lack minimal $J$-unitary factorizations can be found in \cite{AD1} and
\cite{AG1}.
\end{rem}
\subsection{Matrix-unitary rational formal power series}
In this subsection we specialize some of the preceding results to the case
$J=I_q$. We call the corresponding rational formal power series
\emph{matrix-unitary} on ${\mathcal J}_N$.
\begin{thm}
\label{thm:L-un}
Let $F$ be a rational FPS and
$\alpha$ be its minimal GR-realization of the form \eqref{min}. Then $F$ is matrix-unitary on ${\mathcal J}_N$ if and only if the following conditions are fulfilled:

a) $D$ is a unitary matrix, i.e., $DD^*=I_q$;

b) there exists an Hermitian solution $H={\rm diag}(H_1,\ldots,H_N)$, with $H_k\in\mathbb{C}^{\gamma_k\times\gamma_k},\ k=1,\ldots,N$, of the Lyapunov equation
\begin{equation}\label{L-un}
A^*H+HA=-C^*C,
\end{equation}
and
\begin{equation}\label{c-un}
C=-D^{-1}B^*H.
\end{equation}
The property b) is equivalent to

b') there exists an Hermitian solution $G={\rm diag}(G_1,\ldots,G_N)$, with $G_k\in\mathbb{C}^{\gamma_k\times\gamma_k},\ k=1,\ldots,N$, of the Lyapunov equation
\begin{equation}\label{L'-un}
GA^*+AG=-BB^*,
\end{equation}
and
\begin{equation}\label{b-un}
B=-GC^*D^{-1}.
\end{equation}
\end{thm}
\begin{proof}
To obtain Theorem~\ref{thm:L-un} from Theorem~\ref{thm:Lyap} it suffices to show that any structured Hermitian solution to the Lyapunov equation \eqref{L-un} (resp., \eqref{L'-un}) is invertible. Let $H={\rm diag}(H_1,\ldots,H_N)$ be a structured Hermitian solution to  \eqref{L-un}, and $x\in\ker H$, i.e., $x={\rm col}_{1\leq k\leq N}(x_k)$ and $x_k\in\ker H_k,\ k=1,\ldots,N$. Then
$\left\langle HAx,x\right\rangle=\left\langle Ax,Hx\right\rangle=0,$
and equation \eqref{L-un} implies $Cx=0$. In particular, for every $k\in\{ 1,\ldots,N\}$ one can define $\tilde{x}={\rm col}(0,\ldots,0,x_k,0,\ldots,0)$ where $x_k\in\ker H_k$ is on the $k$-th block entry of $\tilde{x}$, and from $C\tilde{x}=0$ get $C_kx_k=0$. Thus, $\ker H_k\subset\ker C_k,\ k=1,\ldots,N$. Consider the  following block representations  with respect to the decompositions $\mathbb{C}^{\gamma_k}=\ker H_k\oplus {\rm ran}\, H_k$:
$$A_{ij}=\begin{pmatrix} A_{ij}^{(11)} &  A_{ij}^{(12)}\\
 A_{ij}^{(21)} &  A_{ij}^{(22)}
 \end{pmatrix},\quad C_k=\begin{pmatrix} 0 & C_k^{(2)}\end{pmatrix},\quad H_k=\begin{pmatrix} 0 & 0\\
 0 & H_k^{(22)}\end{pmatrix},$$
where $i,j,k=1,\ldots,N$. Then \eqref{L-un} implies
$$(A^*H+HA)_{ij}^{(12)}=(A_{ji}^*H_j+H_iA_{ij})^{(12)}={(A_{ji}^{(21)})}^*H_j^{(22)}=0,$$
and $A_{ji}^{(21)}=0,\ i,j=1,\ldots,N$. Therefore, for any $w\in\mathcal{F}_N$ we have $$(C\flat A)^{wg_k}=\begin{pmatrix}
0 & (C^{(2)}\flat A^{(22)})^{wg_k}\end{pmatrix},\ k=1,\ldots,N,$$ where $C^{(2)}={\rm row}_{1\leq k\leq N}(C_k^{(2)}),\ A^{(22)}=(A_{ij}^{(22)})_{i,j=1,\ldots,N}$. If there exists $k\in\{ 1,\ldots,N\}$ such that $\ker H_k\not=\{ 0\}$, then the pair $(C,A)$ is not observable, which contradicts to the assumption on $\alpha$. Thus, $H$ is invertible.

In a similar way one can show that any structured Hermitian solution $G={\rm diag}(G_1,\ldots,G_N)$ of the Lyapunov equation \eqref{L'-un} is invertible.
\end{proof}
A counterpart of Theorem~\ref{thm:inv} in the present case is the following theorem.
\begin{thm}\label{thm:inv-un}
Let $(C,A)$ be an observable pair of matrices $C\in\mathbb{C}^{q\times r},A\in\mathbb{C}^{r\times r}$ in the sense that $\mathbb{C}^r=\bigoplus_{k=1}^N\mathbb{C}^{r_k}$ and $\mathcal{O}_k$ has full column rank for each $k\in\{ 1,\ldots,N\}$. Then there exists a  matrix-unitary on $\mathcal{J}_N$ rational FPS $F$ with a minimal GR-realization $\alpha=(N;A,B,C,D;\mathbb{C}^r=\bigoplus_{k=1}^N\mathbb{C}^{r_k},\mathbb{C}^q)$ if and only if the Lyapunov equation \eqref{L-un} has a structured Hermitian solution $H={\rm diag}(H_1,\ldots,H_N)$. If such a solution $H$ exists, it is invertible, and possible choices of $D$ and $B$ are
\begin{equation}\label{db-un}
D_0=I_q,\quad B_0=-H^{-1}C^*.
\end{equation}
Finally, for a given such $H$, all other choices of $D$ and $B$ differ from $D_0$ and $B_0$ by a right multiplicative unitary constant matrix.
\end{thm}
The proof of Theorem~\ref{thm:inv-un} is a direct application of Theorem~\ref{thm:inv} and Theorem~\ref{thm:L-un}. One can prove analogously the following theorem which is a counterpart of Theorem~\ref{thm:inv'}.
\begin{thm}\label{thm:inv'-un}
Let $(A,B)$ be a controllable pair of matrices $A\in\mathbb{C}^{r\times r},B\in\mathbb{C}^{r\times q}$ in the sense that $\mathbb{C}^r=\bigoplus_{k=1}^N\mathbb{C}^{r_k}$ and $\mathcal{C}_k$ has full row rank for each $k\in\{ 1,\ldots,N\}$. Then there exists a  matrix-unitary on $\mathcal{J}_N$ rational FPS $F$ with a minimal GR-realization $\alpha=(N;A,B,C,D;\mathbb{C}^r=\bigoplus_{k=1}^N\mathbb{C}^{r_k},\mathbb{C}^q)$ if and only if the Lyapunov equation
\eqref{L'-un} has a structured Hermitian  solution $G={\rm diag}(G_1,\ldots,G_N)$. If such a solution $G$ exists, it is invertible, and possible choices of $D$ and $C$ are
\begin{equation}\label{dc-un}
D_0=I_q,\quad C_0=-B^*G^{-1}.
\end{equation}
Finally, for a given such $G$, all other choices of $D$ and $C$ differ from $D_0$ and $C_0$ by a left multiplicative unitary constant matrix.
\end{thm}
Let $\overline{A}=(\overline{A_1},\ldots,\overline{A_N})$ be an $N$-tuple of $r\times r$ matrices. A non-zero vector $x\in\mathbb{C}^r$ is called a \emph{common eigenvector for $\overline{A}$} if there exists $\lambda=(\lambda_1,\ldots,\lambda_N)\in\mathbb{C}^N$ (which is called a \emph{common eigenvalue for $\overline{A}$}) such that
$$\overline{A_k}x=\lambda_kx,\quad k=1,\ldots,N.$$
The following theorem, which is a multivariable non-commutative counterpart of statements a) and b) of Theorem~2.10 in \cite{AG1}, gives a necessary condition on a minimal GR-realization of a matrix-unitary on $\mathcal{J}_N$ rational FPS.
\begin{thm}\label{thm:joint}
Let $F$ be a matrix-unitary on $\mathcal{J}_N$ rational FPS and  $\alpha$ be its minimal GR-realization, with the associated structured Hermitian matrix $H={\rm diag}(H_1,\ldots,H_N)$ and the associated inner products $[\,\cdot\,,\,\cdot\,]_{H_k},\ k=1,\ldots,N$. Let $P_k$ denote the orthogonal projection in $\mathbb{C}^\gamma$ onto the subspace $\{ 0\}\oplus\cdots\oplus\{ 0\}\oplus\mathbb{C}^{\gamma_k}\oplus\{ 0\}\oplus\cdots\oplus\{ 0\}$, and $\overline{A_k}=AP_k,\ k=1,\ldots,N$. If $x\in\mathbb{C}^\gamma$ is a common eigenvector for $\overline{A}$ corresponding to a common eigenvalue $\lambda\in\mathbb{C}^N$ then there exists $j\in\{ 1,\ldots,N\}$ such that ${\rm Re}\,\lambda_j\not=0$ and $[P_jx,P_jx]_{H_j}\not=0$. In particular, $\overline{A}$ has no common eigenvalues on $(i\mathbb{R})^N$.
\end{thm}
\begin{proof}
By \eqref{L-un}, we have for every $k\in\{ 1,\ldots,N\}$,
$$(\overline{\lambda_k}+\lambda_k)[P_kx,P_kx]_{H_k}=-\left\langle CP_kx,CP_kx\right\rangle.$$
Suppose that for all $k\in\{ 1,\ldots,N\}$ the left-hand side of this equality is zero, then $CP_kx=0$. Since for $\emptyset\not=w=g_{i_1}\cdots g_{i_{|w|}}\in\mathcal{F}_N$,
$$(C\flat A)^{wg_k}P_kx=CP_{i_1}\overline{A_{i_2}}\cdots \overline{A_{i_{|w|}}}\cdot\overline{A_k}x=\lambda_{i_2}\cdots\lambda_{i_{|w|}}\lambda_kCP_{i_1}x=0,$$
the observability of the pair $(C,A)$ implies $P_kx=0,\ k=1,\ldots,N$, i.e., $x=0$ which contradicts to the assumption that $x$ is a common eigenvector for $\overline{A}$. Thus, there exists $j\in\{ 1,\ldots,N\}$ such that $(\overline{\lambda_j}+\lambda_j)[P_jx,P_jx]_{H_j}\not=0$, as desired.
\end{proof}

\section{Matrix-$J$-unitary formal power series: A multivariable non-commutative analogue of the circle case}\label{sec:circle}
In this section we study a multivariable non-commutative analogue of
rational ${\mathbb C}^{q\times q}$-valued functions which are
$J$-unitary on the unit circle ${\mathbb T}$.
\subsection{Minimal Givone--Roesser realizations and the Stein equation}
Let $n\in{\mathbb N}$. We denote by ${\mathbb T}^{n\times n}$
the \emph{matrix unit circle}
$$
{\mathbb T}^{n\times n}=\left\{W\in{\mathbb C}^{n\times n}\,\,:\,\,
WW^*=I_n\right\},$$
i.e., the family of unitary $n\times n$ complex matrices. We will call the set
$\left({\mathbb T}^{n\times n}\right)^N$ the \emph{matrix unit
torus}. The set
$$
{\mathcal T}_N=\coprod\limits_{n\in{\mathbb N}}
\left({\mathbb T}^{n\times n}\right)^N$$
serves as a multivariable non-commutative counterpart of the unit circle.
Let $J=J^{-1}=J^*\in{\mathbb C}^{q\times q}$. We will say that a rational FPS $f$ is \emph{matrix-$J$-unitary on
${\mathcal T}_N$} if for every $n\in{\mathbb N}$,
$$f(W)(J\otimes I_n)f(W)^*=J\otimes I_n$$
at all points $W=(W_1,\ldots,W_N)\in\left({\mathbb T}^{n\times n}\right)^N$ where it is defined. In the following theorem we establish the relationship between matrix-$J$-unitary rational FPSs on $\mathcal{J}_N$ and on $\mathcal{T}_N$, their minimal GR-realizations, and the structured Hermitian solutions of the corresponding Lyapunov and Stein equations.
\begin{thm}\label{thm:Cayley}
Let $f$ be a matrix-$J$-unitary on $\mathcal{T}_N$ rational FPS, with a minimal GR-realization $\alpha$ of the form \eqref{min}, and let $a\in\mathbb{T}$ be such that $-\bar{a}\not\in\sigma(A)$. Then
\begin{equation}\label{Cayley}
F(z)=f(a(z_1-1)(z_1+1)^{-1},\ldots,a(z_N-1)(z_N+1)^{-1})
\end{equation}
is well defined as a rational FPS which is matrix-$J$-unitary on $\mathcal{J}_N$, and $F=T_\beta^{\rm nc}$, where $\beta=(N;A_a,B_a,C_a,D_a;\mathbb{C}^\gamma=\bigoplus_{k=1}^N\mathbb{C}^{\gamma_k},\mathbb{C}^q)$, with
\begin{equation}\label{beta}
\begin{array}{cc}
A_a=(aA-I_\gamma)(aA+I_\gamma)^{-1}, & B_a=\sqrt{2}(aA+I_\gamma)^{-1}aB,\\
C_a=\sqrt{2}C(aA+I_\gamma)^{-1}, & D_a=D-C(aA+I_\gamma)^{-1}aB.
\end{array}
\end{equation}
A GR-node $\beta$ is minimal, and its associated structured Hermitian matrix $H={\rm diag}(H_1,\ldots,H_N)$ is the unique invertible structured Hermitian solution of
\begin{equation}\label{H-un}
\begin{pmatrix}
A & B\\
C & D
\end{pmatrix}^*\begin{pmatrix}
H & 0\\
0 & J
\end{pmatrix}\begin{pmatrix}
A & B\\
C & D
\end{pmatrix}=\begin{pmatrix}
H & 0\\
0 & J
\end{pmatrix}.
\end{equation}
\end{thm}
\begin{proof}
For any $a\in\mathbb{T}$ and $n\in\mathbb{N}$ the Cayley transform
$$Z_0\longmapsto W_0=a(Z_0-I_n)(Z_0+I_n)^{-1}$$
maps $i\mathbb{H}^{n\times n}$ onto $\mathbb{T}^{n\times n}$, thus its simultaneous application to each matrix variable maps $(i\mathbb{H}^{n\times n})^N$ onto $(\mathbb{T}^{n\times n})^N$. Since the simultaneous application of the Cayley transform to each formal variable in a rational FPS gives a rational FPS, \eqref{Cayley} defines a rational FPS F. Since $f$ is matrix-$J$-unitary on $\mathcal{T}_N$, $F$ is matrix-$J$-unitary on $\mathcal{J}_N$. Moreover,
\begin{eqnarray*}
F(z) &=& D+C\left(I_\gamma-a(\Delta(z)-I_\gamma)(\Delta(z)+I_\gamma)^{-1}A\right)^{-1}\\
&\times & a(\Delta(z)-I_\gamma)(\Delta(z)+I_\gamma)^{-1}B\\
&=& D+C\left(\Delta(z)+I_\gamma-a(\Delta(z)-I_\gamma)A\right)^{-1}a(\Delta(z)-I_\gamma)B\\
&=& D+C\left(aA+I_\gamma-\Delta(z)(aA-I_\gamma)\right)^{-1}a(\Delta(z)-I_\gamma)B\\
&=& D+C(aA+I_\gamma)^{-1}\left(I_\gamma-\Delta(z)(aA-I_\gamma)(aA+I_\gamma)^{-1}\right)^{-1}\Delta(z)aB\\
&-& C(aA+I_\gamma)^{-1}\left(I_\gamma-\Delta(z)(aA-I_\gamma)(aA+I_\gamma)^{-1}\right)^{-1}aB\\
&=& D-C(aA+I_\gamma)^{-1}aB+C(aA+I_\gamma)^{-1}\\
&\times & \left(I_\gamma-\Delta(z)(aA-I_\gamma)(aA+I_\gamma)^{-1}\right)^{-1}\\
&\times & \Delta(z)\left(I_\gamma-(aA-I_\gamma)(aA+I_\gamma)^{-1}\right)aB\\
&=& D_a+C_a(I_\gamma-\Delta(z)A_a)^{-1}\Delta(z)B_a.
\end{eqnarray*}
Thus, $F=T_\beta^{\rm nc}$. Let us remark that the FPS
$$\varphi_k^a(z)=C_a(I_\gamma-\Delta(z)A_a)^{-1}\big|_{\mathbb{C}^{\gamma_k}}$$
(c.f. \eqref{phi_k}) has the coefficients
$$(\varphi_k^a)_w=(C_a\flat A_a)^{wg_k},\quad w\in\mathcal{F}_N.$$
Remark also that
\begin{eqnarray*}
\lefteqn{\tilde{\varphi}_k(z):= \varphi_k\left(a(z_1-1)(z_1+1)^{-1},\ldots,a(z_N-1)(z_N+1)^{-1}\right)}\\
&=& C\left( I_\gamma-a(\Delta(z)-I_\gamma)(\Delta(z)+I_\gamma)^{-1}A\right)^{-1}\big|_{\mathbb{C}^{\gamma_k}}\\
&=& C\left((\Delta(z)+I_\gamma)-a(\Delta(z)-I_\gamma)A\right)^{-1}(\Delta(z)+I_\gamma)\big|_{\mathbb{C}^{\gamma_k}}\\
&=& C\left((aA+I_\gamma)-\Delta(z)(aA-I_\gamma)\right)^{-1}(\Delta(z)+I_\gamma)\big|_{\mathbb{C}^{\gamma_k}}\\
&=& C(aA+I_\gamma)^{-1}\left(I_\gamma-\Delta(z)(aA-I_\gamma)(aA+I_\gamma)^{-1}\right)^{-1}(\Delta(z)+I_\gamma)\big|_{\mathbb{C}^{\gamma_k}}\\
&=& \frac{1}{\sqrt{2}}\left(C_a(I_\gamma-\Delta(z)A_a)^{-1}\big|_{\mathbb{C}^{\gamma_k}}\right)(z_k+1)\\ &=&
\frac{1}{\sqrt{2}}\left(\varphi_k^a(z)\cdot z_k+\varphi_k^a(z)\right).
\end{eqnarray*}
Let $k\in\{ 1,\ldots,N\}$ be fixed. Suppose that $n\in\mathbb{N},\ n\geq {(q\gamma -1)}^{q\gamma -1}$ (for $q\gamma -1=0$ choose arbitrary $n\in\mathbb{N}$), and $x\in\bigcap_{Z\in\Gamma_n(\varepsilon)}\ker\varphi_k^a(Z)$, where $\Gamma_n(\varepsilon)$ is a neighborhood of the origin of $\mathbb{C}^{n\times n}$ where $\varphi_k^a(Z)$ is well defined, e.g., of the form \eqref{gamma} with $\varepsilon=\| A_a\|^{-1}$. Then, by Theorem~\ref{thm:trunc} and Theorem~\ref{thm:kers}, one has
\begin{eqnarray*}
\bigcap\limits_{Z\in\Gamma_n(\varepsilon)}\ker \varphi_k^a(Z) = \left(\bigcap\limits_{w\in\mathcal{F}_N:\, |w|\leq q\gamma -1}\ker\,( \varphi_k^a)_w\right)\otimes\mathbb{C}^n\\
=\left(\bigcap\limits_{w\in\mathcal{F}_N:\, |w|\leq q\gamma -1}\ker\,(C_a\flat A_a)^{wg_k}\right)\otimes\mathbb{C}^n = \left(\ker\tilde{O}_k(\beta)\right)\otimes\mathbb{C}^n.
\end{eqnarray*}
Thus, there exist $l\in\mathbb{N},\ \{ u^{(\mu)}\}_{\mu=1}^l\subset\ker\tilde{O}_k(\beta),\ \{y^{(\mu)}\}_{\mu=1}^l\subset\mathbb{C}^n$ such that
\begin{equation}\label{x}
x=\sum_{\mu=1}^lu^{(\mu)}\otimes y^{(\mu)}.
\end{equation}
Since $\left(\varphi_k^a(z)\cdot z_k\right)_{wg_k}=(C_a\flat A_a)^{wg_k}$ for $w\in\mathcal{F}_N$, and $\left(\varphi_k^a(z)\cdot z_k\right)_{w'}=0$ for $w'\not=wg_k$ with any $w\in\mathcal{F}_N$, \eqref{x} implies that $\varphi_k^a(Z)(I_{\gamma_k}\otimes Z_k)x\equiv 0$. Thus,
$$\tilde{\varphi}_k(Z)x=\frac{1}{\sqrt{2}}\left(\varphi_k^a(Z)(I_{\gamma_k}\otimes Z_k)+\varphi_k^a(Z)\right)x\equiv 0.$$
Since the Cayley transform
$ a(\Delta(z)-I_\gamma)(\Delta(z)+I_\gamma)^{-1}$ maps an open and dense
subset of the set of matrices of the form
$\Delta(Z)={\rm diag}\,(Z_1,\ldots,Z_N)$, $Z_j\in\mathbb{C}^{\gamma_j\times\gamma_j},\ j=1,\ldots,N$, onto an open and dense subset of the same set,
$$\varphi_k(Z)x=(C\otimes I_n)(I_\gamma-\Delta(Z)(A\otimes I_n))^{-1}x\equiv 0.$$
Since the GR-node $\alpha$ is observable, by Theorem~\ref{thm:com} we get $x=0$. Therefore,
$$\bigcap_{Z\in\Gamma_n(\varepsilon)}\ker\varphi_k^a(Z)=0,\quad k=1,\ldots,N.$$ Applying Theorem~\ref{thm:com} once again, we obtain the observability of the GR-node $\beta$. In the same way one can prove the controllability of $\beta$. Thus, $\beta$ is minimal.

Note that
\begin{eqnarray}\nonumber
\lefteqn{\begin{pmatrix}
A & B\\
C & D
\end{pmatrix}^*
\begin{pmatrix}
H & 0\\
0 & J
\end{pmatrix}
\begin{pmatrix}
A & B\\
C & D
\end{pmatrix}-
\begin{pmatrix}
H & 0\\
0 & J
\end{pmatrix}=}\\
&=& \begin{pmatrix}
A^*HA+C^*JC-H & A^*HB+C^*JD\\
B^*HA+D^*JC & B^*HB+D^*JD-J
\end{pmatrix}.\label{blocks}
\end{eqnarray}
Since $-\bar{a}\notin\sigma(A)$, the matrix $(aA+I_\gamma)^{-1}$ is well defined, as well as $A_a=(aA-I_\gamma)(aA+I_\gamma)^{-1}$, and $I_\gamma-A_a=2(aA+I_\gamma)^{-1}$ is invertible. Having this in mind, one can deduce from \eqref{beta} the following relations:
$$
A^*HA+C^*JC-H =2(I_\gamma-A_a^*)^{-1}(A_a^*H+HA_a+C_a^*JC_a)(I_\gamma-A_a)^{-1}$$
\begin{eqnarray*}
B^*HA+D^*JC &=& \sqrt{2}(B_a^*H+D_a^*JC_a)(I_\gamma-A_a)^{-1}\\
&+& \sqrt{2}B_a^*(I_\gamma-A_a^*)^{-1}(A_a^*H+HA_a+C_a^*JC_a)(I_\gamma-A_a)^{-1}
\end{eqnarray*}
\begin{eqnarray*}
\lefteqn{B^*HB+D^*JD-J}\\ &=& B_a^*(I_\gamma-A_a^*)^{-1}(A_a^*H+HA_a+C_a^*JC_a)(I_\gamma-A_a)^{-1}B_a\\
&+& (B_a^*H+D_a^*JC_a)(I_\gamma-A_a)^{-1}B_a
+ B_a^*(I_\gamma-A_a^*)^{-1}(C_a^*JD_a+HB_a).
\end{eqnarray*}
Thus,  $A,B,C,D,H$ satisfy \eqref{H-un} if and only if  $A_a,B_a,C_a,D_a,H$ satisfy \eqref{L} and \eqref{b} (in the place of $A,B,C,D,H$ therein), which completes the proof.
\end{proof}
We will call the invertible Hermitian solution $H={\rm diag}(H_1,\ldots,H_N)$ of \eqref{H-un}, which is determined uniquely by a minimal GR-realization $\alpha$ of a matrix-$J$-unitary on $\mathcal{T}_N$ rational FPS $f$, the \emph{associated structured Hermitian matrix} (associated with a minimal GR-realization $\alpha$ of $f$). Let us note also that since for the GR-node $\beta$ from Theorem~\ref{thm:Cayley} a pair of the equalities \eqref{L} and \eqref{b} is equivalent to a pair of the equalities \eqref{L'} and \eqref{c}, the equality \eqref{H-un} is equivalent to
\begin{equation}\label{H-un'}
\begin{pmatrix}
A & B\\
C & D
\end{pmatrix}\begin{pmatrix}
H^{-1} & 0\\
0 & J
\end{pmatrix}\begin{pmatrix}
A & B\\
C & D
\end{pmatrix}^*=\begin{pmatrix}
H^{-1} & 0\\
0 & J
\end{pmatrix}.
\end{equation}
\begin{rem}
Equality \eqref{H-un} can be replaced by the following three equalities:
\begin{eqnarray}
H-A^*HA &=& C^*JC,\label{s}\\
D^*JC &=&-B^*HA,\label{adj}\\
J-D^*JD &=& B^*HB,\label{d}
\end{eqnarray}
and equality \eqref{H-un'} can be replaced by
\begin{eqnarray}
H^{-1}-AH^{-1}A^* &=& BJB^*,\label{s'}\\
DJB^* &=&-CH^{-1}A^*,\label{adj'}\\
J-DJD^* &=& CH^{-1}C^*.\label{d'}
\end{eqnarray}
\end{rem}
Theorem~\ref{thm:Cayley} allows to obtain a counterpart of the results from Section~\ref{sec:line} in the setting of rational FPSs which are matrix-$J$-unitary on $\mathcal{T}_N$. We will skip the proofs  when it is clear how to get them.
\begin{thm}\label{thm:j-u}
Let $f$ be a rational FPS and  $\alpha$ be its minimal GR-realization of the form \eqref{min}. Then $f$ is matrix-$J$-unitary on $\mathcal{T}_N$ if and only if there exists an invertible Hermitian matrix $H={\rm diag}(H_1,\ldots,H_N)$, with $H_k\in\mathbb{C}^{\gamma_k\times\gamma_k},\ k=1,\ldots,N$, which satisfies \eqref{H-un}, or equivalently, \eqref{H-un'}.
\end{thm}
\begin{rem}\label{rem:j-un}
In the same way as in \cite[Theorem~3.1]{AG1} one can show that if a rational
FPS $f$ has a (not necessarily minimal) GR-realization \eqref{qq}
which satisfies \eqref{H-un} (resp., \eqref{H-un'}), with an Hermitian
invertible matrix $H={\rm diag}(H_1,\ldots,H_N)$, then for any $n\in\mathbb{N}$,
\begin{eqnarray}
f(Z')^*(J\otimes I_n)f(Z) &=& J\otimes I_n-(B^*\otimes I_n)\left(I_\gamma\otimes I_n-\Delta(Z'^*)(A^*\otimes I_n)\right)^{-1}\nonumber\\
&\times & (H\otimes I_n)(I_\gamma\otimes I_n-\Delta(Z')^*\Delta(Z))\nonumber\\
&\times & \left(I_\gamma\otimes I_n-(A\otimes I_n)\Delta(Z)\right)^{-1}(B\otimes I_n)\label{s-id}
\end{eqnarray}
and respectively,
\begin{eqnarray}
f(Z)(J\otimes I_n)f(Z')^* &=& J\otimes I_n-(C\otimes I_n)\left(I_\gamma\otimes I_n-\Delta(Z)(A\otimes I_n)\right)^{-1}\nonumber\\
&\times & (I_\gamma\otimes I_n-\Delta(Z)\Delta(Z')^*)(H^{-1}\otimes I_n)\nonumber\\
&\times & \left(I_\gamma\otimes I_n-(A^*\otimes I_n)\Delta(Z')^*\right)^{-1}(C^*\otimes I_n),\label{s-id'}
\end{eqnarray}
at all the points $Z,Z'\in {(\mathbb{C}^{n\times n})}^N$ where it is defined,
which implies that $f$ is matrix-$J$-unitary on $\mathcal{T}_N$.
Moreover, the same statement holds true if $H={\rm diag}(H_1,\ldots,H_N)$ in \eqref{H-un} and \eqref{s-id} is not supposed to be invertible, and if $H^{-1}={\rm diag}(H_1^{-1},\ldots,H_N^{-1})$ in \eqref{H-un'} and \eqref{s-id'} is replaced by any Hermitian, not necessarily invertible matrix $Y={\rm diag}(Y_1,\ldots,Y_N)$.
\end{rem}
\begin{thm}\label{thm:c'=o'}
Let $f$ be a matrix-$J$-unitary on $\mathcal{T}_N$ rational FPS, and $\alpha$ be its GR-realization. Let $H={\rm diag}(H_1,\ldots,H_N)$ with $H_k\in\mathbb{C}^{r_k\times r_k},\ k=1,\ldots,N$,  be an Hermitian invertible matrix satisfying \eqref{H-un} or, equivalently, \eqref{H-un'}. Then $\alpha$ is observable if and only if $\alpha$ is controllable.
\end{thm}
\begin{proof}
Let $a\in\mathbb{T},\ -\bar{a}\notin\sigma(A)$. Then $F$ defined by \eqref{Cayley} is a matrix-$J$-unitary on $\mathcal{J}_N$ rational FPS, and \eqref{beta} is its GR-realization. As shown in the proof of Theorem~\ref{thm:Cayley}, $\alpha$ is observable (resp., controllable) if and only if so is $\beta$. Since by Theorem~\ref{thm:Cayley} the GR-node $\beta$ satisfies \eqref{L} and \eqref{b} (equivalently, \eqref{L'} and \eqref{c}), Theorem~\ref{thm:c=o} implies the statement of the present theorem.
\end{proof}
\begin{thm}\label{thm:a-inv}
Let $f$ be a matrix-$J$-unitary on $\mathcal{T}_N$ rational FPS and $\alpha$ be its minimal GR-realization of the form \eqref{min}, with the associated structured Hermitian matrix $H$. If $D=f_\emptyset$ is invertible then so is $A$, and
\begin{equation}\label{a-inv}
A^{-1}=H^{-1}(A^\times)^*H.
\end{equation}
\end{thm}
\begin{proof}
It follows from \eqref{adj} that $C=-JD^{-*}B^*HA$. Then \eqref{s} turns into
$$H-A^*HA=C^*J(-JD^{-*}B^*HA)=-C^*D^{-*}B^*HA,$$
which implies that $H=(A^\times)^*HA$, and \eqref{a-inv} follows.
\end{proof}
The following two lemmas, which are used in the sequel, can be found in \cite{AG1}.
\begin{lem}\label{lem:inv}
Let $A\in\mathbb{C}^{r\times r},\ C\in\mathbb{C}^{q\times r}$, where $A$ is invertible. Let $H$ be an invertible Hermitian matrix and $J$ be a signature matrix such that
$$H-A^*HA=C^*JC.$$
Let $a\in\mathbb{T},\ a\notin\sigma(A)$. Define
\begin{eqnarray}
D_a &=& I_q-CH^{-1}(I_r-aA^*)^{-1}C^*J,\label{da}\\
B_a &=& -H^{-1}A^{-*}C^*JD_a.\label{ba}
\end{eqnarray}
Then
$$\begin{pmatrix}
A & B_a\\
C & D_a
\end{pmatrix}^*
\begin{pmatrix}
H & 0\\
0 & J
\end{pmatrix}
\begin{pmatrix}
A & B_a\\
C & D_a
\end{pmatrix}=
\begin{pmatrix}
H & 0\\
0 & J
\end{pmatrix}.$$
\end{lem}
\begin{lem}\label{lem:inv'}
Let $A\in\mathbb{C}^{r\times r},\ B\in\mathbb{C}^{r\times q}$, where $A$ is invertible. Let $H$ be an invertible Hermitian matrix and $J$ be a signature matrix such that
$$H^{-1}-AH^{-1}A^*=BJB^*.$$
Let $a\in\mathbb{T},\ a\notin\sigma(A)$. Define
\begin{eqnarray}
D_a' &=& I_q-JB^*(I_r-aA^*)^{-1}HB,\label{daprime}\\
C_a' &=& -D_a'JB^*A^{-*}H.\label{caprime}
\end{eqnarray}
Then
$$\begin{pmatrix}
A & B\\
C_a' & D_a'
\end{pmatrix}
\begin{pmatrix}
H^{-1} & 0\\
0 & J
\end{pmatrix}
\begin{pmatrix}
A & B\\
C_a' & D_a'
\end{pmatrix}^*=
\begin{pmatrix}
H^{-1} & 0\\
0 & J
\end{pmatrix}.$$
\end{lem}
\begin{thm}\label{thm:inv-t}
Let $(C,A)$ be an observable pair of matrices $C\in\mathbb{C}^{q\times r},A\in\mathbb{C}^{r\times r}$ in the sense that $\mathbb{C}^r=\bigoplus_{k=1}^N\mathbb{C}^{r_k}$ and $\mathcal{O}_k$ has full column rank for each $k\in\{ 1,\ldots,N\}$. Let $A$ be invertible and  $J\in\mathbb{C}^{q\times q}$ be a signature matrix. Then there exists a  matrix-$J$-unitary on $\mathcal{T}_N$ rational FPS $f$ with a minimal GR-realization $\alpha=(N;A,B,C,D;\mathbb{C}^r=\bigoplus_{k=1}^N\mathbb{C}^{r_k},\mathbb{C}^q)$ if and only if the Stein equation \eqref{s} has a structured solution $H={\rm diag}(H_1,\ldots,H_N)$ which is both Hermitian and invertible. If such a solution $H$ exists, possible choices of $D$ and $B$ are $D_a$ and $B_a$ defined in \eqref{da} and \eqref{ba}, respectively. For a given such $H$, all other choices of $D$ and $B$ differ from $D_a$ and $B_a$ by a right multiplicative $J$-unitary constant matrix.
\end{thm}
\begin{proof}
Let $H={\rm diag}(H_1,\ldots,H_N)$ be a structured solution of the Stein equation \eqref{s} which is both Hermitian and invertible, $D_a$ and $B_a$ are defined as in \eqref{da} and \eqref{ba}, respectively, where $a\in\mathbb{T},\ a\notin\sigma(A)$. Set $\alpha_a=(N;A,B_a,C,D_a;\mathbb{C}^r=\bigoplus_{k=1}^N\mathbb{C}^{r_k},\mathbb{C}^q)$. By Lemma~\ref{lem:inv} and due to Remark~\ref{rem:j-un}, the transfer function $T_\alpha^{\rm nc}$ of $\alpha_a$ is a matrix-$J$-unitary on $\mathcal{T}_N$ rational FPS. Since $\alpha_a$ is observable, by Theorem~\ref{thm:c'=o'} $\alpha_a$ is controllable, and thus, minimal.

Conversely, if $\alpha=(N;A,B,C,D;\mathbb{C}^r=\bigoplus_{k=1}^N\mathbb{C}^{r_k},\mathbb{C}^q)$ is a minimal GR-node whose transfer function is  matrix-$J$-unitary on $\mathcal{T}_N$ then by Theorem~\ref{thm:j-u} there exists a solution $H={\rm diag}(H_1,\ldots,H_N)$ of the Stein equation \eqref{s} which is both Hermitian and invertible. The rest of the proof is analogous to the one of Theorem~\ref{thm:inv}.
\end{proof}

Analogously, one can obtain the following.
\begin{thm}\label{thm:inv-t'}
Let $(A,B)$ be a controllable pair of matrices $A\in\mathbb{C}^{r\times r},\ B\in\mathbb{C}^{r\times q}$ in the sense that $\mathbb{C}^r=\bigoplus_{k=1}^N\mathbb{C}^{r_k}$ and $\mathcal{C}_k$ has full row rank for each $k\in\{ 1,\ldots,N\}$. Let $A$ be invertible and  $J\in\mathbb{C}^{q\times q}$ be a signature matrix. Then there exists a  matrix-$J$-unitary on $\mathcal{T}_N$ rational FPS $f$ with a minimal GR-realization $\alpha=(N;A,B,C,D;\mathbb{C}^r=\bigoplus_{k=1}^N\mathbb{C}^{r_k},\mathbb{C}^q)$ if and only if the Stein equation
\begin{equation}\label{s-eq'}
G-AGA^*=BJB^*
\end{equation}
 has a structured solution $G={\rm diag}(G_1,\ldots,G_N)$ which is both Hermitian and invertible. If such a solution $G$ exists, possible choices of $D$ and $C$ are $D_a'$ and $C_a'$ defined in \eqref{da} and \eqref{ba}, respectively, where $H=G^{-1}$.
 For a given such $G$, all other choices of $D$ and $C$ differ from $D_a'$ and $C_a'$ by a left multiplicative $J$-unitary constant matrix.
\end{thm}

\subsection{The associated structured Hermitian matrix}\label{sec:Herm-circle}
In this subsection we give the  analogue of the results of Section~\ref{sec:Herm-line}. The
proofs are similar and will be omitted.
\begin{lem}
\label{lem:h'}
Let $f$ be a matrix-$J$-unitary on ${\mathcal T}_N$ rational
FPS and
$\alpha^{(i)}=(N;A^{(i)},B^{(i)},C^{(i)},D;{\mathbb C}^\gamma=
\bigoplus_{k=1}^N{\mathbb C}^{\gamma_k},{\mathbb C}^q)$ be its
minimal GR-realizations, with the associated structured Hermitian
matrices $H^{(i)}={\rm diag}(H_1^{(i)},\ldots, H_N^{(i)})$,
$i=1,2$. Then $\alpha^{(1)}$ and $\alpha^{(2)}$ are similar, that is
\begin{equation*}
C^{(1)}=C^{(2)}T,\quad TA^{(1)}=A^{(2)}T,\quad{\rm and}\quad
TB^{(1)}=B^{(2)},
\end{equation*}
for a uniquely defined invertible matrix $T={\rm diag}~(T_1,\ldots, T_N)\in\mathbb{C}^{\gamma\times\gamma}$
and
\begin{equation*}
H_k^{(1)}=T_k^*H_k^{(2)}T_k,\qquad k=1,\ldots, N.
\end{equation*}
In particular, the matrices $H_k^{(1)}$ and $H_k^{(2)}$ have the same
signature.
\end{lem}
\begin{thm}\label{thm:neg'}
Let $f$ be a matrix-$J$-unitary on ${\mathcal T}_N$ rational
FPS, and let
$\alpha$ be its minimal GR-realization of the form \eqref{min}, with the associated structured Hermitian matrix $H={\rm diag}(H_1,\ldots , H_N)$.
Then for each $k\in\left\{1,\ldots, N\right\}$ the number of negative
eigenvalues of the matrix $H_k$ is equal to the number of negative squares
of each of the kernels (on $\mathcal{F}_N$):
\begin{equation}
\label{kers'}
\begin{split}
K_{w,w^\prime}^{f,k}&=(C\flat A)^{wg_k}H_k^{-1}(A^*\sharp C^*)^{g_kw^{
\prime T}},\\
K_{w,{w^\prime}}^{f^*,k}&=(B^*\flat A^*)^{wg_k}
H_k(A\sharp B)^{g_kw'^T}.
\end{split}
\end{equation}
Finally, for $k\in\left\{1,\ldots, N\right\}$ let ${\mathcal K}_k(f)$
(resp., ${\mathcal K}_k(f^*)$) be the span of the functions
$w\mapsto K_{w,w^\prime}^{f,k}c$ (resp.,
$w\mapsto K_{w,w^\prime}^{f^*,k}c$) where $w^\prime\in\mathcal{F}_N$ and
$c\in{\mathbb C}^q$. Then
$$
{\rm dim}~{\mathcal K}_k(f)={\rm dim}~{\mathcal K}_k(f^*)=\gamma_k.
$$
\end{thm}
We will denote by $\nu_k(f)$ the number of negative squares of either of the
functions defined in \eqref{kers'}.
\begin{thm}\label{thm:fact'}
Let $f_i$, $i=1,2$, be two matrix-$J$-unitary on ${\mathcal T}_N$ rational
FPSs, with minimal GR-realizations
$$\alpha^{(i)}=\left(N;A^{(i)},B^{(i)},C^{(i)},D;{\mathbb C}^{\gamma^{(i)}}=
\bigoplus_{k=1}^N{\mathbb C}^{\gamma_k^{(i)}},{\mathbb C}^q\right)$$
and the associated structured Hermitian matrices
$H^{(i)}={\rm diag}(H_1^{(i)},\ldots, H_N^{(i)})$. Assume that the product $\alpha=\alpha^{(1)}\alpha^{(2)}$ is
a  minimal GR-node. Then, for each $k\in\left\{1,\ldots, N\right\}$ the
matrix
\begin{equation}
\label{hk'}
H_k=\begin{pmatrix}H_k^{(1)}&0\\0&H_k^{(2)}\end{pmatrix}\in
{\mathbb C}^{(\gamma_k^{(1)}+\gamma_k^{(2)})\times
(\gamma_k^{(1)}+\gamma_k^{(2)})}
\end{equation}
is the associated $k$-th Hermitian matrix for
$\alpha=\alpha^{(1)}\alpha^{(2)}$.
\end{thm}
\begin{cor}
\label{cor:neg'}
Let $f_1$ and $f_2$ be two matrix-$J$-unitary on ${\mathcal T}_N$ rational
FPSs, and assume that the factorization $f=f_1f_2$ is minimal. Then,
$$\nu(f_1f_2)=\nu(f_1)+\nu(f_2).$$
\end{cor}
\subsection{Minimal matrix-$J$-unitary factorizations}
In this subsection we consider minimal factorizations of  matrix-$J$-unitary on ${\mathcal T}_N$ rational FPSs
 into
two factors, both of which are also matrix-$J$-unitary on ${\mathcal T}_N$ rational
FPSs.
Such factorizations will be called \emph{minimal matrix-$J$-unitary factorizations.}

The following theorem is analogous to its one-variable counterpart
\cite[Theorem 3.7]{AG1} and proved in the same way.
\begin{thm}\label{thm:factor-circle}
Let $f$ be a matrix-$J$-unitary on ${\mathcal T}_N$ rational
FPS and $\alpha$ be its  minimal GR-realization of the form \eqref{min}, with the associated structured Hermitian matrix
$H={\rm diag}(H_1,\ldots, H_N)$, and assume that $D$ is invertible. Let
$\mathcal{M}=\bigoplus_{k=1}^N \mathcal{M}_k$ be an $A$-invariant subspace of $\mathbb{C}^\gamma$, which is non-degenerate in the associated inner product $[\,\cdot\, ,\,\cdot\,]_H$ and
  such that $M_k\subset
{\mathbb C}^{\gamma_k}$, $k=1,\ldots,N$. Let $\Pi={\rm diag}(\Pi_1,\ldots,
\Pi_N)$ be a projection defined by
$$\ker \Pi=M,\quad{and}\quad {\rm ran}\,\Pi=M^{[\perp]},$$
that is
$$
\ker \Pi_k=M_k,\quad{and}\quad {\rm ran}\,\Pi_k=M_k^{[\perp]}\quad for
\quad k=1,\ldots, N.$$
Then $f(z)=f_1(z)f_2(z)$, where
\begin{eqnarray}
\label{f1}
f_1(z)&=&\left[I_q+C(I_\gamma-\Delta(z)A)^{-1}\Delta(z)(I_\gamma-\Pi)
BD^{-1}\right]D_1,\\
f_2(z)&=&D_2\left[I_q+D^{-1}C\Pi
(I_\gamma-\Delta(z)A)^{-1}\Delta(z)B\right],
\label{f2}
\end{eqnarray}
with
$$D_1=I_q-CH^{-1}(I_\gamma-aA^*)^{-1}C^*J,\qquad D=D_1D_2,$$
where $a\in{\mathbb T}$ belongs to the resolvent set of $A_1$, and where
$$C_1=C\big|_\mathcal{M},\quad A_1=A\big|_\mathcal{M},\quad H_1=P_\mathcal{M}H\big|_\mathcal{M}$$
(with $P_\mathcal{M}$ being the orthogonal projection onto $\mathcal{M}$ in the standard
metric of ${\mathbb C}^\gamma$),
is a minimal matrix-$J$-unitary factorization of $f$.

Conversely, any minimal matrix-$J$-unitary factorization of $f$ can be
obtained in such a way, and the correspondence between minimal matrix-$J$-unitary
factorizations of $f$ with $f_1(\overline{a},\ldots, \overline{a})=I_q$
and non-degenerate subspaces of $A$ of the form $\mathcal{M}=\bigoplus_{k=1}^N\mathcal{M}_k$, with
$\mathcal{M}_k\subset{\mathbb C}^{\gamma_k},\ k=1,\ldots,N$, is one--to--one.
\end{thm}
\begin{rem}
In the proof of Theorem~\ref{thm:factor-circle}, as well as of Theorem~\ref{thm:min-fact}, we make use of Theorem \ref{thm:fact} and Corollary \ref{cor:fact}.
\end{rem}
\begin{rem}
Minimal matrix-$J$-unitary factorizations do not always exist, even in the
case $N=1$. See \cite{AG1} for examples in that case.
\end{rem}
\subsection{Matrix-unitary rational formal power series}
In this subsection we specialize some of the results in the present
section to the case $J=I_q$. We shall call corresponding rational
FPSs \emph{matrix-unitary on ${\mathcal T}_N$}.
\begin{thm}
Let $f$ be a rational FPS and
$\alpha$ be
its minimal GR-realization of the form \eqref{min}. Then $f$ is matrix-unitary on
${\mathcal T}_N$ if and only if:

$(a)$ There exists an Hermitian matrix $H={\rm diag}(H_1,\ldots, H_N)$
(with $H_k\in{\mathbb C}^{\gamma_k\times\gamma_k},\ k=1,\ldots,N$) such that
\begin{equation}
\begin{pmatrix}A&B\\ C&D\end{pmatrix}^*\begin{pmatrix} H&0\\ 0&I_q\end{pmatrix}
\begin{pmatrix}A&B\\ C&D\end{pmatrix}=\begin{pmatrix} H&0\\ 0&I_q\end{pmatrix}.
\label{HI-un}
\end{equation}

Condition $(a)$ is equivalent to:

$(a')$ There exists an Hermitian  matrix $G={\rm diag}~(G_1,
\ldots, G_N)$ (with $G_k\in{\mathbb C}^{\gamma_k\times\gamma_k},\ k=1,\ldots,N$) such that
\begin{equation}
\begin{pmatrix}A&B\\ C&D\end{pmatrix}\begin{pmatrix} G&0\\ 0&I_q\end{pmatrix}
\begin{pmatrix}A&B\\ C&D\end{pmatrix}^*=
\begin{pmatrix} G&0\\ 0&I_q\end{pmatrix}.
\label{GI-un}
\end{equation}
\end{thm}
\begin{proof}
The necessity follows from Theorem
\ref{thm:Cayley}. To prove the sufficiency,
suppose that the Hermitian matrix $H={\rm diag}(H_1,\ldots, H_N)$ satisfies
\eqref{HI-un} and let $a\in{\mathbb T}$ be such that $-\overline{a}\not\in
\sigma(A)$. Then, $H$ satisfies conditions \eqref{L-un} and \eqref{c-un} for
the GR-node
$\beta=(N;A_a,B_a,C_a,D_a;{\mathbb C}^{\gamma}=
\bigoplus_{k=1}^N{\mathbb C}^{\gamma_k},{\mathbb C}^q)$
defined by \eqref{beta} (this follows from the proof of Theorem~\ref{thm:Cayley}).
Thus, from Theorem~\ref{thm:L-un} and Theorem~\ref{thm:Cayley} we obtain that $f$ is
matrix-unitary on ${\mathcal T}_N$. Analogously, condition $(a')$
implies that the FPS $f$ is matrix-unitary on ${\mathcal T}_N$.
\end{proof}
A counterpart of Theorem~\ref{thm:inv-un} in the present case is the
following theorem:
\begin{thm}
Let $(C,A)$ be an observable pair of matrices in the sense that
${\mathcal O}_k$ has  full column rank for each $k=1,\ldots, N$.
Assume that $A\in{\mathbb C}^{r\times r}$ is invertible. Then there exists
a matrix-unitary on ${\mathcal T}_N$ rational FPS $f$ with a minimal
GR-realization
$\alpha=(N;A,B,C,D;{\mathbb C}^r=\bigoplus_{k=1}^N
{\mathbb C}^{r_k},{\mathbb C}^q)$
if and only if the Stein equation
\begin{equation}
H-A^*HA=C^*C
\label{SI-eq}
\end{equation}
has an Hermitian solution $H={\rm diag}(H_1,\ldots, H_N)$, with
$H_k\in{\mathbb C}^{r_k\times r_k}$, $k=1,\ldots, N$. If such a matrix $H$
exists, it is invertible, and possible choices of $D$ and $B$ are
$D_a$ and $B_a$ given by \eqref{da} and \eqref{ba} with $J=I_q$. Finally, for a given
$H={\rm diag}(H_1,\ldots, H_N)$, all other choices of $D$ and $B$ differ from
$D_a$ and $B_a$ by a right multiplicative unitary constant.
\end{thm}

A counterpart of Theorem \ref{thm:inv'-un} is the following theorem:

\begin{thm}
Let $(A,B)$ be a controllable pair of matrices, in the sense that
${\mathcal C}_k$ has  full row rank for each $k=1,\ldots, N$.
Assume that $A\in{\mathbb C}^{r\times r}$ is invertible. Then
there exists a matrix-unitary on ${\mathcal T}_N$ rational FPS $f$ with a
minimal GR-realization
$\alpha=(N;A,B,C,D;{\mathbb C}^r=\bigoplus_{k=1}^N
{\mathbb C}^{r_k},{\mathbb C}^q)$
if and only if the Stein equation
\begin{equation}
\label{SI-eqprime}
G-AGA^*=BB^*
\end{equation}
has an Hermitian solution $G={\rm diag}(G_1,\ldots, G_N)$ with
$G_k\in{\mathbb G}^{r_k\times r_k}$, $k=1,\ldots, N$.
If such a matrix $G$ exists, it is invertible, and possible choices of $D$ and
$C$ are $D_a^\prime$ and $C_a^\prime$ given by \eqref{daprime} and
\eqref{caprime} with $H=G^{-1}$ and $J=I_q$. Finally, for a given
$G={\rm diag}(G_1,\ldots, G_N)$, all other choices of $D$ and $C$ differ from
$D_a^\prime$ and $C_a^\prime$ by a left multiplicative unitary constant.
\end{thm}

A counterpart of Theorem \ref{thm:joint} in the present case is the following:
\begin{thm}
Let $f$ be a matrix-unitary on ${\mathcal T}_N$ rational FPS and
$\alpha$
be its minimal GR-realization of the form \eqref{min}, with the associated structured Hermitian matrix
$H={\rm diag}(H_1,\ldots, H_N)$ and the associated $k$-th inner products
$[\cdot,\cdot]_{H_k}$, $k=1,\ldots, N$.
Let $P_k$ denote the orthogonal projection in ${\mathbb C}^\gamma$ onto the
subspace
$\left\{0\right\}\oplus\cdots\oplus\left\{0\right\}
\oplus{\mathbb C}^{\gamma_k}\oplus\left\{0\right\}\oplus\cdots
\oplus\left\{0\right\}$,
and set $\overline{A_k}=AP_k$ for $k=1,\ldots, N$. If $x\in{\mathbb C}^\gamma$
is a common eigenvector for $\overline{A}=\left\{\overline{A_1},
\ldots,\overline{A_N}\right\}$ corresponding to a common
eigenvalue $\lambda=(\lambda_1,\ldots, \lambda_N)\in{\mathbb C}^N$, then there
exists $j\in\left\{1,\ldots, N\right\}$ such that $|\lambda_j|
\not=1$ and $[P_jx,P_jx]_{H_j}\not =0$. In particular
$\overline{A}$ has no common eigenvalues on $\mathbb{T}^N$.
\end{thm}

The proof of this theorem relies on the equality
$$
(1-|\lambda_k|^2)[P_kx,P_kx]_{H_k}=\langle CP_kx,CP_kx\rangle,\quad k=1,\ldots,N,$$
and follows the same argument as the proof of Theorem~\ref{thm:joint}.

\section{Matrix-$J$-inner rational formal power series}\label{sec:inner}
\subsection{A multivariable non-commutative analogue of the halfplane case}
Let $n\in{\mathbb N}$. We define the \emph{matrix open right poly-halfplane} as the
set
$$
\left(\Pi^{n\times n}\right)^N=\left\{
Z=(Z_1,\ldots, Z_N)\in\left({\mathbb C}^{n\times n}\right)^N:\,
Z_k+Z_k^*>0,\ k=1,\ldots, N\right\},
$$
and the \emph{matrix closed right poly-halfplane} as the set
\[
\begin{split}
{\rm clos}\,\left(\Pi^{n\times n}\right)^N&=\left({\rm clos}\,\Pi^{n\times n}
\right)^N\\
&=
\left\{
Z=(Z_1,\ldots, Z_N)\in\left({\mathbb C}^{n\times n}\right)^N:\,
Z_k+Z_k^*\ge 0,\ k=1,\ldots, N\right\}.
\end{split}
\]
We also introduce
$${\mathcal P}_N=\coprod_{n\in{\mathbb N}}\left(\Pi^{n\times n}\right)^N\quad
{\rm and}\quad
{\rm clos}\,{\mathcal P}_N=\coprod_{n\in{\mathbb N}}~
{\rm clos}\,\left(\Pi^{n\times n}\right)^N.$$
It is clear that
$$
\left(i{\mathbb H}^{n\times n}\right)^N\subset {\rm clos}\,\left(
\Pi^{n\times n}\right)^N$$ is the \emph{essential} (or \emph{Shilov}) \emph{boundary} of the
matrix poly-halfplane $\left(\Pi^{n\times n}\right)^N$ (see \cite{Sh})
and that ${\mathcal J}_N\subset {\rm clos}\,\mathcal P_N$
(recall that ${\mathcal J}_N=\coprod_{n\in {\mathbb N}}
\left(i{\mathbb H}^{n\times n}\right)^N$).

Let $J=J^{-1}=J^*\in{\mathbb C}^{q\times q}$. A matrix-$J$-unitary on
${\mathcal J}_N$ rational FPS $F$ is called \emph{matrix-$J$-inner (in $\mathcal{P}_N$)} if for each $n\in{\mathbb N}$:
\begin{equation}
F(Z)(J\otimes I_n)F(Z)^*\le J\otimes I_n
\label{J-contr}
\end{equation}
at those points $Z\in{\rm clos}\,\left(\Pi^{n\times n}\right)^N$ where it is
defined (the set of such points is open and dense, in the relative topology,
in ${\rm clos}\,\left(\Pi^{n\times n}\right)^N$ since $F(Z)$ is a rational
matrix-valued function of the complex variables $(Z_k)_{ij},\ k=1,\ldots, N,\ i,j=1,\ldots, n$).

The following theorem is a counterpart of part $\left.a\right)$
of Theorem 2.16
of \cite{AG1}.

\begin{thm}
\label{thm:J-in}
Let $F$ be a matrix-$J$-unitary on ${\mathcal J}_N$ rational
FPS and
$\alpha$
be its minimal GR-realization of the form \eqref{min}. Then $F$ is
matrix-$J$-inner in ${\mathcal P}_N$ if and only if the associated structured
Hermitian matrix $H={\rm diag}(H_1,\ldots, H_N)$ is strictly positive.
\end{thm}
\begin{proof}
Let $n\in{\mathbb N}$. Equality \eqref{L-id} can be rewritten as
\begin{equation}
\label{L-idprime}
J\otimes I_n-F(Z)(J\otimes I_n){F(Z')}^*=\varphi(Z)\Delta(Z+{Z'}^*)
(H^{-1}\otimes I_n){\varphi(Z')}^*
\end{equation}
where $\varphi$ is a FPS defined by
$$
\varphi(z):=C(I_\gamma-\Delta(z)A)^{-1}\in{\mathbb C}^{q\times\gamma}
\left\langle\left\langle z_1,\ldots, z_N\right\rangle\right\rangle_{\rm rat},
$$
and \eqref{L-idprime} is well defined at all points
$Z,Z'\in({\mathbb C}^{n\times n})^N$ for which $$1\not\in\sigma\left(
\Delta(Z)(A\otimes I_n)\right),\quad 1\not\in\sigma\left(
\Delta(Z')(A\otimes I_n)\right).$$ Set
$
\varphi_k(z):=C(I_\gamma-\Delta(z)A)^{-1}\big|_{{\mathbb C}^{\gamma_k}}
\in{\mathbb C}^{q\times\gamma_k}
\left\langle\left\langle z_1,\ldots, z_N\right\rangle\right\rangle_{\rm rat},
\ k=1,\ldots, N.
$
Then \eqref{L-idprime} becomes:
\begin{equation}
\label{L-idprimeprime}
J\otimes I_n-F(Z)(J\otimes I_n){F(Z')}^*=\sum_{k=1}^N\varphi_k(Z)
(H_k^{-1}\otimes(Z_k+{Z'_k}^*)){\varphi_k(Z')}^*.
\end{equation}
Let $X\in{\mathbb H}^{n\times n}$ be some positive semidefinite matrix, let
$Y\in({\mathbb H}^{n\times n})^N$ be such that $1\not\in\sigma
(\Delta(iY)(A\otimes I_n))$, and set for $k=1,\ldots, N$:
$$
e_k:=(0,0,\ldots, 0,1,0,\ldots, 0)\in{\mathbb C}^N$$
with $1$ at the $k$-th place. Then for $\lambda\in{\mathbb C}$
set
$$
Z^{(k)}_{X,Y}(\lambda):=\lambda X\otimes e_k+iY=
(iY_1,\ldots, iY_{k-1}, \lambda X+iY_k,iY_{k+1},\ldots, iY_N).$$
Now,
\eqref{L-idprimeprime} implies that
\begin{eqnarray}
\lefteqn
{J\otimes I_n-F(Z_{X,Y}^{(k)}(\lambda))(J\otimes I_n){F(Z_{X,Y}^{(k)}(\lambda'))}^*
}\nonumber\\
&=
(\lambda+\overline{\lambda'})\varphi_k(Z_{X,Y}^{(k)}(\lambda))(H_k^{-1}\otimes X){\varphi_k(Z_{X,Y}^{(k)}(\lambda'))}^*.
\label{L-id-s}
\end{eqnarray}
The function $h(\lambda)=F(Z_{X,Y}^{(k)}(\lambda))$ is a rational function
of $\lambda\in{\mathbb C}$. It is easily seen from \eqref{L-id-s} that $h$ is
$(J\otimes I_n)$-inner in the open right halfplane. In particular, it is $(J\otimes I_n)$-contractive
in the closed right halfplane (this also follows directly from \eqref{J-contr}). Therefore
(see e.g. \cite{D}) the function
\begin{equation}
\Psi(\lambda,\lambda')=\frac{J\otimes I_n-F(Z_{X,Y}^{(k)}(\lambda))(J\otimes I_n)
{F(Z_{X,Y}^{(k)})(\lambda')}^*}{\lambda+\overline{\lambda'}}
\label{ker-psi}
\end{equation}
is a positive semidefinite kernel on ${\mathbb C}$: for every choice of
$r\in{\mathbb N}$, of points $\lambda_1,\ldots, \lambda_r\in{\mathbb C}$ for
which the matrices $\Psi(\lambda_j,\lambda_i)$ are well defined,
and vectors $c_1,\ldots, c_r \in
{\mathbb C}^q\otimes {\mathbb C}^n$ one has
$$\sum_{i,j=1}^rc_j^*\Psi(\lambda_j,\lambda_i)c_i\ge 0,$$
i.e., the matrix $\left(\Psi(\lambda_j,\lambda_i)\right)_{i,j=1,\ldots, r}$
is positive semidefinite. Since $\varphi_k(Z_{X,Y}^{(k)}(0))=\varphi_k(iY)$ is
well-defined, we obtain from \eqref{L-id-s} that $\Psi(0,0)$ is also well-defined
and
$$\Psi(0,0)=\varphi_k(iY)(H_k^{-1}\otimes X)\varphi_k(iY)^*\ge 0.$$
This inequality holds for every $n\in{\mathbb N}$, every positive semidefinite
$X\in{\mathbb H}^{n\times n}$ and every $Y\in({\mathbb H}^{n\times n})^N$.
Thus, for an arbitrary $r\in{\mathbb N}$ we can define $\widetilde{n}=nr$,
$\widetilde{Y}=(\widetilde{Y}_1,\ldots,\widetilde{Y}_N)\in({\mathbb H}^{\widetilde{n}\times\widetilde{n}})^N$, where $\widetilde{Y}_k={\rm diag}(Y^{(1)}_k,\ldots, Y^{(r)}_k)$ and $Y^{(j)}_k\in {\mathbb H}^{n\times n},\ k=1,\ldots,N,\ j=1,\ldots,r,$ such that $\varphi_k(i
\widetilde{Y})$ is
well defined,
$$
\widetilde{X}=\begin{pmatrix}I_n&\cdots&I_n\\
\vdots& &\vdots\\
I_n&\cdots&I_n\end{pmatrix}
\in{\mathbb C}^{n\times n}\otimes {\mathbb C}^{r\times r}\cong
{\mathbb C}^{\widetilde{n}\times\widetilde{n}}$$
and get
\[
\begin{split}
0&\le\varphi_k(i\widetilde{Y})(H_k^{-1}\otimes \widetilde{X})
\varphi_k(i\widetilde{Y})^*\\
&=
{\rm diag}(\varphi_k(iY^{(1)}),\ldots, \varphi_k(iY^{(r)}))\times\\
&\times
\left(H_k^{-1}\otimes\begin{pmatrix}I_n\\ \vdots\\ I_n\end{pmatrix}
\begin{pmatrix}I_n&\cdots &I_n\end{pmatrix}\right)
{\rm diag}(\varphi_k(iY^{(1)})^*,\ldots, \varphi_k(iY^{(r)})^*)\\
&=\begin{pmatrix}
\varphi_k(iY^{(1)})\\ \vdots\\ \varphi_k(iY^{(r)})\end{pmatrix}
(H_k^{-1}\otimes I_n)
\begin{pmatrix}
\varphi_k(iY^{(1)})^*& \cdots& \varphi_k(iY^{(r)})^*\end{pmatrix}\\
&=\left(
\varphi_k(iY^{(\mu)})(H_k^{-1}\otimes I_n)\varphi_k(iY^{(\nu)})^*
\right)_{\mu,\nu=1,\ldots, r}.
\end{split}
\]
Therefore, the function
$$K_k(iY,iY^{\prime})=\varphi_k(iY)
(H_k^{-1}\otimes I_n)\varphi_k(iY^{\prime})^*$$
is a positive semidefinite kernel on any subset of
$(i{\mathbb H}^{n\times n})^N$ where it is defined, and in particular in
some neighbourhood of the origin. One can extend this function to
\begin{equation}
\label{hol}
K_k(Z,Z^{\prime})=\varphi_k(Z)
(H_k^{-1}\otimes I_n)\varphi_k(Z^{\prime})^*
\end{equation}
at those points $Z, Z^{\prime}\in({\mathbb C}^{n\times n})^N\times ({\mathbb C}^{n\times n})^N$ where $\varphi_k$ is
defined. Thus, on some neighbourhood $\Gamma$ of the origin in
$({\mathbb C}^{n\times n})^N\times ({\mathbb C}^{n\times n})^N$, the function
$K_k(Z,Z^{\prime})$ is holomorphic in $Z$ and
anti-holomorphic in $Z^{\prime}$. On the other hand, it is well-known
(see e.g. \cite{Aron}) that one can construct a reproducing kernel
Hilbert space (which we will denote by
${\mathcal H}(K_k)$) with reproducing kernel
$K_k(iY,iY^{\prime})$, which is obtained
as the completion of
$${\mathcal H}_0={\rm span}\left\{K_k(\cdot, iY)x\,;\, iY\in
(i{\mathbb H}^{n\times n})^N\cap\Gamma,\,\ x\in{\mathbb C}^q\otimes
{\mathbb C}^n \right\}$$
with respect to the inner product
\begin{eqnarray*}
\lefteqn{\left\langle\sum_{\mu=1}^rK_k(\cdot, iY^{(\mu)})x_\mu,
\sum_{\nu=0}^\ell K_k(\cdot, iY^{(\nu)})x_\nu\right\rangle_{{\mathcal H}_0}}\\
&=&
\sum_{\mu=1}^r\sum_{\nu=1}^\ell
\left\langle K_k(iY^{(\nu)},iY^{(\mu)})x_\mu,x_\nu
\right\rangle_{{\mathbb C}^q\otimes{\mathbb C}^n}.
\end{eqnarray*}
The reproducing kernel property reads:
$$
\langle f(\cdot), K_k(\cdot,iY)x
\rangle_{{\mathcal H}(K_k)}=\langle f(iY),x
\rangle_{{\mathbb C}^q\otimes{\mathbb C}^n},
$$
and thus $K_k(iY, iY^{\prime})=\Phi(iY){\Phi(iY^{\prime})}^*$ where
$$\Phi(iY):\, f(\cdot)\mapsto f(iY)$$
is the evaluation map. In view of \eqref{hol},
the kernel $K_k(\cdot, \cdot)$ is extendable on
$\Gamma\times\Gamma$ to the function $K(Z,
Z^{\prime})$ which is holomorphic in $Z$ and antiholomorphic
in $Z^{\prime}$, all the elements of ${\mathcal H}(K_k)$ have holomorphic
continuations to $\Gamma$, and so has the function $\Phi(\cdot)$. Thus,
$$K_k(Z,Z^{\prime})=\Phi(Z)\Phi(Z^{\prime})^*$$
and so $K_k(Z,Z^{\prime})$ is a positive semidefinite kernel on
$\Gamma$. (We could also use \cite[Theorem 1.1.4, p.10]{ADRS}
to obtain this conclusion.) Therefore, for any choice of $\ell\in\mathbb{N}$ and $Z^{(1)},\ldots,
Z^{(\ell)}\in\Gamma$ the matrix
\begin{equation}
\label{pos}
\begin{split}
\lefteqn{\left(
\varphi_k(Z^{(\mu)})(H_k^{-1}\otimes I_n)\varphi_k(Z^{(\nu)})^*
\right)_{\mu,\nu=1,\ldots, \ell} }\\
&= \begin{pmatrix}
\varphi_k(Z^{(1)})\\ \vdots\\ \varphi_k(Z^{(\ell)})\end{pmatrix}\cdot
(H_k^{-1}\otimes I_n)\cdot
\begin{pmatrix}
\varphi_k(Z^{(1)})^*& \cdots& \varphi_k(Z^{(\ell)})^*\end{pmatrix}
\end{split}
\end{equation}
is positive semidefinite. Since the coefficients of the FPS
$\varphi_k$ are $(\varphi_k)_w=(C\flat A)^{wg_k},\ w\in\mathcal{F}_N$,
and since $\alpha$ is an observable GR-node, we have
$$
\bigcap_{w\in\mathcal{F}_N}\ker (C\flat A)^{wg_k}=\left\{0\right\}.$$
Hence, by Theorem~\ref{thm:kers}
we can chose $n,\ell\in{\mathbb N}$ and
$Z^{(1)},\ldots, Z^{(\ell)}\in{\Gamma}$ such that
$$
\bigcap_{j=1}^\ell\ker \varphi_k(Z^{(j)})=\left\{0\right\}.$$
Thus the matrix ${\rm col}_{j=1,\ldots,\ell}\left(\varphi_k(Z^{(j)})\right)$ has full column
rank. (We could also use Theorem~\ref{thm:com}.) From \eqref{pos} it then follows that $H_k^{-1}>0$. Since this holds for
all $k\in\{1,\ldots, N\}$,  we get $H>0$.

Conversely, if
$H>0$ then it follows from \eqref{L-idprime} that for every $n\in{\mathbb N}$ and
$Z\in\left(\Pi^{n\times n}\right)^N$ for which $1\not \in
\sigma(\Delta(Z)(A\otimes I_n))$, one has
$$
J\otimes I_n-F(Z)(J\otimes I_n)F(Z)^*\ge 0.$$
Therefore $F$ is matrix-$J$-inner in $\mathcal{P}_N$, and
the proof is complete.
\end{proof}
\begin{thm}
\label{thm:cons}
Let $F\in{\mathbb C}^{q\times q}
\left\langle\left\langle z_1,\ldots, z_N\right\rangle\right\rangle_{\rm rat}$ be
 matrix-$J$-inner in ${\mathcal P}_N$. Then $F$ has a minimal
GR-realization of the form \eqref{min}
with the associated structured Hermitian matrix $H=I_\gamma$.
This realization is unique up to a unitary similarity.
\end{thm}
\begin{proof}
Let
$$
\alpha^\circ=(N;A^\circ,B^\circ,C^\circ,D;{\mathbb C}^\gamma
=\bigoplus_{k=1}^N{\mathbb C}^{\gamma_k},
{\mathbb C}^q)$$
be a minimal GR-realization of $F$, with the  associated structured Hermitian matrix
$H^\circ={\rm diag}(H_1^\circ,\ldots, H_N^\circ)$. By Theorem~\ref{thm:J-in} the matrix $H^\circ$ is strictly positive.
Therefore, $(H^\circ)^{1/2}={\rm diag}((H_1^\circ)^{1/2},\ldots,
(H_N^\circ)^{1/2})$ is well defined and strictly positive, and
$$\alpha=(N;A,B,C,D;{\mathbb C}^\gamma=\bigoplus_{k=1}^N{\mathbb C}^{\gamma_k},
{\mathbb C}^q),$$
where
\begin{equation}
\label{abc}
A=(H^\circ)^{1/2}A^\circ(H^\circ)^{-1/2},\quad
B=(H^\circ)^{1/2}B^\circ,\quad
C=C^\circ (H^\circ)^{-1/2},
\end{equation}
is a minimal GR-realization of $F$ satisfying
\begin{eqnarray}
A^*+A&=&-C^*JC,
\label{LI}\\
B&=&-C^*JD,
\label{BI}
\end{eqnarray}
or equivalently,
\begin{eqnarray}
A^*+A &=& -BJB^*,
\label{LI'}\\
C &=& -DJB^*,
\label{CI'}
\end{eqnarray}
and thus having the associated structured Hermitian matrix $H=I_\gamma$. Since in this case the inner product $[ \,\cdot\, ,\,\cdot\,]_H$ coincides
with the standard inner product $\langle \,\cdot\, ,\,\cdot\,\rangle$ of
${\mathbb C}^\gamma$, by Remark~\ref{rem:un} this minimal
GR-realization with the property $H=I_\gamma$ is unique up to
unitary similarity.
\end{proof}
We remark that a one-variable counterpart of the latter result is
essentially contained in \cite{Bro}, \cite{L} (see also
\cite[Section 4.2]{Ar}).
\subsection{A multivariable non-commutative analogue of the disk case}
Let $n\in {\mathbb N}$. We define the \emph{matrix open unit polydisk} as
$$
\left({\mathbb D}^{n\times n}\right)^N=\left\{W=(W_1,\ldots, W_N)\in
\left({\mathbb C}^{n\times n}\right)^N:\  W_kW_k^*<I_n,\
k=1,\ldots, N\right\},$$
and the \emph{matrix closed unit polydisk} as
\begin{eqnarray*}
\lefteqn{{\rm clos}\left({\mathbb D}^{n\times n}\right)^N=\left({\rm clos}\,{\mathbb D}^{n\times n}\right)^N}\\
&=& \left\{W=(W_1,\ldots, W_N)\in
\left({\mathbb C}^{n\times n}\right)^N:\  W_kW_k^*\leq I_n,\
k=1,\ldots, N\right\}.
\end{eqnarray*}
The matrix unit torus $\left({\mathbb T}^{n\times n}\right)^N$ is the
essential (or Shilov) boundary of $\left({\mathbb D}^{n\times n}\right)^N$
(see \cite{Sh}). In our setting, the set
$${\mathcal D}_N=\coprod_{n\in{\mathbb N}}
\left({\mathbb D}^{n\times n}\right)^N\quad\left({\rm resp.,}\quad
{\rm clos}\,{\mathcal D}_N=\coprod_{n\in{\mathbb N}}
{\rm clos}\left({\mathbb D}^{n\times n}\right)^N\right)$$
is a multivariable non-commutative counterpart of the open (resp.,
closed) unit disk.

Let $J=J^{-1}=J^*\in{\mathbb C}^{q\times q}$. A
rational FPS $f$ which is matrix-$J$-unitary on
${\mathcal T}_N$ is called \emph{matrix-$J$-inner in ${\mathcal D}_N$}
if for every $n\in{\mathbb N}$:
\begin{equation}
\label{J-contrprime}
f(W)(J\otimes I_n)f(W)^*\le J\otimes I_n
\end{equation}
at those points $W\in {\rm clos}\left({\mathbb D}^{n\times n}\right)^N$
where it is defined. We note that the set of such points is open and dense (in
the relative topology) in ${\rm clos}\left({\mathbb D}^{n\times n}\right)^N$
since $f(W)$ is a rational matrix-valued function of the complex variables
$(W_k)_{ij},\ k=1,\ldots, N,\ i,j=1,\ldots, n$.
\begin{thm}
\label{thm:j-inprime}
Let $f$ be a rational FPS which is matrix-$J$-unitary
on ${\mathcal T}_N$, and let
$\alpha$  be its minimal GR-realization of the form \eqref{min}. Then $f$ is matrix-$J$-inner in ${\mathcal D}_N$
if and only if the associated structured Hermitian matrix $H={\rm diag}(H_1,\ldots, H_N)$ is
strictly positive.
\end{thm}
\begin{proof} The statement of this theorem follows from Theorem~\ref{thm:J-in} and Theorem~\ref{thm:Cayley},
since the Cayley transform defined in Theorem~\ref{thm:Cayley} maps each open matrix unit polydisk $({\mathbb D}^{n
\times n})^N$ onto the open right matrix poly-halfplane $(\Pi^{n\times n})^N$,
and the inequality \eqref{J-contrprime} turns into \eqref{J-contr}
for the function $F$ defined in \eqref{Cayley}.
\end{proof}
The following theorem is an analogue of Theorem \ref{thm:cons}.
\begin{thm}
\label{thm:consprime}
Let $f$ be a rational FPS which
is matrix-$J$-inner in ${\mathcal D}_N$.  Then there exists its
minimal GR-realization $\alpha$ of the form \eqref{min},
with the associated structured Hermitian matrix $H=I_\gamma$. Such a realization
is unique up to a unitary similarity.
\end{thm}
In the special case of Theorem~\ref{thm:consprime} where $J=I_q$ the FPS $f$ is called \emph{matrix-inner}, and the GR-node $\alpha$ satisfies
$$
\begin{pmatrix}A&B\\C&D\end{pmatrix}^*\begin{pmatrix}A&B\\C&D\end{pmatrix}=
I_{\gamma+q},$$
i.e., $\alpha$ is a \emph{unitary GR-node}, which
has been considered first by J.~Agler in \cite{Agler}. In what
follows
we will show that Theorem \ref{thm:consprime} for $J=I_q$ is a
special case of the theorem of J.~A.~Ball, G.~Groenewald and T.~Malakorn on unitary
GR-realizations of FPSs from
the non-commutative Schur--Agler class \cite{BGM2}, which becomes in
several aspects stronger in this special case.

Let ${\mathcal U}$ and ${\mathcal Y}$ be Hilbert spaces. Denote
by $L({\mathcal U},{\mathcal Y})$ the Banach space of
bounded linear operators from ${\mathcal U}$ into ${\mathcal Y}$.
A GR-node in the general setting of Hilbert spaces is
$$\alpha=(N;A,B,C,D;{\mathcal X}=\bigoplus_{k=1}^N{\mathcal
X}_k,{\mathcal U},{\mathcal Y}),$$
i.e., a collection of Hilbert spaces ${\mathcal X},{\mathcal
X}_1,\ldots, {\mathcal X}_N,{\mathcal U},{\mathcal Y}$ and
operators $A\in L({\mathcal X})= L({\mathcal
X},{\mathcal X})$, $B\in L({\mathcal U},{\mathcal X})$,
$C\in L({\mathcal X},{\mathcal Y})$, and $D\in
L({\mathcal U},{\mathcal Y})$. Such a GR-node $\alpha$ is called
\emph{unitary} if
$$
\begin{pmatrix}A&B\\
C&D\end{pmatrix}^*\begin{pmatrix}A&B\\
C&D\end{pmatrix}=I_{{\mathcal X}\oplus {\mathcal U}},\quad \begin{pmatrix}A&B\\
C&D\end{pmatrix}\begin{pmatrix}A&B\\
C&D\end{pmatrix}^*= I_{{\mathcal X}\oplus {\mathcal Y}},
$$
i.e., {\scriptsize $\begin{pmatrix}A&B\\
C&D\end{pmatrix}$} is a unitary operator from ${\mathcal
X}\oplus{\mathcal U}$ onto ${\mathcal X}\oplus {\mathcal Y}$. The
\emph{non-commutative transfer function of $\alpha$} is
\begin{equation}\label{nctf}
T_\alpha^{\rm nc}(z)=D+C(I-\Delta(z)A)^{-1}\Delta(z)B,
\end{equation}
 where
the expression \eqref{nctf} is understood as a FPS from $L(\mathcal{U,Y})\left\langle \left\langle z_1,\ldots,z_N\right\rangle\right\rangle$ given by
\begin{equation}
\label{tf-ser}
T_\alpha^{\rm nc}(z)=D+\sum_{w\in\mathcal{F}_N\setminus\{\emptyset\}} \left(C\flat
A\sharp B\right)^wz^w=D+\sum_{k=0}^\infty C\left(\Delta(z)A\right)^k\Delta(z)B.
\end{equation}

\emph{The non-commutative Schur--Agler class}
$\mathcal{SA}_N^{\rm
nc}({\mathcal U},{\mathcal Y})$ consists of all FPSs
$f\in L({\mathcal U},{\mathcal Y})
\left\langle\left\langle z_1,\ldots, z_N\right\rangle\right\rangle$
such that for any separable Hilbert space ${\mathcal K}$ and
any $N$-tuple $\delta=(\delta_1,\ldots, \delta_N)$ of strict
contractions in ${\mathcal K}$ the limit in the operator norm topology
$$
f(\delta)=\lim_{m\rightarrow \infty}\sum_{w\in\mathcal{F}_N:\,
|w|\le m}f_w\otimes \delta^w
$$
exists and defines a contractive operator $f(\delta)\in
L({\mathcal U}\otimes{\mathcal K},{\mathcal Y}\otimes {\mathcal
K})$. We note that the non-commutative Schur--Agler class was
defined in \cite{BGM2} also for a more general class of operator $N$-tuples $\delta$.

Consider another  set of non-commuting indeterminates
$z'=(z_1',\ldots, z_N')$. For
$f(z)\in
L({\mathcal V},{\mathcal Y})
\left\langle\left\langle z_1,\ldots, z_N\right\rangle\right\rangle\quad {\rm
and}\quad f'(z')\in L({\mathcal V},{\mathcal
U})\left\langle\left\langle z_1',\ldots, z_N'\right\rangle\right\rangle$
we define a FPS
$$f(z){f'(z')}^*\in L({\mathcal U},{\mathcal
Y})\left\langle\left\langle z_1,\ldots, z_N,z_1',\ldots, z_N'
\right\rangle\right\rangle$$
by
\begin{equation}
\label{fps-prod} f(z){f'(z')}^*=\sum_{w,w'\in\mathcal{F}_N}
f_w {(f'_{w'})}^*
z^w{z'}^{{w'}^T}.
\end{equation}

In \cite{BGM2} the class $\mathcal{SA}_N^{\rm
nc}({\mathcal U},{\mathcal Y})$ was characterized as follows:
\begin{thm}\label{thm:BGM2} Let $f\in L({\mathcal U},{\mathcal
Y})\left\langle\left\langle z_1,\ldots, z_N\right\rangle\right\rangle$.
The following statements are equivalent:

$(1)$ $f\in\mathcal{SA}_N^{\rm nc}({\mathcal U},{\mathcal Y})$;

$(2)$ there exist auxiliary Hilbert spaces $\mathcal{H},{\mathcal
H}_1,\ldots, {\mathcal H}_N$ which are related by
 ${\mathcal
H}=\bigoplus_{k=1}^N{\mathcal H}_k$,  and a FPS $\varphi\in L({\mathcal H},{\mathcal Y})\left\langle \left\langle z_1,\ldots,z_N\right\rangle\right\rangle$ such that
\begin{equation}
\label{Ag} I_{\mathcal Y}-f(z){f(z')}^*=\varphi(z)(I_{\mathcal
H}-\Delta(z) \Delta(z')^*){\varphi(z')}^*;
\end{equation}

$(3)$ there exists a unitary GR-node
$\alpha=(N;A,B,C,D;{\mathcal X}=\bigoplus_{k=1}^N {\mathcal
X}_k,{\mathcal U}, {\mathcal Y})$ such that $f=T_\alpha^{\rm nc}$.
\end{thm}
We now give another characterization of the Schur--Agler class
$\mathcal{SA}_N^{\rm nc}({\mathcal U},{\mathcal Y})$.
\begin{thm}\label{thm:nc-schur}
A FPS $f$ belongs to $\mathcal{SA}_N^{\rm
nc}({\mathcal U},{\mathcal Y})$ if and only if for every
$n\in{\mathbb N}$ and $W\in({\mathbb D}^{n\times n})^N$ the 
limit in the operator norm topology
\begin{equation}
\label{sl}
f(W)=\lim_{m\rightarrow \infty}\sum_{w\in\mathcal{F}_N:\
|w|\le m}f_w\otimes W^w
\end{equation}
exists and $\|f(W)\|\le 1$.
\end{thm}
\begin{proof} The necessity is clear. We prove the sufficiency. We
set
$$f_k(z)=\sum_{w\in\mathcal{F}_N:\ |w|=k}f_wz^w,\quad k=0,1,\ldots\,.$$
Then for every $n\in{\mathbb N}$ and $W\in({\mathbb D}^{n\times
n})^N$, \eqref{sl} becomes
\begin{equation}\label{lim}
 f(W)=\lim_{m\rightarrow \infty}\sum_{k=0}^m
f_k(W),
\end{equation}
where the limit is taken in the operator norm topology.
Let $r\in (0,1)$ and choose $\tau>0$ such that $r+\tau<1$. Let
$W\in({\mathbb D}^{n\times n})^N$ be such that $\|W_j\|\le r$, $j=1,\ldots, N$. Then, for every $x\in{\mathcal U}\otimes{\mathbb
C}^n$ the series
$$
f\left(\frac{r+\tau}{r}\lambda W\right)x=\sum_{k=0}^\infty
\lambda^{k}f_k\left(\frac{r+\tau}{r}W\right)x$$
converges uniformly in $\lambda\in{\rm clos}\,{\mathbb D}$ to a ${\mathcal Y}\otimes{\mathbb C}^n$-valued
function holomorphic on ${\rm clos}\,{\mathbb D}$. Furthermore,
$$
\left\| f_k\left(\frac{r+\tau}{r}W\right)x\right\|=\left\|\frac{1}{2\pi i}\int_{\mathbb T}f\left(\frac{r+\tau}{r}\lambda W
\right)x\lambda^{-k-1}d\lambda\right\|\le\|x\|,$$
and therefore
\begin{equation}
\label{fk}
\|f_k(W)\|=\left\|f_k\left(\frac{r+\tau}{r}W\right)\left(\frac{r}{r+\tau}\right)^k\right\|\le
\left(\frac{r}{r+\tau}\right)^k.
\end{equation}
Thus we have
$$\left\|f(W)-\sum_{k=0}^mf_k(W)\right\|\le \sum_{k=m+1}^\infty
\|f_k(W)\|\le\sum_{k=m+1}^\infty\left(\frac{r}{r+\tau}\right)^k<\infty.
$$
We observe that the limit in \eqref{lim} is uniform in $n\in\mathbb{N}$ and $W\in (\mathbb{D}^{n\times n})^N$ such that $\| W_j\|\leq r,\ j=1,\ldots,N$.
Without loss of generality we may assume that in the definition of
the Schur--Agler class the space ${\mathcal K}$ is taken to be the
space $\ell_2$ of square summable sequences $s=(s_j)_{j=1}^\infty$
of complex numbers indexed by ${\mathbb N}$: $\sum_{j=1}^\infty
|s_j|^2<\infty$. We denote by $P_n$ the orthogonal projection from
$\ell_2$ onto the subspace of sequences for which $s_j=0$ for
$j>n$. This subspace is isomorphic to ${\mathbb C}^n$, and
thus for every $\delta=(\delta_1,\ldots, \delta_N)\in L(\ell_2)^N$ such that $\|\delta_j\|\le r$,
 $j=1,\ldots, N$, we may use \eqref{fk} and write
\begin{equation}
\| f_k(P_n\delta_1P_n,\ldots, P_n\delta_NP_n)\|\le
\left(\frac{r}{r+\tau}\right)^k.
\label{fkprime}
\end{equation}
Since the
sequence $P_n$ converges to $I_{\ell_2}$ in the strong operator topology (see, e.g.,
\cite{AkhGl}), and since strong limits of finite sums and products of operator sequences are equal to the corresponding sums and products of strong limits of these sequences, we obtain that
$$ s-\lim_{n\to\infty}f_k(P_n\delta_1P_n,\ldots,
P_n\delta_NP_n)=f_k(\delta).$$
Thus from \eqref{fkprime} we obtain
$$
\| f_k(\delta)\|\le
\left(\frac{r}{r+\tau}\right)^k.
$$
Therefore, the limit in the operator norm topology
$$
\label{slprime} f(\delta)=\lim_{m\rightarrow \infty}\sum_{k=0}^m
f_k(\delta)
$$
does exist, and
$$\| f(\delta)\|\le \sum_{k=0}^\infty\| f_k(\delta)\|\leq\sum_{k=0}^\infty\left(\frac{r}{r+\tau}\right)^k<\infty.$$
Moreover,
since the limit in \eqref{lim} is uniform  in $n\in\mathbb{N}$ and $W\in (\mathbb{D}^{n\times n})^N$ such that $\| W_j\|\leq r<1,\ j=1,\ldots,N,$ the
rearrangement of limits in the following chain of equalities is justified:
\[
\begin{split}
\lefteqn{\lim_{n\rightarrow\infty} f(P_n\delta_1P_n,\ldots,
P_n\delta_NP_n)h=\lim_{n\rightarrow\infty}\lim_{m\rightarrow\infty}
\sum_{k=0}^m f_k(P_n\delta_1P_n,\ldots, P_n\delta_NP_n)h}\\
&=\lim_{m\rightarrow\infty}\lim_{n\rightarrow\infty} \sum_{k=0}^m
f_k(P_n\delta_1P_n,\ldots, P_n\delta_NP_n)h
=\lim_{m\rightarrow\infty} \sum_{k=0}^m f_k(\delta )h=f(\delta)h
\end{split}
\]
(here $h$ is an arbitrary vector in $\mathcal{U}\otimes\ell_2$ and $\delta\in L(\ell_2)^N$ such that $\|\delta_j\|\le r$,
 $j=1,\ldots, N$).
Thus for every $\delta\in L(\ell_2)^N$ such that $\|\delta_j\|<1$, $j=1,\ldots, N$, we obtain
$\|f(\delta)\|\le 1$, i.e., $f\in\mathcal{SA}_N^{\rm nc}({\mathcal U},{\mathcal Y}).$
\end{proof}
\begin{rem}
One can see from the proof of Theorem~\ref{thm:nc-schur} that for arbitrary $f\in\mathcal{SA}_N^{\rm nc}(\mathcal{U,Y})$ and $r:\,0<r<1$, the series
$$f(\delta)=\sum_{k=0}^\infty f_k(\delta)$$
converges uniformly and absolutely in $\delta\in L(\mathcal{K})^N$  such that $\|\delta_j\|\le r$, $j=1,\ldots, N$, where $\mathcal{K}$ is any separable Hilbert space.
\end{rem}
\begin{cor}
A matrix-inner in ${\mathcal D}_N$ rational FPS
$f$ belongs to the class $\mathcal{SA}_N^{\rm nc}({\mathbb C}^q)=\mathcal{SA}_N^{\rm nc}({\mathbb C}^q,{\mathbb C}^q)$.
\end{cor}
Thus, for the case $J=I_q$, Theorem~\ref{thm:consprime} establishes the
existence of a unitary GR-realization for an
arbitrary matrix-inner rational FPS, i.e., recovers Theorem~\ref{thm:BGM2} for the case of a \emph{matrix-inner rational} FPS. However, it says even more than Theorem~\ref{thm:BGM2} in this case, namely that such a unitary realization can be found minimal, thus
finite-dimensional, and that this minimal unitary realization is unique up to a unitary similarity. The
representation \eqref{Ag} with the rational FPS
$\varphi\in {\mathbb C}^{q\times \gamma}\left\langle\left\langle z_1,\ldots,
z_N\right\rangle\right\rangle_{\rm rat}$
given by
$$\varphi(z)=C(I_\gamma-\Delta(z)A)^{-1}$$
is obtained from \eqref{s-id'} by making use of Corollary~\ref{cor:vanish-fps}.

\section{Matrix-selfadjoint rational formal power series}\label{sec:sa}
\subsection{A multivariable non-commutative analogue of the line case}\label{sec:line-sa}
A rational FPS $\Phi\in
{\mathbb C}^{q\times q}\left\langle\left\langle z_1,\ldots, z_N
\right\rangle\right\rangle_{\rm rat}$ will be called
\emph{matrix-selfadjoint on ${\mathcal J}_N$} if for every $n\in{\mathbb N}$:
$$\Phi(Z)=\Phi(Z)^*$$
at all points $Z\in\left(i{\mathbb H}^{n\times n}\right)^N$ where it is
defined.

The following theorem is a multivariable non-commutative  counterpart
of Theorem 4.1 from \cite{AG1} which was originally proved in \cite{GLR1}.
\begin{thm}
\label{thm:sa}
Let $\Phi\in
{\mathbb C}^{q\times q}
\left\langle\left\langle z_1,\ldots, z_N\right\rangle\right\rangle_{\rm rat}$,  and let $\alpha$
be a minimal GR-realization of $\Phi$ of the form \eqref{min}. Then $\Phi$ is
matrix-selfadjoint on ${\mathcal J}_N$ if and only if the following
conditions hold:

$(a)$ the matrix $D$ is Hermitian, that is, $D=D^*$;

$(b)$ there exists an invertible Hermitian matrix $H={\rm diag}(H_1,\ldots
,H_N)$ with $H_k\in{\mathbb C}^{\gamma_k\times \gamma_k},\ k=1,\ldots,N,$ and such that
\begin{eqnarray}
\label{ls}
A^*H+HA&=&0,\\
\label{cs}
C&=&iB^*H.
\end{eqnarray}
\end{thm}
\begin{proof} We first observe that $\Phi$ is matrix-selfadjoint
on ${\mathcal J}_N$ if and only if the FPS
$F\in{\mathbb C}^{2q\times 2q}\left\langle\left\langle z_1,\ldots, z_N
\right\rangle\right\rangle_{\rm rat}$ given by
\begin{equation}
F(z)=\begin{pmatrix}I_q&i\Phi(z)\\ 0&I_q\end{pmatrix}
\label{F-def}
\end{equation}
is matrix-$J_1$-unitary on ${\mathcal J}_N$, where
\begin{equation}
J_1=\begin{pmatrix}0&I_q\\ I_q&0\end{pmatrix}.
\label{j1}
\end{equation}
Moreover, $F$ admits the GR-realization
$$
\beta=(N;A,\begin{pmatrix}0&B\end{pmatrix},\begin{pmatrix}iC\\0\end{pmatrix},
\begin{pmatrix}I_q&iD\\ 0&I_q\end{pmatrix};{\mathbb C}^\gamma=
\bigoplus_{k=1}^N{\mathbb C}^{\gamma_k},{\mathbb
C}^{2q}).$$
This realization is minimal. Indeed, the $k$-th truncated
observability (resp., controllability) matrix of $\beta$ is equal
to
\begin{equation}
\label{obs}
\widetilde{{\mathcal O}_k}(\beta)=\begin{pmatrix}
i\widetilde{{\mathcal O}_k}(\alpha)\\0\end{pmatrix}
\end{equation}
and, resp.,
\begin{equation}
\label{contr}
\widetilde{{\mathcal C}_k}(\beta)=\begin{pmatrix}
0&\widetilde{{\mathcal C}_k}(\alpha)\end{pmatrix},
\end{equation}
and therefore has full column (resp., row) rank. Using Theorem~\ref{thm:Lyap} of the present paper we see that $\Phi$ is matrix-selfadjoint on
${\mathcal J}_N$ if and only if:

$(1)$ the matrix {\scriptsize $\begin{pmatrix}I_q&iD\\ 0&I_q\end{pmatrix}$} is
$J_1$-unitary;

$(2)$ there exists an invertible Hermitian matrix $H={\rm diag}(H_1,\ldots,H_N)$, with $H_k\in\mathbb{C}^{\gamma_k\times\gamma_k},\ k=1,\ldots,N$,
such that
\begin{eqnarray*}
A^*H+HA&=&-\begin{pmatrix}iC\\0\end{pmatrix}^*J_1\begin{pmatrix}iC\\0\end{pmatrix},\\
\begin{pmatrix}0&B\end{pmatrix}&=&-H^{-1}\begin{pmatrix}iC\\0\end{pmatrix}^*J_1
\begin{pmatrix}I_q&iD\\ 0&I_q\end{pmatrix}.
\end{eqnarray*}
These conditions are in turn readily seen to be equivalent to
conditions $(a)$ and $(b)$ in the statement of the theorem.
\end{proof}

From Theorem \ref{thm:Lyap} it follows that the matrix $H={\rm
diag}(H_1,\ldots, H_N)$ is uniquely determined by the given
minimal GR-realization of $\Phi$. In a similar way
as in Section~\ref{sec:line}, it can be shown that $H_k,\ k=1,\ldots,N$, are given by the
formulas
\[
\begin{split}
H_k&=-\left({\rm col}_{w\in{\mathbb
F}_N:\,|w|\le q \gamma-1}\left(B^*\flat(-A^*)\right)^{w
g_k}\right)^+ \left({\rm col}_{w\in{\mathbb
F}_N:\,|w|\le q \gamma-1}\left(C\flat A\right)^{w
g_k}\right)\\
&=\left({\rm row}_{w\in\mathcal{F}_N:\,|w|\le q
\gamma-1}\left((-A^*)\sharp C^*\right)^{g_kw^T}\right)
\left({\rm row}_{w\in\mathcal{F}_N:\,|w|\le q
\gamma-1}\left(A\sharp B\right)^{g_kw^T}\right)^\dag .
\end{split}
\]
The matrix $H={\rm diag}(H_1,\ldots,H_N)$ is called in this case \emph{the associated structured Hermitian matrix} (associated with a minimal GR-realization of the FPS $\Phi$).

It follows from \eqref{ls} and \eqref{cs} that for $n\in{\mathbb
N}$ and $Z,Z'\in \left(i{\mathbb H}^{n\times n}\right)^N$ we have:
\begin{eqnarray}
\lefteqn{
\Phi(Z)-\Phi(Z')^*=i(C\otimes I_n)\left(I_\gamma\otimes
I_n-\Delta(Z)(A\otimes
I_n)\right)^{-1}}\label{m}\\
&\times\Delta(Z+{Z'}^*)\left(H^{-1}\otimes
I_n\right)\left(I_\gamma\otimes I_n-(A^*\otimes
I_n)\Delta({Z'}^*)\right)^{-1}(C^*\otimes I_n),\nonumber \\ \lefteqn{
\Phi(Z)-\Phi(Z')^*=i(B^*\otimes I_n)\left(I_\gamma\otimes
I_n-\Delta({Z'}^*)(A^*\otimes
I_n)\right)^{-1}}\label{mprime}\\
& \times  \Delta(Z+{Z'}^*)\left(H\otimes
I_n\right)\left(I_\gamma\otimes I_n-(A\otimes
I_n)\Delta(Z)\right)^{-1}(B\otimes I_n). \nonumber
\end{eqnarray}
Note that if $A,B$ and $C$ are matrices which satisfy \eqref{ls}
and \eqref{cs} for some (not necessarily invertible) Hermitian
matrix $H$, and if $D$ is Hermitian, then
$$\Phi(z)=D+C(I-\Delta(z)A)^{-1}\Delta(z)B$$
is a rational FPS which is matrix-selfadjoint on
${\mathcal J}_N$. This follows from the fact that \eqref{mprime}
is still valid in this case (the corresponding GR-realization of $\Phi$ is, in general, not minimal).

If $A,B$ and $C$ satisfy the equalities
\begin{eqnarray}
\label{lsprime}
GA^*+AG&=&0,\\
B&=&iGC^*
\label{bsprime}
\end{eqnarray}
for some (not necessarily invertible) Hermitian matrix $G={\rm
diag}(G_1,\ldots, G_N)$ then  \eqref{m} is valid with
$H^{-1}$ replaced by $G$ (the diagonal structures of $G$,
$\Delta(Z)$ and $\Delta(Z')$ are compatible), and hence $\Phi$ is
matrix-selfadjoint on $\mathcal{J}_N$.

As in Section~\ref{sec:line}, we can solve inverse problems using Theorem~\ref{thm:sa}. The proofs are easy and omitted.
\begin{thm}
Let $(C,A)$ be an observable pair of matrices, in the sense that
$\mathcal{O}_k$ has a full column rank for all
$k\in\left\{1,\ldots, N\right\}$. Then there exists a rational
FPS which is matrix-selfadjoint on  ${\mathcal
J}_N$ with a minimal GR-realization $\alpha$ of the form \eqref{min} if and only if the equation
$$A^*H+HA=0$$
has a solution $H={\rm diag}(H_1,\ldots, H_N)$ (with
$H_k\in{\mathbb C}^{\gamma_k\times\gamma_k},\ k=1,\ldots,N$) which is both
Hermitian and invertible. When such a solution exists, $D$ can be
any Hermitian matrix and $B=iH^{-1}C^*$.
\end{thm}
\begin{thm}
Let $(A,B)$ be a controllable pair of matrices, in the sense that $
{\mathcal C}_k$ has a full row rank for all $k\in\{1,\ldots, N\}$. Then there exists
a rational FPS which is matrix-selfadjoint on ${\mathcal J}_N$ with a
minimal GR-realization
$\alpha$ of the form \eqref{min} if and only if the equation
$$GA^*+AG=0$$
has a solution $G={\rm diag}(G_1,\ldots, G_N)$ (with $G_k\in{\mathbb C}^{
\gamma_k\times\gamma_k}$, $k=1,\ldots, N$) which is both Hermitian and
invertible. When such a solution exists, $D$ can be any Hermitian
matrix and $C=iB^*G^{-1}$.
\end{thm}
From \eqref{obs} and \eqref{contr} obtained in Theorem \ref{thm:sa}, and from Theorem~\ref{thm:c=o} we obtain the
following result:
\begin{thm}\label{thm:cosa}
Let $\Phi$ be a matrix-selfadjoint on ${\mathcal J}_N$ rational
FPS with a  GR-realization
$\alpha$ of the form \eqref{qq}.
Let $H={\rm diag}(H_1,\ldots, H_N)$ (with $H_k\in{\mathbb C}^{r_k\times r_k}$, $k=1,\ldots, N$) be both Hermitian and
invertible and satisfy \eqref{ls} and \eqref{cs}. Then the GR-node $\alpha$ is
observable if and only if it is controllable.
\end{thm}
The following Lemma is an analogue of Lemma~\ref{lem:h}. It is easily
proved by applying  Lemma \ref{lem:h} to the  matrix-$J_1$-unitary on ${\mathcal J}_N$ function $F$ defined in
\eqref{F-def}.
\begin{lem}\label{lem:h-cosa}
Let $\Phi\in{\mathbb C}^{q\times q}\left\langle\left\langle z_1,\ldots, z_N
\right\rangle\right\rangle_{\rm rat}$ be matrix-selfadjoint on ${\mathcal J}_N$, and
let
$\alpha^{(i)}=(N;A^{(i)},B^{(i)},C^{(i)},D;{\mathbb C}^{\gamma}=
\bigoplus_{k=1}^N{\mathbb C}^{\gamma_k},{\mathbb C}^q)$
be two minimal GR-realizations of $\Phi$, with the associated structured Hermitian matrices
$H^{(i)}={\rm diag}(H_1^{(i)},\ldots, H_N^{(i)})$, $i=1,2$.
Then these two realizations and associated matrices $H^{(i)}$ are linked by
\eqref{sim} and \eqref{h-sim}. In particular, for each $k\in\{ 1,\ldots, N\}$ the matrices $H_k^{(1)}$ and $H_k^{(2)}$ have the
same signature.
\end{lem}

For $n\in{\mathbb N}$, points $Z,Z'\in({\mathbb C}^{n\times n})^N$ where
$\Phi(Z)$ and $\Phi(Z')$ are well-defined, $F$ given by
\eqref{F-def}, and $J_1$ defined by \eqref{j1} we have:
\begin{equation}
\label{u-sa}
J_1\otimes I_n-F(Z)(J_1\otimes I_n){F(Z')}^*=\begin{pmatrix}
  \frac{\Phi(Z)-{\Phi(Z')}^*}{i}&0\\0&0\end{pmatrix}
\end{equation}
and
\begin{equation}
\label{u-saprime}
J_1\otimes I_n-{F(Z')}^*(J_1\otimes I_n)F(Z)=\begin{pmatrix}
0&0\\0&  \frac{\Phi(Z)-{\Phi(Z')}^*}{i}\end{pmatrix}.
\end{equation}
Combining these equalities with \eqref{m} and \eqref{mprime} and using Corollary~\ref{cor:vanish-fps} we obtain the
following analogue of Theorem~\ref{thm:neg}.
\begin{thm}
\label{thm:neg-sa}
Let $\Phi$ be a matrix-selfadjoint on $\mathcal{J}_N$ rational FPS, and let
$\alpha$
be its minimal GR-realization of the form \eqref{min}, with the associated structured Hermitian matrix
$H={\rm diag}(H_1,\ldots, H_N)$. Then for each $k\in\left\{1,\ldots,
N\right\}$ the number of negative eigenvalues of the matrix $H_k$ is equal to
the number of negative squares of the kernels
\begin{equation}
\label{ker-sa}
\begin{split}
K^{\Phi,k}_{w,w^\prime}
&=(C\flat A)^{wg_k}H_k^{-1}(A^*\sharp C^*)^{g_kw^{\prime T}}\\
K^{\Phi^*,k}_{w,w^\prime}
   &=(B^*\flat A^*)^{wg_k}H_k(A\sharp B)^{g_kw^{\prime T}},
\end{split}\quad w,w^\prime
\in\mathcal{F}_N.
\end{equation}
Finally, for $k\in\left\{1,\ldots, N\right\}$, let ${\mathcal K}_k(\Phi)$
(resp., ${\mathcal K}_k(\Phi^*)$) denote the span of the functions $w\mapsto
K_{w,w^\prime}^{\Phi,k}$ (resp.,
$w\mapsto K_{w,w^\prime}^{\Phi^*,k}$) where $w^\prime\in\mathcal{F}_N$ and
$c\in{\mathbb C}^q$. Then,
$${\rm dim}~{\mathcal K}_k(\Phi)=
{\rm dim}~{\mathcal K}_k(\Phi^*)=\gamma_k.
$$
\end{thm}
Let $\Phi_1$ and $\Phi_2$ be two FPSs from
${\mathbb C}^{q\times q}\left\langle\left\langle z_1,\ldots, z_N
\right\rangle\right\rangle_{\rm rat}$. The \emph{additive decomposition}
$$\Phi=\Phi_1+\Phi_2$$
is called \emph{minimal} if
$$\gamma_k(\Phi)=\gamma_k(\Phi_1)+\gamma_k(\Phi_2),\quad k=1,\ldots, N,$$
where $\gamma_k(\Phi), \gamma_k(\Phi_1)$ and $\gamma_k(\Phi_2)$ denote the
dimensions of the $k$-th component of the state space of a minimal
GR-realization of $\Phi,\Phi_1$ and $\Phi_2$, respectively.
The following theorem is an analogue of Theorem \ref{thm:factor}.
\begin{thm}
\label{thm:sa-decomp}
Let $\Phi_i,\ i=1,2$, be  matrix-selfadjoint on $\mathcal{J}_N$ rational
FPSs, with minimal GR-realizations $\alpha^{(i)}=(N;A^{(i)},B^{(i)},C^{(i)},D^{(i)};{\mathbb C}^{\gamma^{(i)}}=
\bigoplus_{k=1}^N{\mathbb C}^{\gamma_k^{(i)}},{\mathbb C}^q)$
and the associated structured Hermitian matrices
$H^{(i)}={\rm diag}(H_1^{(i)},\ldots, H_N^{(i)})$. Assume that the additive decomposition
$\Phi=\Phi_1+\Phi_2$ is minimal. Then the GR-node
$\alpha=(N;A,B,C,D;{\mathbb C}^\gamma=\bigoplus_{k=1}^N
{\mathbb C}^{\gamma_k},{\mathbb C}^q)$ defined by
$$D=D^{(1)}+D^{(2)},\qquad \gamma_k=\gamma_k^{(1)}+\gamma_k^{(2)},\quad
k=1,\ldots, N,$$
and with respect to the decomposition ${\mathbb C}^\gamma=
{\mathbb C}^{\gamma^{(1)}}\oplus{\mathbb C}^{\gamma^{(2)}}$,
\begin{equation}
A=\begin{pmatrix}A^{(1)}&0\\
0&A^{(2)}\end{pmatrix},\quad
B=\begin{pmatrix}B^{(1)}\\
B^{(2)}\end{pmatrix},\quad
C=\begin{pmatrix}C^{(1)}&
C^{(2)}\end{pmatrix},
\end{equation}
is a minimal GR-realization of $\Phi$, with the associated structured Hermitian matrix
$H={\rm diag}(H_1,\ldots, H_N)$ such that for each $k\in\left\{1,\ldots, N
\right\}$:
$$H_k=\begin{pmatrix}H_k^{(1)}&0\\0&H_k^{(2)}\end{pmatrix}.
$$
\end{thm}
Let $\nu_k(\Phi)$ denote the number of negative squares of either of the
functions defined in \eqref{ker-sa}. In view of Theorem \ref{thm:neg-sa} and
Theorem \ref{thm:sa} these numbers are uniquely determined by $\Phi$.
\begin{cor}
Let $\Phi_1$ and $\Phi_2$ be matrix-selfadjoint on
${\mathcal J}_N$ rational FPSs, and assume that the additive decomposition
$\Phi=\Phi_1+\Phi_2$ is minimal. Then
$$
\nu_k(\Phi)=\nu_k(\Phi_1)+\nu_k(\Phi_2),\quad k=1,2,\ldots, N.$$
\end{cor}
An additive decomposition of a matrix-selfadjoint  on
${\mathcal J}_N$ rational FPS $\Phi$ is called a \emph{minimal matrix-selfadjoint decomposition} if it is minimal and both $\Phi_1$ and $\Phi_2$ are matrix-selfadjoint   on
${\mathcal J}_N$ rational FPSs. The set of all minimal matrix-selfadjoint decompositions of a
matrix-selfadjoint on ${\mathcal J}_N$ rational FPS is given by the
following theorem, which is a multivariable non-commutative counterpart of
\cite[Theorem 4.6]{AG1}. The proof uses Theorem~\ref{thm:min-fact} applied to
the FPS $F$ defined by \eqref{F-def}, and follows the same argument as one in the proof of Theorem~4.6 in \cite{AG1}.
\begin{thm}
Let $\Phi$ be a matrix-selfadjoint on ${\mathcal J}_N$ rational
FPS, with a minimal GR-realization
$\alpha$ of the form \eqref{min}
and the associated structured Hermitian matrix $H={\rm diag}(H_1,\ldots, H_N)$. Let
$\mathcal{M}=\bigoplus_{k=1}^N\mathcal{M}_k$ be an $A$-invariant subspace, with $\mathcal{M}_k\subset
{\mathbb C}^{\gamma_k},\ k=1,\ldots,N$, and assume that $\mathcal{M}$ is non-degenerate in the
associated inner product $[\,\cdot\, ,\,\cdot\,]_H$. Let
$\Pi={\rm diag}(\Pi_1,\ldots, \Pi_N)$ be the projection defined by
$$\ker\Pi=\mathcal{M},\qquad {\rm ran}\,\Pi=\mathcal{M}^{[\perp]},$$
that is,
$$\ker\Pi_k=\mathcal{M}_k,\qquad {\rm ran}\,\Pi_k=\mathcal{M}_k^{[\perp]},\quad k=1,\ldots, N.$$
Let $D=D_1+D_2$ be a decomposition of $D$ into two Hermitian matrices.
Then the decomposition $\Phi=\Phi_1+\Phi_2$, where
\[
\begin{split}
\Phi_1(z)&=D_1+C(I_\gamma-\Delta(z)A)^{-1}\Delta(z)(I_\gamma-\Pi)B,\\
\Phi_2(z)&=D_2+C\Pi(I_\gamma-\Delta(z)A)^{-1}\Delta(z)B,
\end{split}
\]
is a minimal matrix-selfadjoint decomposition of $\Phi$.

Conversely, any minimal matrix-selfadjoint decomposition of $\Phi$ can be
obtained in such a way, and with a fixed decomposition $D=D_1+D_2$,
the correspondence between
minimal matrix-selfadjoint decompositions of $\Phi$ and non-degenerate
$A$-invariant subspaces of the form
$\mathcal{M}=\bigoplus_{k=1}^N \mathcal{M}_k$,
where $\mathcal{M}_k\subset {\mathbb C}^{\gamma_k}$, $k=1,\ldots, N$, is
one-to-one.
\end{thm}
\begin{rem}
Minimal matrix-selfadjoint decompositions do not always exist, even in the
case $N=1$. For counterexamples see \cite{AG1}.
\end{rem}
\subsection{A multivariable non-commutative analogue of the circle case}\label{sec:circle-sa}
In this subsection we briefly review some analogues of the theorems presented
in Section~\ref{sec:line-sa}.
\begin{thm}
Let $\Psi$ be a rational FPS and
$\alpha$ be its minimal GR-realization of the form \eqref{min}. Then $\Psi$ is
matrix-selfadjoint on ${\mathcal T}_N$ (that is, for all $n\in\mathbb{N}$ one has $\Psi(Z)=\Psi(Z)^*$ at all
points $Z\in(\mathbb{T}^{n\times n})^N$ where $\Psi$ is defined) if and only if there
exists an invertible Hermitian matrix $H={\rm diag}(H_1,\ldots, H_N)$, with
$H_k\in{\mathbb C}^{\gamma_k\times \gamma_k}$, $k=1,\ldots, N$, such that
\begin{equation}
A^*HA=H,\quad D-D^*=iB^*HB,\quad C=iB^*HA.
\label{saprime}
\end{equation}
\end{thm}
\begin{proof}
Consider the FPS $f\in
{\mathbb C}^{2q\times 2q}\left\langle\left\langle z_1,\ldots, z_N
\right\rangle\right\rangle_{\rm rat}$ defined by
\begin{equation}
f(z)=\begin{pmatrix}I_q&i\Psi(z)\\ 0&I_q\end{pmatrix}.
\label{f-defprime}
\end{equation}
Using Theorem~\ref{thm:j-u}, we see that $f$ is matrix-$J_1$-unitary
on ${\mathcal T}_N$, with
\begin{equation}
\label{j1prime}
J_1=\begin{pmatrix}0&I_q\\ I_q&0\end{pmatrix},
\end{equation}
if and only if its GR-realization
$$\beta=(N;A,\begin{pmatrix}0&B\end{pmatrix},
\begin{pmatrix}iC\\0\end{pmatrix},
\begin{pmatrix}I_q&iD\\0&I_q\end{pmatrix};
{\mathbb C}^\gamma=\oplus_{j=1}^N{\mathbb C}^{\gamma_j},
{\mathbb C}^{2q})$$
(which turns out to be minimal, as can be shown in the same way as in Theorem~\ref{thm:sa}) satisfies the following
condition: there exists an Hermitian invertible matrix
$H={\rm diag}(H_1,\ldots, H_N)$, with $H_k\in{\mathbb C}^{
\gamma_k\times\gamma_k}$, $k=1,\ldots, N$, such that
$$
\begin{pmatrix}A&0&B\\
iC&I_q&iD\\
0&0&I_q\end{pmatrix}^*
\begin{pmatrix}
H&0&0\\
0&0&I_q\\
0&I_q&0\end{pmatrix}
\begin{pmatrix}A&0&B\\
iC&I_q&iD\\
0&0&I_q\end{pmatrix}
=\begin{pmatrix}
H&0&0\\
0&0&I_q\\
0&I_q&0\end{pmatrix},
$$
which is equivalent to the condition stated in the theorem.
\end{proof}
For a given minimal GR-realization of $\Psi$ the matrix $H$ is unique,
as follows from Theorem~\ref{thm:Cayley}. It is called
\emph{the associated structured Hermitian matrix} of $\Psi$.

The set of all minimal matrix-selfadjoint additive decompositions of a given
matrix-selfadjoint on ${\mathcal T}_N$ rational FPS is described by the following
theorem, which is a multivariable non-commutative counterpart of
\cite[Theorem 5.2]{AG1}, and is proved by applying Theorem~\ref{thm:factor-circle}
to the matrix-$J_1$-unitary
on ${\mathcal T}_N$ FPS $f$ defined by \eqref{f-defprime}, where $J_1$ is defined by \eqref{j1prime}. (We omit
the proof.)
\begin{thm} Let $\Psi$ be a matrix-selfadjoint on ${\mathcal T}_N$ rational
FPS and
$\alpha$ be its minimal GR-realization of the form \eqref{min}, with the associated structured Hermitian matrix
$H={\rm diag}(H_1,\ldots, H_N)$. Let
$\mathcal{M}=\bigoplus_{k=1}^N\mathcal{M}_k$
be an
$A$-invariant subspace, with $\mathcal{M}_k\subset {\mathbb C}^{\gamma_k}$, $k=1,\ldots, N$,  and assume that $\mathcal{M}$ is non-degenerate
in the associated inner product $[\cdot,\cdot]_H$. Let
$\Pi={\rm diag}(\Pi_1,\ldots, \Pi_N)$ be the projection defined by
$$\ker\Pi=\mathcal{M},\qquad {\rm ran}\,\Pi=\mathcal{M}^{[\perp]},$$
that is,
$$\ker\Pi_k=\mathcal{M}_k,\qquad {\rm ran}\,\Pi_k=\mathcal{M}_k^{[\perp]},\quad k=1,\ldots, N.$$
Then the decomposition $\Psi=\Psi_1+\Psi_2$, where
\[
\begin{split}
\Psi_1(z)&=D_1+C(I_\gamma-\Delta(z)A)^{-1}\Delta(z)(I_\gamma-\Pi)B,\\
\Psi_2(z)&=D_2+C\Pi(I_\gamma-\Delta(z)A)^{-1}\Delta(z)B,
\end{split}
\]
with $D_1=\frac{i}{2}B_1^*H^{(1)}B_1+S$, the matrix $S$ being an arbitrary
Hermitian matrix, and
$$B_1=P_\mathcal{M}B,\qquad H^{(1)}=P_\mathcal{M}H\big|_\mathcal{M},$$
is a minimal matrix-selfadjoint additive decomposition of $\Psi$ (here
$P_\mathcal{M}$ denotes the orthogonal projection onto $\mathcal{M}$ in the standard
metric of ${\mathbb C}^\gamma$).

Conversely, any minimal matrix-selfadjoint additive decomposition of $\Psi$ is
obtained in such a way, and for a fixed $S$, the correspondence
between minimal matrix-selfadjoint additive decompositions of $\Psi$ and
non-degenerate $A$-invariant subspaces of the form $\mathcal{M}=\bigoplus_{k=1}^N \mathcal{M}_k$,
where $\mathcal{M}_k\subset {\mathbb C}^{\gamma_k}$, $k=1,\ldots, N$, is
one-to-one.
\end{thm}
\section{Finite-dimensional de Branges--Rovnyak spa\-ces and backward shift realizations:
The multivariable non-commutative setting}\label{sec:spec}
In this section we describe certain model realizations of matrix-$J$-unitary rational FPSs. We restrict ourselves to the case of
FPSs which are matrix-$J$-unitary on ${\mathcal J}_N$.
Analogous realizations can be constructed for rational FPSs
which are matrix-$J$-unitary on ${\mathcal T}_N$ or matrix-selfadjoint
either on ${\mathcal J}_N$ or ${\mathcal T}_N$.
\subsection{Non-commutative formal reproducing kernel
Pontryagin spaces}
Let $F$ be a matrix-$J$-unitary on ${\mathcal J}_N$ rational FPS and $\alpha$ be its minimal GR-realization of the form \eqref{min},
with the associated structured Hermitian matrix $H={\rm diag}(H_1,\ldots, H_N)$. Then by
Theorem \ref{thm:neg}, for each $k\in\left\{1,\ldots, N\right\}$ the kernel \eqref{kerns}
has the number $\nu_k(F)$ of negative eigenvalues equal to the number of negative squares of $H_k$. Lemma \ref{lem:h} implies that the
kernel $K_{w,w'}^{F,k}$ from \eqref{kerns} does not depend on the choice of a
minimal realization of $F$.
Theorem \ref{thm:neg} also asserts that the span of the functions
$$
w\mapsto K_{w,w'}^{F,k}c,\quad{\rm where}\quad
w'\in\mathcal{F}_N\quad{\rm and}\quad c\in{\mathbb C}^q,
$$
is the space ${\mathcal K}_k(F)$ with $\dim{\mathcal K}_k(F)=
\gamma_k$, $k=1,\ldots, N$. One can introduce a new metric on
each of the spaces ${\mathcal K}_k(F)$ as follows. First, define an Hermitian
form $[\,\cdot\, ,\,\cdot\,]_{F,k}$ by:
$$
[K_{\cdot,w'}^{F,k}c',K_{\cdot,w}^{F,k}c]_{F,k}=c^*K^{F,k}_{
w,w'}c'.
$$
This form is easily seen to be well defined on the whole space
${\mathcal K}_k(F)$, that is, if $f$ and $h$ belong to ${\mathcal K}_k(F)$
and
$$
f_w=\sum_jK_{w,w_j}^{F,k}c_j=\sum_\ell
K_{w,w'_\ell}^{F,k}
c'_\ell$$
and
$$
h_w=\sum_s K_{w,v_s}^{F,k}d_s=\sum_t
K_{w,v'_t}^{F,k}
d'_t,$$
where all the sums are finite, then
\[
[f,h]_{F,k}=\left[\sum_j K_{\cdot,w_j}^{F,k}c_j,
\sum_s K_{\cdot,v_s}^{F,k}d_s\right]_{F,k}=\left[\sum_\ell K_{\cdot,w'_\ell}^{F,k}c'_\ell,
\sum_t K_{\cdot,v'_t}^{F,k}d'_t\right]_{F,k}.
\]
Thus, the space ${\mathcal K}_k(F)$ endowed with this new (indefinite) metric
is a finite dimensional reproducing kernel Pontryagin space (RKPS) of functions
 on $\mathcal{F}_N$ with the
reproducing kernel $K_{w,w'}^{F,k}$. We refer to
\cite{Sor,AD1,ADRS} for more information on the theory of
reproducing kernel Pontryagin spaces.
In a similar way, the space ${\mathcal K}(F)=\bigoplus_{k=1}^N{\mathcal K}_k(F)$
endowed with the indefinite inner product
$$[f,h]_F=\sum_{k=1}^N[f_k,h_k]_{F,k}.$$
where $f={\rm col}\,(f_1,\ldots, f_N)$ and $h={\rm col}\,(h_1,\ldots, h_N)$,
becomes a reproducing kernel Pontryagin space with the reproducing kernel
$$K_{w,w'}^F={\rm diag}(K_{w,w'}^{F,1},\ldots,
K_{w,w'}^{F,N}),\quad w,w'\in\mathcal{F}^N.$$

Rather than the kernels $K_{w,w'}^{F,k}$, $k=1,\ldots N$,
and $K_{w,w'}^{F}$ we prefer to use the FPS
 kernels
\begin{eqnarray}
\label{kfk}
K^{F,k}(z,z')&=&\sum_{w,w'\in\mathcal{F}^N}
K_{w,w'}^{F,k}z^w {z^\prime}^{{w'}^T},\quad
k=1,\ldots, N,\\
\label{kf}
K^{F}(z,z')&=&\sum_{w,w'\in\mathcal{F}^N}
K_{w,w'}^{F}z^w {z^\prime}^{{w'}^ T},
\end{eqnarray}
and instead of the reproducing kenrel Pontryagin spaces
${\mathcal K}_k(F)$ and ${\mathcal K}(F)$ we will use the notion of
\emph{non-commutative formal reproducing kernel Pontryagin spaces}
(\emph{NFRKPS} for short; we will use the same notations for these spaces) which we
introduce below in a way analogous  to the way J.~A.~Ball and
V.~Vinnikov introduce
non-commutative formal reproducing kernel Hilbert spaces
(NFRKHS for short) in \cite{BV}.

Consider a FPS $$K(z,z')=\sum_{w,w'\in\mathcal{F}_N}
K_{w,w'} z^w {z^{\prime}}^{{w'}^T}\in
L({\mathcal C})\left\langle\left\langle z_1,\ldots,z_N,z'_1,\ldots,z'_N\right\rangle\right\rangle_{\rm rat},$$
where ${\mathcal C}$ is a Hilbert space. Suppose that
$$K(z',z)={K(z,z')}^*=
\sum_{w,w'\in\mathcal{F}_N} K^*_{w,w'}
z^{\prime{w'}}z^{w^T}.$$
Then $K^*_{w,w'}=K_{w',w}$ for
all $w,w'\in\mathcal{F}_N$. Let $\kappa\in{\mathbb N}$.
We will say that the \emph{FPS} $K(z,z')$ is a \emph{kernel
with $\kappa$ negative squares} if
$K_{w,w'}$ is a kernel on $\mathcal{F}_N$ with $\kappa$ negative squares, i.e. for every integer $\ell$ and
every choice of $w_1,\ldots, w_\ell\in\mathcal{F}_N$ and
$c_1,\ldots, c_\ell\in{\mathcal C}$  the $\ell\times \ell$ Hermitian matrix
with $(i,j)$-th entry equal to $c_i^*K_{w_i,w_j}c_j$
has at most $\kappa$ strictly negative eigenvalues, and exactly $\kappa$ such
eigenvalues for some choice of $\ell,w_1,\ldots, w_\ell,
c_1,\ldots, c_\ell$.

Define  on the space ${\mathcal G}$ of finite sums of
FPSs of the form $$K_{w'}(z)c=\sum_{w\in\mathcal{F}_N}K_{w,w'}z^wc,$$ where $w'
\in\mathcal{F}_N$ and $c\in{\mathcal C}$, the inner product as follows:
$$
\left[\sum_{i}K_{w_i}(z)c_i,\sum_{j}K_{w'_j}(z)c'_j\right]_\mathcal{G}=
\sum_{i,j}\langle K_{w'_j,w_i}c_i,c'_j\rangle_{{\mathcal C}}.$$
It is easily seen to be well defined. The space ${\mathcal G}$ endowed with
this inner product can be completed in a unique way to a Pontryagin space
${\mathcal P}(K)$ of FPSs, and in ${\mathcal P}(K)$ the
reproducing kernel property is
\begin{equation}
[f,K_w(\cdot)c]_{{\mathcal P}(K)}=\langle f_w,c\rangle_{\mathcal C}.
\label{repr}
\end{equation}
See \cite[Theorem 6.4]{AD1} for more details on such completions.

Define  the pairings $[\cdot,\cdot]_{\mathcal{P}(K)\times\mathcal{P}(K)\left\langle \left\langle z_1,\ldots,z_N\right\rangle\right\rangle}$ and $\left\langle \cdot,\cdot\right\rangle_{\mathcal{C}\left\langle \left\langle z_1,\ldots,z_N\right\rangle\right\rangle\times\mathcal{C}}$ as mappings
$\mathcal{P}(K)\times\mathcal{P}(K)\left\langle \left\langle z_1,\ldots,z_N\right\rangle\right\rangle\rightarrow\mathbb{C}\left\langle \left\langle z_1,\ldots,z_N\right\rangle\right\rangle$ and $\mathcal{C}\left\langle \left\langle z_1,\ldots,z_N\right\rangle\right\rangle\times\mathcal{C}\rightarrow\mathbb{C}\left\langle \left\langle z_1,\ldots,z_N\right\rangle\right\rangle$
by
\begin{eqnarray*}
\left[f,\sum_{w\in\mathcal{F}_N}g_wz^w\right]_{\mathcal{P}(K)\times\mathcal{P}(K)\left\langle \left\langle z_1,\ldots,z_N\right\rangle\right\rangle} &=& \sum_{w\in\mathcal{F}_N}\left[ f,g_w\right]_{\mathcal{P}(K)}z^{w^T},\\
\left\langle \sum_{w\in\mathcal{F}_N}f_wz^w,c\right\rangle_{\mathcal{C}\left\langle \left\langle z_1,\ldots,z_N\right\rangle\right\rangle\times\mathcal{C}} &=& \sum_{w\in\mathcal{F}_N}\left\langle f_w,c\right\rangle_\mathcal{C}z^w.
\end{eqnarray*}
Then the reproducing kernel property \eqref{repr} can be rewritten as
\begin{equation}\label{repr'}
\left[f,K(\cdot,z)c\right]_{\mathcal{P}(K)\times\mathcal{P}(K)\left\langle \left\langle z_1,\ldots,z_N\right\rangle\right\rangle}=\left\langle f(z),c\right\rangle_{\mathcal{C}\left\langle \left\langle z_1,\ldots,z_N\right\rangle\right\rangle\times\mathcal{C}}.
\end{equation}
The space $\mathcal{P}(K)$ endowed with the metric $[\cdot,\cdot]_{\mathcal{P}(K)}$ will be said to be a \emph{NFRKPS}
associated with the FPS kernel $K(z,z')$. It is clear that this space is isomorphic to the RKPS associated with the kernel $K_{w,w'}$ on $\mathcal{F}_N$, and this isomorphism is well defined by
$$K_{w'}(\cdot)c\mapsto K_{\cdot,w'}c,\quad w'\in\mathcal{F}_N,c\in\mathcal{C}.$$
Let us now come back to the kernels \eqref{kfk} and \eqref{kf} (see also \eqref{kerns}). Clearly, they can be rewritten as
\begin{eqnarray}\label{kfk-phi}
K^{F,k}(z,z') &=& \varphi_k(z)H_k^{-1}\varphi_k(z')^*,\quad k=1,\ldots,N,\\
 K^{F}(z,z') &=& \varphi(z)H^{-1}\varphi(z')^*,\label{kf-phi}
\end{eqnarray}
where rational FPSs $\varphi_k,\ k=1,\ldots,N,$ and $\varphi$ are determined by a given minimal GR-realization $\alpha$ of the FPS $F$ as
\begin{eqnarray*}
\varphi(z) &=& C(I_\gamma-\Delta(z)A)^{-1},\\
\varphi_k(z) &=& \varphi(z)\big|_{\mathbb{C}^{\gamma_k}},\quad  k=1,\ldots,N.
\end{eqnarray*}
For a model minimal GR-realization of $F$, we will start, conversely, with establishing an explicit formula for the kernels \eqref{kfk} and \eqref{kf} in terms of $F$ and then define a minimal GR-realization via these kernels.

Suppose that for a fixed $k\in\{ 1,\ldots,N\}$, \eqref{kfk-phi} holds with some rational FPS $\varphi_k$. Recall that
\begin{equation}\label{start}
J-F(z)JF(z')^*=\sum_{k=1}^N\varphi_k(z)H_k^{-1}(z_k+(z_k')^*)\varphi_k(z')^*
\end{equation}
(note that $(z_k')^*=z_k'$). Then for any $n\in\mathbb{N}$ and $Z,Z'\in\mathbb{C}^{n\times n}$:
\begin{equation}\label{start-z}
J\otimes I_n-F(Z)(J\otimes I_n)F(Z')^*=\sum_{k=1}^N\varphi_k(Z)(H_k^{-1}\otimes(Z_k+(Z_k')^*))\varphi_k(Z')^*.
\end{equation}
Therefore, for $\lambda\in\mathbb{C}$:
\begin{eqnarray}
\nonumber \lefteqn{J\otimes I_{2n}-F(\Lambda_{Z,Z'}(\lambda))(J\otimes I_{2n})
 F({\rm diag}(-Z^*,Z'))^*}\\
&=& \lambda\varphi_k(\Lambda_{Z,Z'}(\lambda))\left\{H_k^{-1}\otimes\begin{pmatrix} I_n & I_n\\
I_n & I_n\end{pmatrix}\right\}\varphi_k({\rm diag}(-Z^*,Z'))^*,\label{id-lambda}
\end{eqnarray}
where
\begin{multline*}
\Lambda_{Z,Z'}(\lambda) := \lambda\begin{pmatrix} I_n & I_n\\
I_n & I_n\end{pmatrix}\otimes e_k+\begin{pmatrix} Z & 0\\
0 & -{Z'}^*\end{pmatrix} \\
 = \left(\begin{pmatrix} Z_1 & 0\\
0 & -{(Z'_1)}^*\end{pmatrix},\ldots,\begin{pmatrix} Z_{k-1} & 0\\
0 & -{(Z'_{k-1})}^*\end{pmatrix},\begin{pmatrix} \lambda I_n+Z_k & \lambda I_n\\
\lambda I_n & \lambda I_n-{(Z'_k)}^*\end{pmatrix},\right. \\
 \left.
\begin{pmatrix} Z_{k+1} & 0\\
0 & -{(Z'_{k+1})}^*\end{pmatrix},\ldots,\begin{pmatrix} Z_N & 0\\
0 & -{(Z'_N)}^*\end{pmatrix}\right),
\end{multline*}
$${\rm diag}(-Z^*,Z'):=\left(\begin{pmatrix} -Z^*_1 & 0\\
0 & Z'_1\end{pmatrix},\ldots,\begin{pmatrix} -Z^*_N & 0\\
0 & Z'_N\end{pmatrix}\right),$$
and, in particular,
$$\Lambda_{Z,Z'}(0)={\rm diag}(Z,-{Z'}^*).$$
For $Z$ and $Z'$ where both $F$ and $\varphi_k$ are holomorphic,
$\varphi_k\left(\Lambda_{Z,Z'}(\lambda)\right)$
is continuous in $\lambda$, and
$F\left(\Lambda_{Z,Z'}(\lambda)\right)$
is holomorphic in $\lambda$ at $\lambda=0$. Thus, dividing by $\lambda$ the expressions in both sides of \eqref{id-lambda} and passing to the limit as $\lambda\to 0$, we get
\begin{eqnarray*}
\lefteqn{-\frac{d}{d\lambda}\left\{F\left(\Lambda_{Z,Z'}(\lambda)\right)\right\}\big|_{\lambda=0}(J\otimes I_{2n})
 F({\rm diag}(-Z^*,Z'))^*}\\
 &=& \varphi_k\left({\rm diag}(Z,-{Z'}^*)\right)\left\{H_k^{-1}\otimes\begin{pmatrix} I_n & I_n\\
I_n & I_n\end{pmatrix}\right\}\varphi_k({\rm diag}(-Z^*,Z'))^*\\
 &=& \begin{pmatrix} \varphi_k(Z)\\
\varphi_k({-Z'}^*)\end{pmatrix}(H_k^{-1}\otimes I_n)\begin{pmatrix} {\varphi_k(-Z^*)}^* & \varphi_k(Z')^*\end{pmatrix}.
\end{eqnarray*}
Taking the $(1,2)$-th entry of the $2\times 2$ block matrices in this equality, we get:
\begin{equation}\label{kfk-expl}
K^{F,k}(Z,Z')=-\frac{d}{d\lambda}\left\{F\left(\Lambda_{Z,Z'}(\lambda)\right)_{12}\right\}\big|_{\lambda=0}(J\otimes I_{n})
 F(Z')^*.
\end{equation}
Using the FPS representation for $F$ we obtain from
\eqref{kfk-expl} the representation
$$
K^{F,k}(Z, Z^{\prime})=\sum_{w,w'\in\mathcal{F}_N}\left(
\sum_{v,v'\in\mathcal{F}_N:\,vv'=w'}(-1)^{|v'|+1}F_{w{g_k}v^{\prime T}}
JF_v\right)\otimes Z^w\left(Z^{\prime *}\right)^{w'^T}.$$
From Corollary~\ref{cor:vanish-fps} we get the expression for a
FPS $K^{F,k}(z,z')$, namely:
\begin{equation}
\label{kfk-fps}
K^{F,k}(z, z^{\prime})=\sum_{w,w'\in\mathcal{F}_N}\left(
\sum_{v,v'\in\mathcal{F}_N:\,vv'=w'}(-1)^{|v'|+1}F_{w{g_k}v^{\prime T}}
JF_v\right) z^w{z'}^{w'^T}.\end{equation}
Using formal differentiation with respect to $\lambda$ we can also represent
this  kernel as
\begin{equation}
\label{kfk-formal}
K^{F,k}(z,z')=-\frac{d}{d\lambda}\left\{F\left(\Lambda_{z,z'}(\lambda)\right)_{12}\right\}\big|_{\lambda=0}J
 F(z')^*.
\end{equation}
We note that one gets \eqref{kfk-fps} and \eqref{kfk-formal}
from \eqref{start} using the same argument applied to FPSs.

Let us now consider the NFRKPSs ${\mathcal K}_k(F)$, $k=1,\ldots, N$, and ${\mathcal K}(F)=
\bigoplus_{k=1}^N{\mathcal K}_k(F)$. They are finite dimensional and isomorphic
to the reproducing kernel Pontryagin spaces on $\mathcal{F}_N$ which were denoted above with the same notation. Thus
\begin{equation}
\label{dim-k}
\begin{split}
{\rm dim}\,{\mathcal K}_k(F)&=\gamma_k,\qquad k=1,\ldots, N,\\
{\rm dim}\,{\mathcal K}(F)&=\gamma.
\end{split}
\end{equation}
The space ${\mathcal K}(F)$ is a multivariable non-commutative analogue
of a certain de Branges--Rovnyak space (see \cite[p.~24]{dBR}, \cite[Section~6.3]{AD1}, and \cite[p.~217]{AG1}).
\subsection{Minimal realizations in
non-commutative de Branges--Rovnyak spaces}
Let us define for every $k\in\left\{1,\ldots, N\right\}$ the backward shift operator
$$R_k:\ {\mathbb C}^q\left\langle\left\langle z_1,\ldots,z_N\right\rangle\right\rangle_{\rm rat}\longrightarrow{\mathbb C}^q\left\langle\left\langle z_1,\ldots,z_N\right\rangle\right\rangle_{\rm rat}$$
by
$$
R_k\,:\ \sum_{w\in\mathcal{F}_N}f_w z^w\longmapsto
\sum_{w\in\mathcal{F}_N}f_{w g_k} z^w.
$$
(Compare with the one-variable backward shift operator $R_0$ considered in Section~\ref{sec:intr}.)
\begin{lem}
\label{lem:shifts}
Let $F$ be a matrix-$J$-unitary on ${\mathcal J}_N$ rational FPS.
Then for every $k\in\left\{1,\ldots, N\right\}$ the following
is true:
\begin{enumerate}
\item $R_kF(z)c\in{\mathcal K}_k(F)$ for every $c\in{\mathbb C}^q$;
\item $R_k{\mathcal K}_j(F)\subset {\mathcal K}_k(F)$ for every $j\in\left\{1,\ldots, N\right\}$.
\end{enumerate}
\end{lem}
\begin{proof} From \eqref{start} and the $J$-unitarity of $F_\emptyset$ we
get
\begin{eqnarray*}
J-F(z)JF_\emptyset^* &=& (F_\emptyset-F(z))JF_\emptyset^*
=-\sum_{k=1}^NR_kF(z)z_kJF_\emptyset^*\\
&=& \sum_{k=1}^N\varphi_k(z)H_k^{-1}z_k\left(\varphi_k\right)_\emptyset^*,
\end{eqnarray*}
and therefore for every $k\in\left\{1,\ldots, N\right\}$
and every $c\in{\mathbb C}^q$ we get
$$
R_kF(z)c=-\varphi_k(z)H_k^{-1}(\varphi_k)_\emptyset^*JF_\emptyset
c=K_\emptyset^{F,k}(z)\left(-JF_\emptyset c\right)\in {\mathcal K}_k(F).
$$
Thus, the first statement of this Lemma is true. To prove the second statement we start again from
\eqref{start} and get for a fixed $j\in\left\{1,\ldots, N\right\}$ and
$w\in\mathcal{F}_N$:
$$-F(z)JF_{w g_j}^*=\varphi_j(z)H_j^{-1}{(\varphi_j)}^*_w+\sum_{k=1}^N\varphi_k(z)H_k^{-1}z_k{(\varphi_k)}_{w g_j}^*,
$$
and therefore for any $c\in{\mathbb C}^q$:
\[
-\sum_{k=1}^N\left(R_kF(z)JF_{w g_j}^*c\right)z_k=\sum_{k=1}^N
\left(R_kK_w^{F,j}(z)c\right)z_k+\sum_{k=1}^N\left(
K_{w g_j}^{F,k}(z)c\right)z_k.
\]
Hence, one has for every $k\in\left\{1,\ldots, N\right\}$:
\begin{equation}
R_kK_w^{F,j}(z)c=-R_kF(z)JF_{w g_j}^*c-K_{w g_j}^{F,k}(z)c,
\label{rkk}
\end{equation}
and from the first statement of this Lemma we obtain that the right-hand side of this equality
belongs to ${\mathcal K}_k(F)$. Thus, the second statement is true, too.
\end{proof}
We now define operators
$A_{kj}:\,{\mathcal K}_j(F)\rightarrow{\mathcal K}_k(F),\
A:\,{\mathcal K}(F)\rightarrow{\mathcal K}(F),\
B:\,{\mathbb C}^q\rightarrow{\mathcal K}(F),\
C:\,{\mathcal K}(F)\rightarrow {\mathbb C}^q,\
D:\,{\mathbb C}^q\rightarrow {\mathbb C}^q
$
by
\begin{eqnarray}
A_{kj}&=&R_k\big|_{{\mathcal K}_j(F)},\quad k,j=1,\ldots, N,
\label{akj}\\
A&=&(A_{kj})_{k,j=1,\ldots, N},\label{a}\\
\label{bb}
B&:& c\longmapsto\begin{pmatrix}R_1F(z)c\\ \vdots \\ R_NF(z)c\end{pmatrix},\\
\label{cc}
C&:&\begin{pmatrix}f_1(z)\\  \vdots\\ f_N(z)\end{pmatrix}\longmapsto \sum_{k=1}^N
(f_k)_\emptyset,\\
D&=&F_\emptyset.
\label{dd}
\end{eqnarray}
These definitions make sense in view of Lemma~\ref{lem:shifts}.
\begin{thm}
\label{thm:realiz} Let $F$ be a matrix-$J$-unitary on ${\mathcal
J}_N$ rational FPS. Then the GR-node $\alpha=(N;A,B,C,D;{\mathcal
K}(F)=\bigoplus_{k=1}^N {\mathcal K}_k(F),{\mathbb C}^q)$,
with operators defined by
\eqref{akj}--\eqref{dd}, is a minimal GR-realization of $F$.
\end{thm}
\begin{proof} We first check that for every $w\in\mathcal{F}_N:\ w\not=
\emptyset$ we have
\begin{equation}
F_w=\left(C\flat A\sharp B\right)^w.
\label{coefs}
\end{equation}
Let $w=g_k$ for some $k\in\{ 1,\ldots,N\}$. Then for $c\in{\mathbb C}^q$:
$$
\left(C\flat A\sharp B\right)^w c=C_kB_kc=\left(R_kF(z)c\right)_\emptyset=
\left(\sum_{w\in\mathcal{F}_N}F_{w g_k}z^w c\right)_\emptyset=
F_{g_k}c.$$
Assume now that $|w|>1$, $w=g_{j_1}\ldots g_{j_{|w|}}$.
Then for $c\in{\mathbb C}^q$:
\[
\begin{split}
\left(C\flat A\sharp B\right)^w c&=C_{j_1}A_{j_1,j_2}\cdots
A_{j_{|w|-1},j_{|w|}}B_{j_{|w|}}c
=\left(R_{j_1}\cdots R_{j_{|w|}}F(z)c\right)_\emptyset\\
&=\left(\sum_{w'\in\mathcal{F}_N} F_{w' g_{j_1}
\cdots g_{j|w|}}z^{w'}c\right)_\emptyset=F_{g_{j_1}\cdots g_{j_{|w|}}}c=F_w c.
\end{split}
\]
Since $F_\emptyset=D$, we obtain that
$$
F(z)=D+C(I-\Delta(z)A)^{-1}\Delta(z)B,$$
that is, $\alpha$ is a GR-realization of $F$. The minimality of
$\alpha$ follows from \eqref{dim-k}.
\end{proof}
Let us now show how the associated structured Hermitian matrix $H={\rm diag}(H_1,\ldots,
H_N)$ arises from this special realization. Let
$$
h={\rm col}_{1\le j\le N}(K_{w_j}^{F,j}(\cdot)c_j) \quad
{\rm and}\quad
h'={\rm col}_{1\le j\le N}(K_{w_j'}^{F,j}(\cdot)c_j').
$$
Using \eqref{rkk}, we obtain
\begin{eqnarray}
\nonumber
\lefteqn{[A_{kj}h_j,h'_k]_{F,k}+[h_j,A_{jk}h'_k]_{F,j} }\\
&=&
[R_kK_{w_j}^{F,j}(\cdot)c_j,K_{w_k'}^{F,k}(\cdot)c_k']_{F,k}+
[K_{w_j}^{F,j}(\cdot)c_j,
R_jK_{w_k'}^{F,k}(\cdot)c_k']_{F,j}\nonumber\\
&=&
\label{kk}
(c_k^{\prime})^*\left(K_{{w_k' g_k},w_j}^{F,j}+
K^{F,k}_{{w_k'},w_{j}g_j}\right)c_j.
\end{eqnarray}
Let $\stackrel{\circ}{\alpha}=(N;\stackrel{\circ}{A},\stackrel{\circ}{B},
\stackrel{\circ}{C},\stackrel{\circ}{D};
{\mathbb C}^\gamma=\bigoplus_{k=1}^N {\mathbb C}^{\gamma_k},{\mathbb C}^q)$
be any minimal GR-realization of $F$, with the associated structured
Hermitian matrix $\stackrel{\circ}{H}={\rm diag}(\stackrel{\circ}{H_1},
\ldots, \stackrel{\circ}{H_N})$. Then the right-hand side of \eqref{kk}
can be rewritten as
\[
\begin{split}
\lefteqn{(c_k^{\prime })^*\left(K_{{w_k' g_k},w_j}^{F,j}+
K^{F,k}_{{w_k'},w_{j}g_j}\right)c_j}\\
&=
{(c_k^{\prime })}^*\left(\left(\stackrel{\circ}{C}\flat\stackrel{\circ}{A}\right)^{
w_k' g_kg_j}\left(\stackrel{\circ}{H_j}\right)^{-1}
\left(\stackrel{\circ}{A^*}\sharp \stackrel{\circ}{C^*}\right)^{g_jw_j^T}\right.\\
&\left.+
\left(\stackrel{\circ}{C}\flat\stackrel{\circ}{A}\right)^{
w_k' g_k}\left(\stackrel{\circ}{H_k}\right)^{-1}
\left(\stackrel{\circ}{A^*}\sharp \stackrel{\circ}{C^*}\right)^{g_kg_jw_j^T}\right)
c_j\\
&={(c_k^{\prime })}^*\left(\stackrel{\circ}{C}\flat\stackrel{\circ}{A}\right)^{
w_k' g_k}\left(\stackrel{\circ}{A_{kj}}\left(\stackrel{\circ}{H_j}\right)^{-1}+
\left(\stackrel{\circ}{H_k}\right)^{-1}\left(\stackrel{\circ}{A_{kj}}\right)^*\right)
\left(\stackrel{\circ}{A^*}\sharp \stackrel{\circ}{C^*}\right)^{g_jw_j^T}c_j\\
&=-{(c_k^{\prime })}^*\left(\stackrel{\circ}{C}\flat\stackrel{\circ}{A}\right)^{
w_k' g_k}
\stackrel{\circ}{B_k}J\left(\stackrel{\circ}{B_j}\right)^*
\left(\stackrel{\circ}{A^*}\sharp \stackrel{\circ}{C^*}\right)^{g_jw_j^T}c_j\\
&=-{(c_k^{\prime })}^*\left(\stackrel{\circ}{C}\flat\stackrel{\circ}{A}\right)^{
w_k' g_k}
\left(\stackrel{\circ}{H_k}\right)^{-1}\left(\stackrel{\circ}{C_k}\right)^*J
\stackrel{\circ}{C_j}\left(\stackrel{\circ}{H_j}\right)^{-1}
\left(\stackrel{\circ}{A^*}\sharp \stackrel{\circ}{C^*}\right)^{g_jw_j^T}c_j\\
&=-{(c_k^{\prime })}^*K^{F,k}_{w_k',\emptyset}JK_{\emptyset,
w_j}^{F,j}c_j\\
&=-{(c_k^{\prime })}^*\left(K^{F,k}_{\emptyset,w_k'}\right)^*
JK_{\emptyset,w_j}^{F,j}c_j\\
&=-(h'_k)^*_\emptyset J(h_j)_\emptyset.
\end{split}
\]
In this chain of equalities we have exploited the relationship between
$\stackrel{\circ}{A},\stackrel{\circ}{B},\stackrel{\circ}{C},
\stackrel{\circ}{D},J$
and $\stackrel{\circ}{H}$ from Theorem~\ref{thm:Lyap} applied to a GR-node $\stackrel{\circ}{\alpha}$.
Thus we have for all $k,j\in\left\{1,\ldots, N\right\}$:
\begin{equation}
\label{form}
[A_{kj}h_j,h'_k]_{F,k}+[h_j,A_{jk}h'_k]_{F,j}=-(h_k')^*C_k^* JC_jh_j.
\end{equation}
Since this equality holds for generating elements of the spaces
${\mathcal K}_k(F)$, $k=1,\ldots, N$) it extends by linearity to arbitrary
elements
$h={\rm col}(h_1,\ldots, h_N)$ and $h'={\rm col}(h'_1,\ldots, h'_N)$
in ${\mathcal K}(F)$. For $k=1,\ldots, N,$ let
$\langle\,\cdot ,\,\cdot\,\rangle_{F,k}$ be any inner product for which
${\mathcal K}_k(F)$ is a Hilbert space. Thus, ${\mathcal K}(F)$ is a
Hilbert space with respect to the inner product
$$\left\langle h,h'\right\rangle_F:=\sum_{k=1}^N\left\langle h_k,h'_k\right\rangle_{F,k}.$$
Then there exist uniquely defined linear operators
$H_k:\,{\mathcal K}_k(F)\rightarrow{\mathcal K}_k(F)$ such that:
$$
[h_k,h'_k]_{F,k}=\langle H_kh_k, h'_k\rangle_{F,k},\quad k=1,\ldots N,
$$
and so with $H:={\rm diag}(H_1,\ldots, H_N):\,
{\mathcal K}(F)\rightarrow{\mathcal K}(F)$ we have:
$$[h,h']_F=\langle Hh,h'\rangle_F.$$
Since the spaces ${\mathcal K}_k(F)$ are non-degenerate (see \cite{AD1}), the operators
$H_k$ are invertible and \eqref{form} can be rewritten as:
$$
(A^*)_{kj}H_j+H_kA_{kj}=-C_k^*JC_j,\quad k,j=1,\ldots N,$$
which is equivalent to \eqref{L}.

Now, for arbitrary $c,c'\in{\mathbb C}^q$ and $w\in\mathcal{F}_N$
we have:
$$
\langle H_kB_kc, K_{w'}^{F,k}(\cdot)c'\rangle_{F,k}=
[R_kF(\cdot)c, K_{w'}^{F,k}(\cdot)c']_{F,k}={c'}^*F_{w' g_k}
c.$$
On the other hand,
\begin{multline*}
-\langle C_k^*JDc, K_{w'}^{F,k}(\cdot)c'\rangle_{F,k}=-
\langle JF_\emptyset c, C_kK_{w'}^{F,k}(\cdot)c'\rangle_{F,k}=-
\langle JF_\emptyset c, K_{\emptyset,w'}^{F,k}c'
\rangle_{{\mathbb C}^q}\\
=-{c'}^*K_{w',\emptyset}^{F,k}JF_\emptyset c=-{c'}^*(\stackrel{\circ}{C}\flat\stackrel{\circ}{A})^{w' g_k}
\left(\stackrel{\circ}{H_k}\right)^{-1}\left(\stackrel{\circ}{C_k}\right)^*J\stackrel{\circ}{D}c\\
=
{c'}^*\left(\stackrel{\circ}{C}\flat\stackrel{\circ}{A}\right)^{w' g_k}
\stackrel{\circ}{B_k}c={c'}^*\left(\stackrel{\circ}{C}\flat\stackrel{\circ}{A}\sharp
\stackrel{\circ}{B}\right)^{w' g_k}c={c'}^*F_{w' g_k}c.
\end{multline*}
Here we have used the relation \eqref{b} for an arbitrary minimal GR-realization
$\stackrel{\circ}{\alpha}=(N;\stackrel{\circ}{A},\stackrel{\circ}{B},
\stackrel{\circ}{C},\stackrel{\circ}{D};
{\mathbb C}^\gamma=\bigoplus_{k=1}^N {\mathbb C}^{\gamma_k},{\mathbb C}^q)$ of
$F$, with the associated structured Hermitian matrix $
\stackrel{\circ}{H}={\rm diag}(
\stackrel{\circ}{H_1},\ldots,
\stackrel{\circ}{H_N})$.
Thus,
$H_kB_k=-C_k^*JD,\ k=1,\ldots, N$, that is,
$B=-H^{-1}C^*JD,$
and \eqref{b} holds for the GR-node $\alpha$. Finally, by Theorem~\ref{thm:Lyap}, we may conclude that $H={\rm diag}(H_1,\ldots, H_N)$ is the
associated structured Hermitian matrix of the special GR-realization $\alpha$.
\subsection{Examples}
In this subsection we give certain examples of matrix-inner rational FPSs
on ${\mathcal J}_2$ with scalar coefficients (i.e., $N=2$, $q=1$, and $J=1$). We also present the
corresponding non-commutative positive kernels $K^{F,1}(z,z^{\prime})$ and
$K^{F,2}(z,z^{\prime })$ computed using formula \eqref{kfk-formal}.
\begin{ex} $F(z)=(z_1+1)^{-1}(z_1-1)(z_2+1)^{-1}(z_2-1)$.
\[
\begin{split}
K^{F,1}(z,z^{\prime })&=2(z_1+1)^{-1}(z_1^{\prime }+1)^{-1},\\
K^{F,2}(z,z^{\prime })&=2(z_1+1)^{-1}(z_1-1)(z_2+1)^{-1}
(z_2^{\prime }+1)^{-1}(z_1^{\prime }-1)(z_1^{\prime }+1)^{-1}.
\end{split}
\]
\end{ex}
\begin{ex} $F(z)=(z_1+z_2+1)^{-1}(z_1+z_2-1).$
\[
K^{F,1}(z,z^{\prime })=K^{F,2}(z,z^{\prime })=2(z_1+z_2+1)^{-1}
(z_1^{\prime}+z_2^{\prime}+1)^{-1}.
\]
\end{ex}
\begin{ex}
\[
\begin{split}
F(z)&=\left(z_1+(z_2+i)^{-1}+1\right)^{-1}\left(z_1+(z_2+i)^{-1}-1\right)\\
&=
\left((z_2+i)(z_1+1)+1\right)^{-1}\left((z_2+i)(z_1-1)+1\right).
\end{split}
\]
\[
\begin{split}
K^{F,1}(z,z^{\prime })&=
2\left((z_2+i)(z_1+1)+1\right)^{-1}(z_2+i)(z_2^{\prime }-i)
\left((z_1^{\prime }+1)(z_2^{\prime }-i)+1\right)^{-1},\\
K^{F,2}(z,z^{\prime })&=
2\left((z_2+i)(z_1+1)+1\right)^{-1}
\left((z_1^{\prime }+1)(z_2^{\prime }-i)+1\right)^{-1}.
\end{split}
\]
\end{ex}

%\bibliographystyle{amsplain}
%\bibliography{/users/faculty/math/dany/bib/all}
%\bibliography{/home/user/bib/all}
%\bibliography{/root/Desktop/dany/bib/all}
%\bibliography{C:/dany/bib/all}
%\bibliography{C:/WIN98/Desktop/dany/bib/all}
%\bibliography{matrix-J-unitary}
% ------------------------------------------------------------------------
% ------------------------------------------------------------------------
\end{document}